\documentclass[11pt]{amsart}
\usepackage{amssymb,amsmath,txfonts,mathrsfs,mathtools,titletoc, tikz, stmaryrd}
\usepackage{graphicx}
\usepackage{float}
\usepackage[alphabetic]{amsrefs}
\usepackage{verbatim}

\newtheorem{theorem}{Theorem}[section]
\newtheorem{prop}[theorem]{Proposition}
\newtheorem{lemma}[theorem]{Lemma}
\newtheorem{remark}[theorem]{Remark}

\newtheorem{definition}[theorem]{Definition}
\newtheorem{cor}[theorem]{Corollary}

\numberwithin{equation}{section}

\def\pf{{\it Proof:}~}

\begin{document}
\title[The construction of splitting map]{The construction of $\epsilon$-splitting map}
\author{Guoyi Xu, Jie Zhou}
\address{Guoyi Xu\\ Department of Mathematical Sciences\\Tsinghua University, Beijing\\P. R. China}
\email{guoyixu@tsinghua.edu.cn}
\date{\today}
\address{ Jie Zhou\\ Department of Mathematical Sciences\\Tsinghua University;School of Mathematical Sciences\\ Capital Normal University, Beijing\\ P.R. China }
\email{zhoujiemath@cnu.edu.cn}

\begin{abstract}
For a geodesic ball with non-negative Ricci curvature and almost maximal volume, without using compactness argument, we construct an $\epsilon$-splitting map on a concentric geodesic ball with uniformly small radius. There are two new technical points in our proof. The first one is the way of finding $n$ directional points by induction and stratified almost Gou-Gu Theorem. The other one is the error estimates of projections, which guarantee the $n$ directional points we find really determine $n$ different directions.
\\[3mm]
\end{abstract}
\thanks{The first author was partially supported by Research was partially supported by Beijing Natural Science Foundation Z190003, NSFC 11771230 and NSFC 12141103}

\maketitle

\tableofcontents

\section{Introduction}

For a compact $n$-dimensional Riemannian manifold with $Rc\geq (n- 1)$, if its volume is close to the volume of unit round sphere $\mathbb{S}^{n}$, Colding \cite{Colding-shape} proved that the manifold is Gromov-Hausdorff close to $\mathbb{S}^{n}$. Analogue to the positive Ricci curvature case, for a geodesic ball with $Rc\geq 0$ and almost maximal volume, Colding \cite[Theorem $0.8$]{Colding-volume} pointed out that such geodesic ball is Gromov-Hausdorff close to the Euclidean ball, which can be proved in similar way as \cite{Colding-shape}.

Later, Cheeger \cite[Theorem $9.69$]{Cheeger-note} gave a complete proof of \cite[Theorem $0.8$]{Colding-volume}. The proof of Cheeger is by the compactness argument, which relies on the almost cone property and the almost splitting property established in the Ricci limit space theory (see\cite{CC-Ann}, \cite{CC1}, \cite{CC2} and \cite{CC3}).

We recall the definition of $\epsilon$-Gromov-Hausdorff approximation.  Let $(\mathbf{X}, d_\mathbf{X})$ and $(\mathbf{Y}, d_\mathbf{Y})$ be two metric spaces. For $\epsilon> 0$, a map $f: \mathbf{X}\rightarrow \mathbf{Y}$ is called an \textbf{$\epsilon$-Gromov-Hausdorff approximation} if
\begin{equation}\nonumber
\left\{
\begin{array}{rl}
&\mathbf{Y}\subset \mathbf{U}_{\epsilon}\big(f(\mathbf{X})\big) \\
&\Big|d_\mathbf{Y}\big(f(x_1), f(x_2)\big)- d_\mathbf{X}(x_1, x_2)\Big|< \epsilon \ , \quad \quad \quad \quad \forall x_1, x_2\in \mathbf{X}
\end{array} \right.
\end{equation}
where $\mathbf{U}_{\epsilon}\big(f(\mathbf{X})\big)= \Big\{z\in \mathbf{Y}: d\big(z, f(\mathbf{X})\big)\leq \epsilon\Big\}$.  For simplicity reason, we also use G-H approximation instead of Gromov-Hausdorff approximation in the rest of the paper.

In $Rc\geq (n- 1)$ case \cite{Colding-shape}, the corresponding G-H approximation is explicitly constructed, although some topological results are used. On the other hand, the proof of Cheeger for $Rc\geq 0$ case did not provide the way of constructing the corresponding G-H approximation. We are interested in constructing the explicit G-H approximation for $Rc\geq 0$ case. By looking at the model space, it seems the exponential map, the linear harmonic map or some optimal transportation map are possible candidates. But it is not easy to judge which one is canonical from the compactness argument, since we do not know whether the estimates related to these maps are stable under the Gromov-Hausdorff topology. We note that the stable version of splitting property of Ricci limit space is the quantitative version of almost splitting property in \cite{CC-Ann} (see Lemma \ref{lem property 2 of G-H appr}), which is well-known to the experts in this field. This result is our starting point.

Now we recall the definition of $\epsilon$-splitting map introduced in \cite{CN} as follows.  If $Rc(M^n)\geq 0$ and $B_r(p)\subseteq M^n$ is a geodesic ball, then the harmonic map $\mathbf{b}= \big\{\mathbf{b}_i\big\}_{i= 1}^{n}: B_{r}(p)\rightarrow \mathbb{R}^n$ is called an \textbf{$\epsilon$-splitting map}, if $\mathbf{b}$ satisfies
\begin{align}
\sup_{B_{r}(p)\atop i= 1, \cdots, n} |\nabla \mathbf{b}_i|\leq 1+ \epsilon ,\quad \quad and \quad \quad \fint_{B_{r}(p)} \big|\langle \nabla \mathbf{b}_i, \nabla \mathbf{b}_j\rangle- \delta_{ij}\big|^2\leq \epsilon  .\nonumber
\end{align}

From \cite{CC-Ann} and \cite{CC1} (also see \cite{CJN}), the existence of $\epsilon$-G-H approximation from such geodesic ball to the corresponding Euclidean ball, is equivalent to the existence of an $\epsilon'$-splitting map between them.  In fact,  from \cite{Ding} and \cite{Cheeger},  an $\epsilon'$-splitting map comes from pulling back the harmonic functions on the Euclidean space by one $\epsilon$-G-H approximation.  We can view $\epsilon'$-splitting map as a `canonical' G-H approximation, although there is no uniqueness restriction. Hence we can reduce the construction of $\epsilon$-G-H approximation for $Rc\geq 0$ case to the construction of $\epsilon'$-splitting map.

In this paper we construct an $\epsilon$-splitting map on a concentric geodesic ball with uniformly small radius, and establish the corresponding quantitative estimate in term of the volume ratio between geodesic ball and Euclidean ball. Although the quantitative version of almost splitting property in \cite{CC-Ann} (see Lemma \ref{lem property 2 of G-H appr}) is our starting point, it is not enough for our purpose. One new ingredient of our construction is finding $n$ direction points by combining the stratified Gou-Gu Theorem with the induction method, which is the content of Theorem \ref{thm dist of points in geodesic balls induc-dim}.  The $(k+ 1)$th direction point $q_{k+ 1}$ is chosen from the image of $\mathscr{P}_k$, where $\mathscr{P}_k$ is the composition of the first $k$ projection maps with respect to the first $k$ direction points in order.  To establish the stratified Gou-Gu Theorem for $(k+ 1)$ direction points,  the distance between $q_{k+ 1}$ and the origin $p$ needs to be in a larger scale than the scale of the error estimate in the stratified Gou-Gu Theorem. In other words, we need to get a lower bound of the radius of the image of $\mathscr{P}_k$ with respect to the origin $p$.  This lower bound is proved by the stratified Gou-Gu Theorem for $k$ direction points and the volume comparison theorem (see Step (1) of the proof of Theorem \ref{thm dist of points in geodesic balls induc-dim}).

To prove that $n$ direction points we find really determine $n$ almost orthogonal directions,  we need to show that upper bound of discrete Lipschitz constant of distance map is close to $1$ and use the integral Toponogov Theorem.  The first one is from distance aspect, the second one is from angle aspect.

One technical difficulty we overcome is, proving the corresponding error estimates of projections, which yield that the upper bound of discrete Lipschitz constant of distance map is close to $1$. These estimates of projections are trivial in Euclidean case,  to obtain them on manifolds, we need to use the stratified Gou-Gu Theorem carefully,  and the induction method will be used again.  These estimates are provided in Section \ref{sec quasi-isom}. After these preparation work, we prove the distance map determined by the $n$ direction points is in fact a quasi-isometry.

The reason that we can only deal with the geodesic ball of uniformly small radius, is that we do not have the quantitative version of the almost cone property established in \cite{CC-Ann}. Lack of the almost cone property, we can not relate the cases with variant scales in the expected quantitative way.

In Section \ref{sec existence of splitting map}, we will use the fact that the upper bound of discrete Lipschitz constant is close to $1$, combining with the integral Toponogov Theorem for $Rc\geq 0$ to establish the almost orthogonality of distance maps.

The integral Toponogov Theorem for $Rc\geq 0$ was proved in \cite{Colding-volume} by approximating the distance function by harmonic functions and covering argument. Motivated by the cosine law for triangle (the terms with distance functions are in the form of the square of distance functions), we consider the model function $f$ defined in Subsection \ref{subsec integral est of diff} relating with the square of the distance function, present a different proof of this result here (Corollary \ref{cor Colding-1.28}), which avoids the covering argument. And our argument is consistent with the former argument for almost Gou-Gu Theorem of distance functions. Finally we show that the harmonic map constructed from the distance map is an $\epsilon$-splitting map.

\section{Almost Gou-Gu Theorem of distance functions}\label{sec general Pytha}

We always assume $\beta\geq 3$ in this section unless otherwise mentioned.

\subsection{Distance estimate by multiple line integral of difference}

For $B_r(p)\subseteq (M^n, g)$ and $\beta\geq 3$, let $d(p, q)= \beta r$. For $x\in  B_r(p)$, let  $\rho(x)= d(q, x)$.  If there exists $\theta_x\in S(T_qM)$ such that $x=\exp_q(\rho(x)\theta_x)$ and $\exp_q(t\theta_x)$ is a segment up to $t=\beta r$, then define  $\pi(x)=\exp_q(\beta r \theta_x)\in \partial B_{\beta r}(q)$.

In the following, we firstly estimate the distance $d(x,y)$ of $x,y\in B_r(p)$ when $\pi(x), \pi(y)$ are well defined.  In this case, we define $\displaystyle \sigma_x(t)=\gamma_{\pi(x),x}(t)$,  for  $t\in [0,d(x,\pi(x))]$, where $\gamma_{\pi(x),x}$ is a geodesic segment from $\pi(x)$ to $x$; and the relative velocity $r_x^y$ (for $\rho(x)\neq \beta r$) by $\displaystyle r_x^y:=\frac{d(\pi(y),y)}{d(\pi(x),x)}= \Big|\frac{\beta r-\rho(y)}{\beta r-\rho(x)}\Big|$. For $t\in [0,d(x,\pi(x))]$, denote
\begin{align*}
\tilde{\sigma}_y(t)=\sigma_y(r_x^yt),\ \ \  \tau_t(s)= \gamma_{\sigma_x(t), \tilde{\sigma}_y(t)}(s) \text{ for } s\in [0, l_t],
\end{align*}
where $\displaystyle l_t:= d\big(\sigma_x(t), \tilde{\sigma}_y(t)\big)$. Denote $\mathfrak{S}(x)\vcentcolon = \mathcal{S}(\rho(x))$, where
\begin{align*}
\mathcal{S}(t)=\left\{
\begin{aligned}
 1&, & \text{ if }t-\beta r>0, \\
-1&, & \text{ if } t-\beta r<0.
\end{aligned}
\right.
\end{align*}

One general philosophy of distance estimate and derivative estimate of distance function on manifolds is, to reduce the estimates to the $C^0$ estimate, integral gradient estimate of error functions between original function and make-up function $f$ and the Hessian estimate of $f$.  Although the make-up function $f$ is freely chosen, to get the Hessian estimate of $f$, we usually need to choose $f$ as the solution of suitable elliptic PDE. One typical technical tool is the following estimate of distance functions.

\begin{lemma}\label{lem property 2 of G-H appr}
{For $x, y\in B_r(p)$ with well-defined $\pi(x), \pi(y)$,  we have
\begin{align}
&\quad \quad \Big|d(x, y)^2- \Big(\big[\rho(x)- \rho(y)\big]^2+ d(\pi(x),  \pi(y))^2\Big)\Big| \nonumber\\
&\leq C(n)\Big\{\int_0^{d(x,\pi(x))} \big|\nabla (\rho^2- f)\big|(\tilde{\sigma}_y(t))+ \big|\nabla (\rho^2- f)\big|(\sigma_x(t))dt \nonumber \\
&\quad  + \sup_{B_{2r}(p)} |\rho^2- f|+ \int_0^{d(x,\pi(x))} dt\int_{\gamma_{\sigma_x(t), \tilde{\sigma}_y(t)}(s)} |\nabla^2 f- 2g| ds+ r^2\beta^{-1} \Big\}.\nonumber
\end{align}
}
\end{lemma}

\begin{remark}\label{rem almost G-G in Rn}
{Note if $M=\mathbb{R}^n$, then it is straightforward to verify that
$$\sup_{x,y\in B_r(p)}|d(x,y)^2-d(\pi(x),\pi(y))^2-|\rho(x)-\rho(y)|^2|\le \frac{12r^2}{\beta}.$$
The above lemma is the Riemannian manifolds version of the almost Gou-Gu inequality in $\mathbb{R}^n$.
}
\end{remark}

\begin{figure}[H]
\begin{center}
\includegraphics{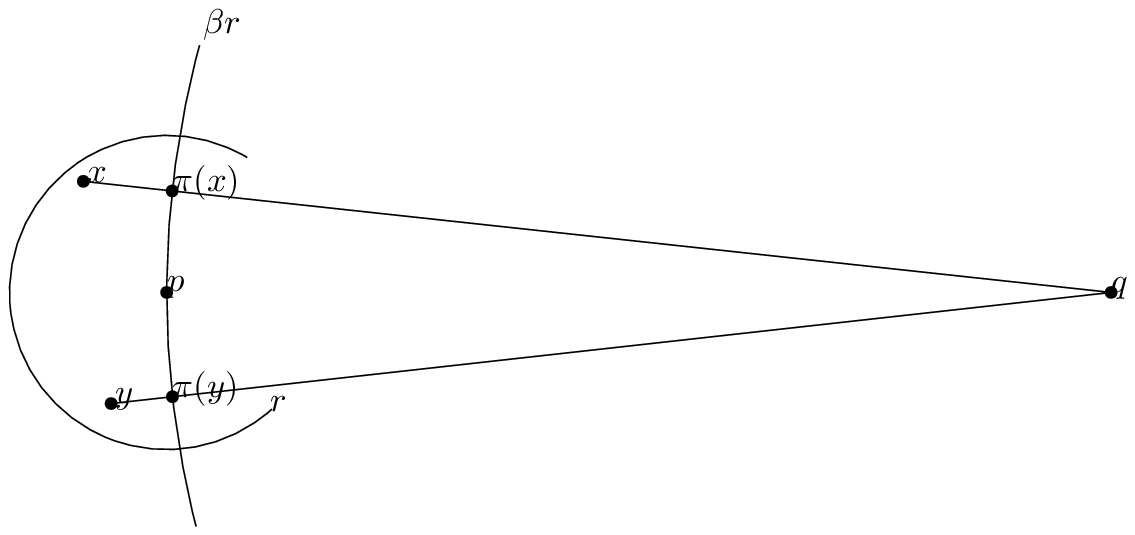}
\caption{Lemma \ref{lem property 2 of G-H appr}}
\label{figure: lemma2.1}
\end{center}
\end{figure}

\pf
{\textbf{Step (1)}. Since
\begin{align*}
\int_0^{d(x,\pi(x))} \big|\nabla (\rho^2- f)\big|(\tilde{\sigma}_y(t))+ \big|\nabla (\rho^2- f)\big|(\sigma_x(t))dt =\int_{\gamma_{\pi(y),y}\cup \gamma_{\pi(x),x}} \big|\nabla (\rho^2- f)\big|,
\end{align*}
 we know the right hand term is symmetric with respect to $x$ and $y$. Thus we can assume $d(x,\pi(x))\ge d(y,\pi(y))$ without loss of generality.

For $\zeta\in (0, 1)$ to be determined later, we firstly assume $\displaystyle d(x,\pi(x))=|\rho(x)- \beta r| \geq \zeta r$. Define $\displaystyle \alpha= \frac{\rho(y)- \rho(x)}{d(x,\pi(x))}$. Then $|\alpha|\le \frac{2}{\zeta}$ and
 \begin{align}
\rho\big(\tau_t(l_t)\big)=\beta r+\big(\alpha+\mathfrak{S}(x)\big)t \ ,\quad \quad
\rho\big(\tau_t(0)\big)=\beta r+\mathfrak{S}(x)t. \nonumber
\end{align}
Let $\mathcal{U}_t$ be the solution of
\begin{equation}\nonumber
\left\{
\begin{array}{rl}
\mathcal{U}_t''&= 1 \\
\mathcal{U}_t(0)&= \frac{1}{2}\big(\beta r+\mathfrak{S}(x)t\big)^2, \quad \quad and \quad \quad \mathcal{U}_t(l_t)= \frac{1}{2}\bigg(\beta r+\big(\alpha+\mathfrak{S}(x)\big)t\bigg)^2 . \\
\end{array} \right.
\end{equation}
Then for any $s_1, s_2\in [0, l_t]$, we have
\begin{align}
\big|(\frac{1}{2} f(\tau_t)- \mathcal{U}_t)'(s_2)- (\frac{1}{2} f(\tau_t)- \mathcal{U}_t)'(s_1)\big|\leq \int_0^{l_t} |\nabla^2\frac{1}{2} f- g|(\tau_t(s))ds \label{1st deri diff}
\end{align}
and
\begin{align}
\mathcal{U}_t'(l_t)- \mathcal{U}_t'(0)= l_t \label{diff of 1st derivative}
\end{align}
From the mean value theorem, there is some $\xi\in [0, l_t]$ such that
\begin{align}
(\frac{1}{2} f(\tau_t)- \mathcal{U}_t)(l_t)- (\frac{1}{2} f(\tau_t)- \mathcal{U}_t)(0)= l_t\cdot (\frac{1}{2} f(\tau_t)- \mathcal{U}_t)'(\xi), \nonumber
\end{align}
which implies
\begin{align}
\big|(\frac{1}{2} f(\tau_t)- \mathcal{U}_t)'(\xi)\big|\leq \frac{1}{l_t}\sup_{\sigma_x[0, t]\cup \tilde{\sigma}_y[0, t]} |\rho^2- f| . \label{one point deri bound}
\end{align}
From (\ref{1st deri diff}) and (\ref{one point deri bound}), we get
\begin{align}
\sup_{s\in [0, l_t]}\big|(\frac{1}{2} f(\tau_t)- \mathcal{U}_t)'(t)\big|\leq \frac{\sup\limits_{\sigma_x[0, t]\cup \tilde{\sigma}_y[0, t]} |\rho^2- f|}{l_t}+ \int_0^{l_t} |\nabla^2\frac{1}{2} f- g|(\tau_t(s))ds . \label{crucial tilde U est}
\end{align}

\textbf{Step (2)}. Since $l_t$ is almost everywhere differentiable on $[0, d(x, \pi(x))]$,  from the extension of the first variation formula (for example, see \cite{Liu}) and (\ref{diff of 1st derivative}), we have
\begin{align}
\frac{d}{dt}l_t&= (\alpha+ \mathfrak{S}(x))\langle \nabla \rho, \tau_t'\rangle(l_t)- \mathfrak{S}(x)\langle \nabla \rho, \tau_t'\rangle(0) \nonumber \\
 &= \frac{(\alpha+ \mathfrak{S}(x))t+ \mathfrak{S}(x)\beta r+ \frac{1}{2}\alpha \beta r}{\big[\beta r+ (\alpha+ \mathfrak{S}(x))t\big](\beta r+\mathcal{ S}(x)t)}l_t+ \frac{\alpha^2\beta r t\big[\beta r+ \frac{1}{2}(\alpha+ 2\mathfrak{S}(x))t\big]}{\big[\beta r+ (\alpha+ \mathfrak{S}(x))t\big](\beta r+ \mathfrak{S}(x)t)}\frac{1}{l_t}
 \nonumber \\
 &\quad  + (I)+ (II) , \label{eq need 1.1.1}
\end{align}
where
\begin{align}
(I)& = \frac{\alpha+ \mathfrak{S}(x)}{\beta r+ (\alpha+ \mathfrak{S}(x))t}\Big[(\frac{1}{2} \rho^2\circ \tau_t)'(l_t)- (\frac{1}{2} f\circ \tau_t)'(l_t)\Big] \nonumber \\
&\quad - \frac{\mathfrak{S}(x)}{\beta r+ \mathfrak{S}(x)t}\Big[(\frac{1}{2} \rho^2\circ \tau_t)'(0)- (\frac{1}{2} f\circ \tau_t)'(0)\Big] \nonumber \\
(II)& = \frac{\alpha+ \mathfrak{S}(x)}{\beta r+ (\alpha+ \mathfrak{S}(x))t}\big[(\frac{1}{2} f\circ \tau_t)'- \mathcal{U}_t'\big](l_t)- \frac{\mathfrak{S}(x)}{\beta r+ \mathfrak{S}(x)t}\big[(\frac{1}{2} f\circ \tau_t)'- \mathcal{U}_t'\big](0) .\nonumber
\end{align}

Simplifying (\ref{eq need 1.1.1}), let $h(t)= \frac{\beta r+ \mathfrak{S}(x)t}{\beta r+ (\alpha+\mathfrak{S}(x))t}+ \frac{\beta r+ (\alpha+ \mathfrak{S}(x))t}{\beta r+\mathfrak{S}(x)t}$,
we have
\begin{align}
&\frac{\big[\beta r+ (\alpha+ \mathfrak{S}(x))t\big](\beta r+\mathfrak{S}(x)t)}{2l_t}\cdot \frac{d}{dt}\Big\{\frac{l_t^2}{\big[\beta r+ (\alpha+ \mathfrak{S}(x))t\big](\beta r+ \mathfrak{S}(x)t)}-h(t)\Big\} \nonumber \\
&= (I)+ (II) \label{first diff equa 1.1}
\end{align}

By $d(\pi(x),x)=|\rho(x)-\beta r|=|\rho(x)-\rho(p)|\le d(x,p)\le r$, we know $\pi(x)\in B_{2r}(p)$ and hence $\sigma_x(t)\in B_{2r}(p)$ for any $t\in [0,d(x,\pi(x))]$. Similarly, $d(\tilde{\sigma}_y(t),p)\le 2r$. From (\ref{first diff equa 1.1}) and (\ref{crucial tilde U est}), note $\beta \ge 3$,  $|(\alpha+\mathfrak{S}(x))t|, |\mathfrak{S}(x)t|\le 2r$, and $l_t\le 4r$.
\begin{align}
&\quad \Big|\frac{d}{dt}\Big\{\frac{l_t^2}{\big[\beta r+ (\alpha+ \mathfrak{S}(x))t\big](\beta r+\mathfrak{S}(x) t)}-h(t)\Big\}\Big| \nonumber \\
&\leq \frac{216(|\alpha|+ 1)r}{(\beta r)^3}\Big\{\big|\nabla (\frac{1}{2} \rho^2- \frac{1}{2} f)\big|(\tau_t(l_t))+ \big|\nabla (\frac{1}{2} \rho^2- \frac{1}{2} f)\big|(\tau_t(0))\Big\}\nonumber \\
&\quad + \frac{216(|\alpha|+ 1)}{(\beta r)^3}\Big\{2\sup_{\sigma_x[0, t]\cup \tilde{\sigma}_y[0, t]} |\rho^2- f|+ 4r\int_0^{l_t} |\nabla^2\frac{1}{2} f- g|(\tau_t(s))ds\Big\}
\end{align}

\textbf{Step (3)}. Since for any $t_1, t_2\in [0,d(x,\pi(x))]$,
 $$|l_{t_1}-l_{t_2}|\le d(\sigma_x(t_1),\sigma_x(t_2))+d(\tilde{\sigma}_{y}(t_1),\tilde{\sigma}_{y}(t_2))\le (1+r_x^y)|t_2-t_1|,$$ we know $l_t$ is a Lipschitz function. Then we get
\begin{align}
&\quad \Big|\Big\{\frac{l_t^2}{\big[\beta r+ (\alpha+ \mathfrak{S}(x))t\big](\beta r+ \mathfrak{S}(x)t)}-h(t)\Big\}- \Big\{\frac{l_0^2}{(\beta r)^2}-h(0)\Big\}\Big| \nonumber \\
&\leq \frac{C(n)}{(\beta r)^3\zeta}\Big\{r\int_0^t \big|\nabla (\rho^2- f)\big|(\tau_v(l_v))+ \big|\nabla (\rho^2- f)\big|(\tau_v(0))dv \nonumber \\
&\quad \quad \quad + t\sup_{\sigma_x[0, t]\cup \tilde{\sigma}_y[0, t]} |\rho^2- f|+ r\int_0^t \int_0^{l_v} |\nabla^2 f- 2g|(\tau_v(s))dsdv\Big\}.\nonumber
\end{align}

Let $t= d(x,\pi(x))\leq r$ in the above, note $l_{d(x,\pi(x))}= d(x, y)$ and $\displaystyle \Big(\sigma_x[0, d(x,\pi(x))]\cup \tilde{\sigma}_y[0, d(x,\pi(x))]\Big)\subseteq B_{2r}(p)$, we have
\begin{align}
&\quad \Big|\Big\{\frac{d(x, y)^2}{\rho(y)\rho(x)}-\frac{\rho(x)}{\rho(y)}-\frac{\rho(y)}{\rho(x)}\Big\}- \Big\{\frac{l_0^2}{(\beta r)^2}-2\Big\}\Big|\nonumber  \\
&\leq \frac{C(n)}{\beta^3 r^2\zeta}\Big\{\int_0^{d(x,\pi(x))} \big|\nabla (\rho^2- f)\big|(\tilde{\sigma}_y(t))+ \big|\nabla (\rho^2- f)\big|(\sigma_x(t))dt \nonumber \\
&\quad  + \sup_{B_{2r}(p)} |\rho^2- f|+ \int_0^{d(x,\pi(x))} \int_{\gamma_{\sigma_x(t), \tilde{\sigma}_y(t)}} |\nabla^2 f- 2g| dt\Big\} . \nonumber
\end{align}

Note $\displaystyle |(\frac{\rho(x)\rho(y)}{(\beta r)^2})- 1|\ell_0^2\leq C(n)\frac{r^2}{\beta}$, we obtain
\begin{align}
&\quad \quad \Big|d(x, y)^2- \Big(\big[\rho(x)- \rho(y)\big]^2+ l_0^2\Big)\Big| \nonumber\\
&\leq \frac{C(n)}{\beta\zeta}\Big\{\int_0^{d(x,\pi(x))} \big|\nabla (\rho^2- f)\big|(\tilde{\sigma}_y(s))+ \big|\nabla (\rho^2- f)\big|(\sigma_x(s))ds \nonumber \\
&\quad  + \sup_{B_{2r}(p)} |\rho^2- f|+ \int_0^{d(x,\pi(x))} \int_{\gamma_{\sigma_x(s), \tilde{\sigma}_y(s)}} |\nabla^2 f- 2g| ds\Big\}\nonumber \\
&\quad + C(n)r^2(\zeta + \beta^{-1}) .\label{needed in lemma}
\end{align}

If $|\rho(x)-\beta r|< \zeta r$, then by assumption,  $|\rho(y)- \beta r|=d(y,\pi(y))\le d(x,\pi(x))< \zeta r$. Thus $|\ell_0- d(x, y)|\leq 2\zeta r$ and we get
\begin{align}
\Big|d(x, y)^2- \Big(\big[\rho(x)- \rho(y)\big]^2+ l_0^2\Big)\Big|\leq C(n)\zeta r^2. \nonumber
\end{align}

Combining this with (\ref{needed in lemma}) and taking $\zeta= \beta^{-1}$, the conclusion follows.
}
\qed

\subsection{Integral estimates of model functions on balls}\label{subsec integral est of diff}

In this subsection, we define the volume of the geodesic sphere $V(s)= V(\rho^{-1}(s))$, where $\rho(x)= d(x, q)$.
\begin{definition}\label{def F and G}
{We define $f$ satisfying
\begin{equation}\nonumber
\left\{
\begin{array}{rl}
&\Delta  f= 2n, \quad \quad \quad \quad on\ B_r(q), \\
&f\big|_{\partial B_r(q)}= r^2 . \nonumber
\end{array} \right.
\end{equation}
}
\end{definition}
\begin{remark}This definition is similar to the one in Cheeger-Colding's almost metric cone theorem. And the following Propositions \ref{prop L2 diff of gadient of F}-\ref{prop integral Hessian small} are, in
some sense, an adaption of the proof of almost metric cone theorem in the case of almost maximal volume growth. The equation above
 is the one satisfied by  the function occurs in the cosine law of the model space, and hence contains much geometric information.
\end{remark}
\begin{lemma}\label{lem C0 bound of G}
{We have
\begin{align}
&-C(n)r^2\leq f(x)\leq \rho^2(x) \ , \quad \quad \quad \quad \forall x\in B_r(q) .\nonumber
\end{align}
}
\end{lemma}

\pf
{By Laplacian Comparison Theorem, $\Delta(f- \rho^2)\geq 0$. By the Maximum Principle, we have $f(x)\leq \rho^2(x)$. Let $\phi(x)= d(z, x)$, where $d(z, q)= 2r$. Then from Laplacian Comparison Theorem, $\Delta(\phi^{4- 2n})\geq 2(n- 2)^2\phi^{2- 2n}\geq 2(n- 2)^2 3^{2- 2n}r^{2- 2n}$ for any $x\in B_r(q)$. Now we get
\begin{align}
\Delta\Big(\frac{n \phi^{4- 2n}}{(n- 2)^23^{2- 2n}r^{2- 2n}}- f\Big)\geq 0. \nonumber
\end{align}

Apply the Maximum principle again, we obtain
\begin{align}
\sup_{B_r(q)} \Big(\frac{n \phi^{4- 2n}}{(n- 2)^23^{2- 2n}r^{2- 2n}}- f\Big)\leq \max_{\partial B_r(q)} \Big(\frac{n \phi^{4- 2n}}{(n- 2)^23^{2- 2n}r^{2- 2n}}- f\Big)\leq \frac{n r^{4- 2n}}{(n- 2)^23^{2- 2n}r^{2- 2n}}- r^2, \nonumber
\end{align}
which implies the conclusion by $\displaystyle \sup_{B_r(q)}\phi\leq 3r$.
}
\qed

\begin{prop}\label{prop L2 diff of gadient of F}
{If $\displaystyle Rc\geq 0$ and $\displaystyle \frac{V(\partial B_r(q))}{r^{n- 1}}\geq (1- \omega)n\omega_n$, where $\omega\in (0, 1)$, then
\begin{align}
\fint_{B_r(q)} |\nabla (f- \rho^2)|^2 \leq C(n)\omega r^2. \nonumber
\end{align}
}
\end{prop}

\pf
{From Lemma \ref{lem C0 bound of G}, we get
\begin{align}
\int_{B_r(q)} |\nabla (f- \rho^2)|^2&= \int_{B_r(q)}(-\Delta f+ \Delta \rho^2)(f- \rho^2)\leq \sup_{B_r(q)}|\rho^2- f|\int_{B_r(q)}(2n-  \Delta \rho^2) \nonumber \\
&\leq C(n)r^2\Big(2n\mathrm{V}(B_r(q))- \int_{\partial B_r(q)} 2\rho\Big)= C(n)r^2\mathrm{V}(B_r(q))\Big(n- r\frac{\mathrm{V}(\partial B_r(q))}{\mathrm{V}(B_r(q))}\Big) \nonumber \\
&\leq C(n)\omega r^2. \nonumber
\end{align}
which implies the conclusion.
}
\qed

\begin{prop}\label{prop C0 diff of F}
{Assume $\displaystyle Rc\geq 0$ and $\displaystyle \frac{V(\partial B_r(q))}{n\omega_nr^{n- 1}}\geq (1- \omega)$, where $\omega\in (0, 1)$, then
\begin{align}
\sup_{B_{(1- 2\omega^{\frac{1}{2n+ 4}})r}(q)}\Big|f- \rho^2\Big|\leq C(n)r^2\omega^{\frac{1}{2n+ 4}} . \nonumber
\end{align}
}
\end{prop}

\pf
{By the scaling invariant property of $Rc\geq 0$, we only need to show the conclusion for $r= 1$. From the Poincare inequality and Proposition \ref{prop L2 diff of gadient of F}, there is
\begin{align}
\fint_{B_1(q)} \big|f- \rho^2\big|^2\leq C(n)\cdot \fint_{B_1(q)} \big|\nabla(f- \rho^2)\big|^2 \leq C(n)\omega .\label{L2 diff of F}
\end{align}

For any $B_{2\ell}(x)\subset B_{2\sqrt{\ell}}(x)\subset B_1(q)$, where $\ell\in (0, 1)$ is to be determined later, then from (\ref{L2 diff of F}),
\begin{align}
\frac{V\big(B_{\ell}(x)\big)}{V(B_1(q))}\min_{B_{\ell}(x)} \big|f- \rho^2\big|^2\leq C(n)\omega . \label{4.50.1}
\end{align}

From Bishop-Gromov Comparison Theorem, we have
\begin{align}
V\big(B_{\ell}(x)\big)\geq \big(\frac{\ell}{2}\big)^n\cdot V\big(B_{2}(x)\big)\geq \big(\frac{\ell}{2}\big)^n\cdot V\big(B_1(q)\big) \label{4.50.2}
\end{align}

By (\ref{4.50.1}) and (\ref{4.50.2}),
\begin{align}
\min_{B_{\ell}(x)} \big|f- \rho^2\big|&\leq C(n)\big(\frac{1}{\ell}\big)^{\frac{n}{2}}\omega^{\frac{1}{2}}. \label{4.50.3}
\end{align}

From Cheng-Yau's gradient estimate and Lemma \ref{lem C0 bound of G}, we get
\begin{align}
\sup_{B_{\ell}(x)}|\nabla f|\leq \sup_{B_{\sqrt{\ell}}(x)}|\nabla f|\leq C(n)(\sqrt{\ell})^{-1}(\sup_{B_{2\sqrt{\ell}}(x)} |f|+ 2n) \leq \frac{C(n)}{\sqrt{\ell}}. \label{gradient of b has a bound}
\end{align}

Now from (\ref{4.50.3}) and (\ref{gradient of b has a bound}),
\begin{align}
&\quad \max_{B_{\ell}(x)} \Big|f- \rho^2\Big|\leq \min_{B_{\ell}(x)} \big|f- \rho^2\big|+ 2\ell \sup_{B_{\ell}(x)} \big|\nabla (f- \rho^2)\big| \nonumber \\
&\leq C(n)\big(\frac{1}{\ell}\big)^{\frac{n}{2}}\omega^{\frac{1}{2}}+ 2\ell\cdot \big(2+ \sup_{B_{\ell}(x)} \big|\nabla f\big|\big) \leq C(n)\big(\frac{1}{\ell}\big)^{\frac{n}{2}}\omega^{\frac{1}{2}}+ C(n)\sqrt{\ell} .\nonumber
\end{align}
Let $\ell= \omega^{\frac{1}{n+ 2}}$ in the above, from the free choice of $x\in B_1(q)$, we have
\begin{align}
\sup_{B_{1- 2\omega^{\frac{1}{2n+ 4}}}(q)}\Big|f- \rho^2\Big|\leq C(n)\omega^{\frac{1}{2n+ 4}} . \nonumber
\end{align}
}
\qed

Now we get the integral Hessian estimate of difference.
\begin{prop}\label{prop integral Hessian small}
{If $\displaystyle Rc\geq 0$ and $\displaystyle \frac{V(\partial B_r(q))}{r^{n- 1}}\geq (1- \omega)n\omega_n$, where $\omega\in (0, \frac{1}{2})$, then
\begin{align*}
\fint_{B_{(1- \omega^{\frac{1}{32}})r}(q)} \big|\nabla^2 f- 2g\big|^2 \leq C(n)\omega^{\frac{1}{4}}.
\end{align*}
}
\end{prop}

\pf
{From \cite[Lemma $2.3$]{Xu-group}, we can choose $\phi\in C^\infty(M^n; [0, 1])$ with $\mathrm{supp}(\phi)\subseteq B_r(q)$ and $\phi\big|_{B_{(1- \tau)r}}\equiv 1$ for some $\tau\in (0, \frac{1}{2})$, furthermore $|\Delta \phi|\leq C(n)\tau^{-8}r^{-2}$. Then from the Bochner formula and $\Delta f= 2n$, we have
\begin{align}
&\quad \int_{B_r}\Delta\phi\cdot (|\nabla f|^2- 4f)= \int_{B_r}\phi\cdot \Delta(|\nabla f|^2- 4f) \nonumber \\
&= \int_{B_r}\phi\cdot \Big(2|\nabla^2 f|^2+ 2\langle\nabla \Delta f, \nabla f\rangle+ Rc(\nabla f, \nabla f)- 8n\Big)\nonumber \\
&\geq 2\int_{B_r}\phi\cdot(|\nabla^2 f|^2- 4n) = 2\int_{B_r}\phi\cdot \big|\nabla^2 f- 2g\big|^2 ,\nonumber \\
&\geq 2\int_{B_{(1- \tau)r}} \big|\nabla^2 f- 2g\big|^2 .\nonumber
\end{align}

On the other hand, from (\ref{L2 diff of F}) and Proposition \ref{prop L2 diff of gadient of F},
\begin{align}
&\int_{B_r}\Delta\phi\cdot (|\nabla f|^2- 4f)\leq \sup_{B_r(q)}|\Delta\phi|\cdot \Big\{\int_{B_r} \big||\nabla f|^2- 4\rho^2\big|+ 4\int_{B_r}|\rho^2- f|\Big\} \nonumber \\
&\leq C(n)\tau^{-8}r^{-2}\cdot \Big\{\int_{B_r} |\nabla (f- \rho^2)|\cdot |\nabla(f+ \rho^2)|+ C(n)r^2\sqrt{\omega}\mathrm{V}(B_r)\Big\}  \nonumber \\
&\leq C(n)\tau^{-8}r^{-2}\cdot \Big\{\sqrt{r^2\omega\mathrm{V}(B_r)}\cdot \Big(\sqrt{\int_{B_r} |\nabla(f- \rho^2)|^2}+ \sqrt{\int_{B_r} |\nabla(2\rho^2)|^2}\Big)+ r^2\sqrt{\omega}\mathrm{V}(B_r)\Big\}  \nonumber \\
&\leq C(n)\tau^{-8}\cdot \sqrt{\omega}\mathrm{V}(B_r)  \nonumber
\end{align}

Let $\tau= \omega^{\frac{1}{32}}$, then we have
\begin{align}
\int_{B_r}\Delta\phi\cdot (|\nabla f|^2- 4f)\leq C(n) \omega^{\frac{1}{4}}\mathrm{V}(B_r)\leq C(n)\omega^{\frac{1}{4}}\mathrm{V}(B_{(1- \tau)r}). \nonumber
\end{align}

The conclusion follows from the above.
}
\qed

\begin{lemma}\label{lem volume conv imply level set conv}
{If $Rc(M^n)\geq 0$ and $V\big(B_r(p)\big)\geq (1- \delta)V\big(B_r(0)\big)$, then $\displaystyle V(s)\geq (1- \sqrt{\delta})n\omega_n s^{n- 1}$ for any $s\in (0, (1- \sqrt{\delta})^{\frac{1}{n}}r)$.
}
\end{lemma}

\pf
{Let $s_0= (1- \sqrt{\delta})^{\frac{1}{n}}r$, if $V(s_0)< (1- \sqrt{\delta})V(\mathbb{S}^{n- 1})s_0^{n- 1}$, then from the Bishop-Gromov Comparison Theorem,
\begin{align}
V(s)< (1- \sqrt{\delta})V(\mathbb{S}^{n- 1})s^{n- 1} \ , \quad \quad \quad \quad \forall s\in [s_0, r] \nonumber
\end{align}

Then we have
\begin{align}
V\big(B_r(p)\big)&= V\big(B_{s_0}(p)\big)+ V(A_{s_0, r})\leq V\big(B_{s_0}(0)\big)+ \int_{s_0}^r V(t)dt \nonumber \\
&< \frac{V(\mathbb{S}^{n- 1})}{n}s_0^n+ (1- \sqrt{\delta})V(\mathbb{S}^{n- 1})\int_{s_0}^r t^{n- 1} dt \nonumber \\
&= \frac{V(\mathbb{S}^{n- 1})}{n}\big[r^n- (r^n- s_0^n)\sqrt{\delta}\big]\leq \frac{V(\mathbb{S}^{n- 1})}{n}(1- \delta)r^n \nonumber \\
&= (1- \delta)V\big(B_r(0)\big) \nonumber
\end{align}
which is contradicting to the assumption. Then $\displaystyle V(s_0)\geq (1- \sqrt{\delta})V(\mathbb{S}^{n- 1})s_0^{n- 1}$.
}
\qed

\begin{definition}\label{def choice of f wrt q p r}
{Assume there is $q\in M^n$ with $d(p, q)= \beta r_1$, we define $f$ by
\begin{equation}\nonumber
\left\{
\begin{array}{rl}
&\Delta  f= 2n, \quad \quad \quad \quad on \ B_{(\beta+ 8)r_1}(q), \\
&f\big|_{\partial B_{(\beta+ 8)r_1}(q)}= (\beta+ 8)^2r_1^2\nonumber
\end{array} \right.
\end{equation}
We call $f$ the \textbf{model function with respect to $\{p, q, \beta\}$}.
}
\end{definition}

Given the volume ratio lower bound and the existence of $q$,  combining the former results in this section,  we obtain the integral estimate of difference as follows.

\begin{lemma}\label{lem difference controlled by volume ratio}
{For $\tau\in \big(0,\frac{1}{6(n+1)}\big)$, $3\leq \beta\leq \tau^{-c(n)}$ and $r>0$, assume $B_{r}(p)\subset (M,g)$ satisfying $Rc\ge 0$ and $\frac{V(B_{r}(p))}{V(B_{r})}\ge 1-\tau$. Assume there is $q\in M^n$ with $r_q\vcentcolon = \beta^{-1} d(p, q), d(p, q)\leq \beta^{-C(n)}r$, and define $\rho(x)= d(q, x)$. Then the model function $f$ with respect to $\{p, q, \beta\}$ satisfies
\begin{align}
r_q^{-1}\fint_{B_{2r_q}(p)} |\nabla (f- \rho^2)|+ r_q^{-2}\sup_{B_{2r_q}(p)}\Big|f- \rho^2\Big|+ \fint_{B_{4r_q}(p)} \big|\nabla^2 f- 2g\big| \leq C(n)\beta^{-C(n)} . \nonumber
\end{align}
}
\end{lemma}

\pf
{Note $r\geq (2\beta+ 16) r_q$, from the Bishop-Gromov volume comparison Theorem, we obtain
\begin{align}
\frac{V(B_{(2\beta+ 16) r_q}(q))}{V(B_{(2\beta+ 16) r_q}(0))}&\geq \frac{V(B_{r}(q))}{V(B_{r}(0))}\geq \frac{V(B_{r- \beta r_q}(p))}{V(B_{r}(0))}\geq (1- \tau)(1- \frac{\beta r_q}{r})^n \nonumber \\
&\geq  (1- \beta^{-C(n)})^{n+ 1}.\nonumber
\end{align}
In the last inequality above, we use $\displaystyle d(p, q)\leq \beta^{-C(n)}r$ and $\displaystyle \beta\leq \tau^{-c(n)}$.

From the above inequality, we have
\begin{align}
\frac{V(B_{(2\beta+ 16) r_q}(q))}{V(B_{(2\beta+ 16) r_q}(0))}\geq 1- (n+ 1)\beta^{-C(n)}. \nonumber
\end{align}

Let $\delta= (n+ 1)\beta^{-C(n)}$, apply Lemma \ref{lem volume conv imply level set conv} on $B_{(2\beta+ 16)r_q}(q)$, using $\delta\leq (1- 2^{-n})^2$, we obtain
\begin{align}
\frac{V(\partial B_s(q))}{s^{n- 1}}\geq (1- \sqrt{\delta})n\omega_n, \quad \quad \quad \quad \forall s\in (0, (\beta+ 8)r_q) . \label{boundary volume ineq}
\end{align}

From (\ref{boundary volume ineq}), we can apply Proposition \ref{prop L2 diff of gadient of F}, Proposition \ref{prop C0 diff of F} and Proposition \ref{prop integral Hessian small} on $B_{(\beta+ 8)r_q}(q)$ to get
\begin{align}
&\fint_{B_{(\beta+ 8)r_q}(q)} |\nabla (f- \rho^2)|^2 \leq C(n)\sqrt{\delta} \big((\beta+8)r_q\big)^2\le C(n)\beta^{- C(n)}r_q^2, \nonumber \\
&\sup_{B_{(1- 2\delta^{\frac{1}{4n+ 8}})(\beta+ 8)r_q}(q)}\Big|f- \rho^2\Big|\leq C(n)\big((\beta+8)r_q\big)^2\delta^{\frac{1}{4n+ 8}}\le C(n)\beta^{- C(n)}r_q^2 , \nonumber \\
&\fint_{B_{(1- \delta^{\frac{1}{64}})(\beta+ 8)r_q}(q)} \big|\nabla^2 f- 2g\big|^2 \leq C(n)\delta^{\frac{1}{8}}\leq C(n)\beta^{- C(n)}.  \nonumber
\end{align}

Note $\delta= (n+ 1)\beta^{-C(n)}$, then we have $\displaystyle \frac{\beta+ 4}{\beta+ 8}\leq \min\{1- 2\delta^{\frac{1}{4n+ 8}}, 1- \delta^{\frac{1}{64}}\}$. Now we obtain
\begin{align}
B_{4r_q}(p)\subseteq B_{(\beta+ 4)r_q}(q)\subseteq B_{(1- 2\delta^{\frac{1}{4n+ 8}})(\beta+ 8)r_q}(q)\cap B_{(1- \delta^{\frac{1}{64}})(\beta+ 8)r_q}(q) . \nonumber
\end{align}
Hence by Bishop-Gromov volume comparison Theorem, we get
\begin{align}
&\fint_{B_{2r_q}(p)} |\nabla (f- \rho^2)|^2 \leq C(n)\beta^{-C(n)} r_q^2\frac{V(B_{(\beta+ 8)r_q}(q))}{V(B_{2r_q}(p))}\leq C(n)\beta^{-C(n)}r_q^2, \nonumber  \\
&\sup_{B_{2r_q}(p)}\Big|f- \rho^2\Big|\leq C(n)\beta^{-C(n)}r_q^2\leq C(n)\beta^{-C(n)}r_q^2, \nonumber  \\
&\fint_{B_{4r_q}(p)} \big|\nabla^2 f- 2g\big|^2 \leq C(n)\beta^{-C(n)} \frac{V(B_{(1- \delta^{\frac{1}{64}})(\beta+ 8)r_q}(q))}{V(B_{4r_q}(p))}\leq C(n)\beta^{-C(n)} .\nonumber
\end{align}

From the H\"older inequality, we get
\begin{align}
\fint_{B_{2r_q}(p)} |\nabla (f- \rho^2)| \leq C(n)\beta^{-C(n)}r_q, \quad \quad \quad \fint_{B_{4r_q}(p)} \big|\nabla^2 f- 2g\big| \leq C(n)\beta^{-C(n)}. \nonumber
\end{align}
}
\qed

\subsection{The measure of points with integral control}

For $d(p, q)= \beta r> r$, define $\mathcal{B}(q)\subset B_{r}(p)$ is the set of points $x$, where $\pi(x)$ is not well defined and $\pi(x)= \overline{q,x}\cap \partial B_q(\beta r)$. When the context is clear, we use $\mathcal{B}$ instead of $\mathcal{B}(q)$ for simplicity.

Let $r_q=\beta^{-1}d(p,q)>0$, recall $\mathcal{B}(q)\subset B_{r_q}(p)$ be the set of points $x$, where $\pi(x)$ is not well defined and $\pi(x)= \overline{q,x}\cap \partial B_q(\beta r_q)$.

\begin{lemma}\label{lem dist between points and segment}
For $\tau\in (0,\frac{1}{6(n+1)}), 3\leq \beta\leq \tau^{-c(n)}, \epsilon>0$ and $r>0$, if $B_{r}(p)\subset (M,g)$ satisfies $Rc\ge 0$ and $\displaystyle \frac{V(B_{r}(p))}{V(B_r(0))}\ge 1-\tau$. Then for $q\in B_{\beta^{-C(n)} r}(p)\backslash \{p\}$ and $x\in B_{r_q}(p)$, we have
\begin{align*}
\frac{V(B_{\epsilon r_q}(x)\cap \mathcal{B}(q))}{V(B_{\epsilon r_q}(x))}\le C(n)\frac{\beta^{-C(n)}}{\epsilon^{n-1}}.
\end{align*}
\end{lemma}

\begin{proof}
Let $\delta= (n+1)\beta^{-C(n)}$, then we have
\begin{align}
\frac{V(B_{r}(q))}{V(B_{r}(0))}&\geq \frac{V(B_{(1- \beta^{-C(n)})r}(p))}{V(B_{r}(0))}\geq (1- \tau)(1- \beta^{-C(n)})^n\ge 1-(n+1)\beta^{-C(n)}=1-\delta. \nonumber
\end{align}
For $\theta\in S(T_qM)$, we define
$\displaystyle l_\theta=\sup\{l\ |\ \exp_q(t\theta)\text{ is a segment on } [0,l]\}$, let $\mathcal{A}=\{\theta\in S(T_{q} M)| \ l_{\theta}\le 3\beta r_q\}$ and set
$\alpha:=\frac{\mathcal{H}^{n-1}(\mathcal{A})}{\mathcal{H}^{n-1}(S^{n-1})}=\frac{\mathcal{H}^{n-1}(\mathcal{A})}{n\omega_n}.$ Then by Bishop-Gromov volume comparison theorem, we have
\begin{align*}
(1- \delta)\omega_nr^{n}&\le V(B_{r}(q))\le \int_{0}^{r}\int_{S(T_{q}M)\backslash \mathcal{A}}t^{n-1}dtd\theta+\int_{0}^{3\beta r_q}\int_{\mathcal{A}}t^{n-1}dtd\theta\\
&\le \omega_nr^{n}-\alpha \omega_n(r^{n}- (3\beta r_q)^n).
\end{align*}
So by $\beta r_q=d(p,q)\le \beta^{-C(n)} r$, we have
  $$\alpha\le \frac{\delta}{1-(3\beta^{-C(n)} )^n}\le \frac{(n+1)\beta^{-C(n)}}{1-3\beta^{-C(n)} }\le C(n)\beta^{-C(n)} .$$

Now, For $r_2\geq r_1\geq 0$, we define
\begin{align*}
&A_{r_1, r_2}(q)= \{y\in M^n|\ r_1< d(y, q)< r_2\} ,\nonumber \\
&\mathcal{C}=\{y\in A_{(\beta-\epsilon)r_q, (\beta+\epsilon)r_q}(q)|  \exists \theta(y)\in S(T_{q}M) s.t.  y=\exp_q(\rho(y)\theta(y)), l_{\theta(y)}\ge 3\beta r_q\}.
  \end{align*}

  Note $\rho(x)\le d(x,p)+d(p,q)\le (\beta+1)r_q\le 2\beta r_q$.  Thus by Bishop-Gromov volume comparison theorem again, we get
  \begin{align*}
  \frac{V(A_{\rho(x)-\epsilon r_q,\ \rho(x)+\epsilon r_q}(q)\backslash \mathcal{C})}{V(B_{\epsilon r_q}(x))}&\leq \frac{\int_{\rho(x)-\epsilon r_q}^{\rho(x)+\epsilon r_q}\int_{\mathcal{A}}t^{n-1}dtd\theta}{V(B_{\epsilon r_q}(x))}\nonumber \\
& \le C(n)\alpha\omega_n \frac{(\rho(x)+\epsilon r_q)^n-(\rho(x)-\epsilon r_q)^n}{(1-(n+1)\delta)\omega_n \epsilon^n}\\
  &\le C(n)\beta^{-C(n)}  \frac{n\cdot 2^n \big(\frac{\rho(x)}{r_q}\big)^{n-1}}{\epsilon^{n-1}}\le C(n)\frac{\beta^{-C(n)}}{\epsilon^{n-1}}.
  \end{align*}

Now the conclusion follows from the fact $\displaystyle B_{\epsilon r_q}(x)\cap \mathcal{B}(q)\subseteq A_{\rho(x)-\epsilon r_q,\ \rho(x)+\epsilon r_q}(q)\backslash \mathcal{C}$.
\end{proof}

For $d(p, q)= \beta r> r$, we define the points of the ball $B_r(p)\subseteq B_{(\beta+ 1)r}(q)$ with integral control as follows.

\begin{definition}\label{def local integral has good control inner ball}
{For $0< \eta< \frac{1}{2}$ and $L^\infty$ function $h\geq 0$ in $B_{4r_q}(p)$, we define
\begin{align}
Q_{\eta}^{r_q}(h)&= \Big\{x\in B_{r_q}(p)\backslash \mathcal{B}: \int_
{0}^{|\rho(x)- \beta r_q|} h\big(\sigma_x(s)\big)ds\leq \eta r_q^2\Big\}\ , \nonumber \\
\check{Q}_{\eta}^{r_q}(h)&= \Big\{x\in B_{r_q}(p)\backslash \mathcal{B}: \int_
{0}^{|\rho(x)- \beta r_q|} h\big(\sigma_x(s)\big)ds> \eta r_q^2\Big\}\ , \nonumber \\
T_{\eta}^{r_q}(h)&= \Big\{x\in B_{r_q}(p)\backslash \mathcal{B}: \frac{1}{V(B_{r_q}(p))}\int_{B_{r_q}(p)\backslash \mathcal{B}} dy\Big(\int_0^{|\rho(x)- \beta r_q|} \big(\int_{\gamma_{\sigma_x(s), \tilde{\sigma}_y(s)}} h\big) ds\Big)\leq \eta r_q^2\Big\} \ , \nonumber \\
\check{T}_{\eta}^{r_q}(h)&= \Big\{x\in B_{r_q}(p)\backslash \mathcal{B}: \frac{1}{V(B_{r_q}(p))}\int_{B_{r_q}(p)\backslash \mathcal{B}} dy\Big(\int_0^{|\rho(x)- \beta r_q|} \big(\int_{\gamma_{\sigma_x(s), \tilde{\sigma}_y(s)}} h\big) ds\Big)> \eta r_q^2\Big\} \ , \nonumber \\
T_{\eta}^{r_q}(h, x)&= \Big\{y\in B_{r_q}(p)\backslash \mathcal{B}: \int_
0^{|\rho(x)- \beta r_q|} \big(\int_{\gamma_{\sigma_x(s), \tilde{\sigma}_y(s)}} h\big) ds\leq \sqrt{\eta}r_q^2 \Big\}\ , \quad \quad \quad \quad \forall x\in T_{\eta}^{r_q}(h) .\nonumber\\
\check{T}_{\eta}^{r_q}(h, x)&= \Big\{y\in B_{r_q}(p)\backslash \mathcal{B}: \int_
0^{|\rho(x)- \beta r_q|} \big(\int_{\gamma_{\sigma_x(s), \tilde{\sigma}_y(s)}} h\big) ds> \sqrt{\eta}r_q^2 \Big\}\ , \quad \quad \quad \quad \forall x\in T_{\eta}^{r_q}(h) .\nonumber
\end{align}
}
\end{definition}

We recall the following segment inequality due to Cheeger and Colding \cite{CC-Ann}.
\begin{lemma}[Segment Inequality]\label{lem segment ineq}
{Assume $(M^n, g)$ is a complete Riemannian manifold with $Rc\geq 0$, then for any nonnegative function $h$ defined on $B_{2r}(q)\subset M^n$,
\begin{align}
\int_{B_r(q)\times B_r(q)} \Big(\int_0^{d(y_1, y_2)} h\big(\gamma_{y_1, y_2}(s)\big) ds\Big) dy_1 dy_2\leq 2^{n+ 1}r\cdot V\big(B_r(q)\big)\cdot \int_{B_{2r}(q)} h \nonumber
\end{align}
}
\end{lemma}\qed

Now we use the segment inequality to show that the set of points with integral control occupies a large part of the inner ball.

\begin{lemma}\label{lem lower bound of good local integra set inner ball}
{For $\beta\ge 3$ and $r>0$, assume $B_{ r}(p)\subset (M,g)$ satisfies $Rc\ge 0$ and $r_q=\beta^{-1}d(p,q)$, then we have
\begin{align}
&\frac{V(\check{Q}_{\eta}^{r_q}(h))}{V\big(B_{r_q}(p)\big)}\leq \frac{C(n)}{\eta r_q}\fint_{B_{2r_q}(p)} h  , \quad \quad \quad
\frac{V(\check{T}_{\eta}^{r_q}(h))}{V\big(B_{r_q}(p)\big)}\leq \frac{C(n)}{\eta} \cdot \fint_{B_{4r_q}(p)}h  , \nonumber \\
&\frac{V\big(\check{T}_{\eta}^{r_q}(h, x)\big)}{V\big(B_{r_q}(p)\big)}\leq \sqrt{\eta} , \quad \quad \quad \quad \forall x\in T_{\eta}^{r_q}(h)\nonumber
\end{align}
}
\end{lemma}

\pf
{\textbf{Step (1)}. Let $B= B_{r_q}(p)$, we have

\begin{align}
\int_{ \check{Q}_\eta^{r_q}(h)} dx \int_0^{|\rho(x)- \beta r_q|} h\big(\sigma_x(s)\big) ds \geq \eta r_q^2\cdot V\big( \check{Q}_\eta^{r_q}(h)\big) .\nonumber
\end{align}

Assume $\theta_s(x)$ is the gradient flow of $\rho(\cdot)$ starting from $\pi(x)$ at time $s$ and denote
\begin{align*}
\mathcal{S}(t)=\left\{
\begin{aligned}
 1&, & \text{ if }t-\beta r_q>0, \\
-1&, & \text{ if } t-\beta r_q<0.
\end{aligned}
\right.
\end{align*}
 Using Co-Area formula and Bishop-Gromov volume comparison theorem, we get
\begin{align}
\quad \int_{ \check{Q}_\eta^{r_q}(h)} dx &\int_0^{|\rho(x)- \beta r_q|} h\big(\sigma_x(s)\big) ds
\leq \int_{(\beta- 1) r_q}^{(\beta+ 1) r_q} dt \int_0^{|t- \beta r_q|} ds \int_{\rho^{-1}(t)\cap B\backslash \mathcal{B}}  h\big(\theta_{\mathcal{S}(t)s}(x)\big) d\mathcal{H}^{n-1}(x)\nonumber\\
& \leq \big(\frac{\beta+ 1}{\beta- 1}\big)^{n- 1} \int_{(\beta- 1) r_q}^{(\beta+ 1) r_q} dt \int_0^{|t- \beta r_q|} ds \int_{\theta_{\mathcal{S}(t)s}\big(\rho^{-1}(t)\cap B\backslash \mathcal{B}\big)}   h\big(\tilde{x}\big) d\mathcal{H}^{n-1}(\tilde{x}) \nonumber \\
&\leq C(n) \int_{(\beta- 1) r_q}^{(\beta+ 1) r_q} dt  \int_{B_{2r_q}(p)}   h\big(\tilde{x}\big) d\tilde{x} \leq C(n)r_qV(B_{r_q}(p))\cdot \fint_{B_{2r_q}(p)}   h .\nonumber
\end{align}
In the above we use $(\frac{\beta+ 1}{\beta- 1})\leq 2$ by $\beta \geq 3$.  From the above, we have
\begin{align*}
\frac{V\big(\check{Q}_\eta^{r_q}(h)\big)}{V(B)}\leq \frac{C(n)}{\eta r_q}\fint_{B_{2r_q}(p)}h.
\end{align*}

To prove the $3$rd inequality of the conclusion, we note
\begin{align}
\int_{ \check{T}_{\eta}^{r_q}(h, x)} dy \int_0^{|\rho(x)- \beta r_q|}\Big(\int_{\gamma_{\sigma_x(s), \tilde{\sigma}_y(s)}} h\Big) ds\geq \sqrt{\eta}r_q^2 V\big( \check{T}_{\eta}^{r_q}(h, x)\big) \nonumber
\end{align}

On the other hand, note $x\in T_{\eta}^{r_q}(h)$, we have
\begin{align}
\int_{\check{T}_{\eta}^{r_q}(h, x)} dy \int_0^{|\rho(x)- \beta r_q|}\Big(\int_{\gamma_{\sigma_x(s), \tilde{\sigma}_y(s)}} h\Big) ds
\leq  V\big(B\big)\cdot \eta r_q^2 .\nonumber
\end{align}
Hence we obtain the $3$rd inequality.

\textbf{Step (2)}. Finally we prove the $2$nd inequality. Note we have

\begin{align}
\int_{ \check{T}_{\eta}^{r_q}(h)} dx\frac{1}{V(B)}\int_{B\backslash\mathcal{B}} dy \int_0^{|\rho(x)- \beta r_q|}\Big(\int_{\gamma_{\sigma_x(s), \tilde{\sigma}_y(s)}} h\Big) ds \geq \eta r_q^2 \cdot V\big( \check{T}_{\eta}^{r_q}(h)\big) .\nonumber
\end{align}

Let $\displaystyle \Omega_0= \big(\rho^{-1}(t_1)\cap B\backslash\mathcal{B}\big)\times \big(\rho^{-1}(t_2)\cap B\backslash\mathcal{B}\big)$, and
\begin{align}
\Omega_1= \theta_{\mathcal{S}(t_1)s}\big(\rho^{-1}(t_1)\cap B\backslash\mathcal{B}\big)\times \theta_{\frac{t_2- \beta r_q}{|t_1- \beta r_q|}\cdot s}\big(\rho^{-1}(t_2)\cap B\backslash\mathcal{B}\big). \nonumber
\end{align}
From the Co-area formula, we have
\begin{align}
&\quad\frac{1}{V(B)} \int_{\check{T}_{\eta}^{r_q}(h)} dx\int_{B\backslash\mathcal{B}} dy \int_0^{|\rho(x)- \beta r_q|}\Big(\int_{\gamma_{\sigma_x(s), \tilde{\sigma}_y(s)}} h\Big) ds \nonumber\\
&\leq \frac{1}{V\big(B\big)}\int_{(\beta- 1) r_q}^{(\beta+ 1) r_q} dt_1\int_{(\beta- 1) r_q}^{(\beta+ 1) r_q}  dt_2 \int_{\Omega_0} d\mathcal{H}^{2n-2}(x,y) \int_0^{|t_1- \beta r_q|}\Big(\int_{\gamma_{\sigma_x(s), \tilde{\sigma}_y(s)}} h\Big) ds \nonumber \\
&\leq \frac{C(n)}{V\big(B\big)}\int_{(\beta- 1) r_q}^{(\beta+ 1) r_q} dt_1\int_{(\beta- 1) r_q}^{(\beta+ 1) r_q}  dt_2\int_0^{|t_1- \beta r_q|} \Big(\int_{\Omega_1} d\mathcal{H}^{2n-2}(\tilde{x},\tilde{y}) \Big(\int_{\gamma_{\tilde{x}, \tilde{y}}} h\Big) \Big) ds \nonumber \\
&\leq \frac{C(n)}{V\big(B\big)}\int_{0}^{r_q} \int_{B_{2r_q}(p)\times B_{2r_q}(p)} d\tilde{x}d\tilde{y} \Big(\int_{\gamma_{\tilde{x}, \tilde{y}}} h\Big)ds.  \nonumber
\end{align}

Now from  Lemma \ref{lem segment ineq} and the Bishop-Gromov comparison Theorem, we get
\begin{align}
&\quad\frac{1}{V(B)} \int_{\check{T}_{\eta}^{r_q}(h)} dx\int_{B\backslash\mathcal{B}} dy \int_0^{|\rho(x)- \beta r_q|}\Big(\int_{\gamma_{\sigma_x(s), \tilde{\sigma}_y(s)}} h\Big) ds \leq C(n)r_q^2 \int_{B_{4r_q}(p)}h \nonumber \\
&\leq  C(n)r_q^2 V(B_{r_q}(p))\fint_{B_{4r_q}(p)}h,  \nonumber
\end{align}
and hence the $2$nd inequality follows.
}
\qed

\subsection{Stratified Gou-Gu Theorem for multiple points}

In the rest of the paper, we assume $\nu\geq 4^{n+ 1}$ is some fixed constant. To establish the stratified Gou-Gou Theorem with small error estimate, we have to choose points in the ball carefully, which satisfy the corresponding integral estimates. For this reason, we define a model $\epsilon r$-dense net in $B_r(p)$ as follows.
\begin{definition}\label{def epsilon r dense model sets}
{For $3\leq \beta$, assume $B_{r}(p)\subset (M,g)$. For some fixed $k$ with $1\leq k\leq n$, assume there is $\{q_i\}_{i= 1}^{k- 1}$ such that
\begin{align}
C(n)^{-1}\beta^{\nu^{n+ 1- i}}\leq \frac{d(p, q_i)}{d(p, q_{i+ 1})}\leq C(n)\beta^{\nu^{n+ 1- i}}, \quad \quad \quad \forall 1\leq i\leq k- 2. \nonumber
\end{align}
Define $r_{q_i}= \beta^{-\nu^{n+ 1- i}}d(p, q_i)$, $\rho_i(x)= d(q_i, x)$ and $\pi_i(x)= \overline{q_i x}\cap \partial B_{q_i}(d(p, q_i))$ for any $x\in B_{r_{q_i}}(p)$. For the model functions $f_i$ with respect to $\{p, q_i, \beta^{\nu^{n+ 1- i}}\}$ where $1\leq i\leq k- 1, \eta\in (0, \frac{1}{2})$, if there is an  $(\epsilon r_{q_i})$-dense subset  $\mathfrak{N}_i\subset B_{r_{q_i}}(p)$ and the family $\{\mathfrak{N}_i\}_{i=1}^{k-1}$ satisfies
\begin{enumerate}
\item[(a)]. For any $x\in \mathfrak{N}_i$ and $1\le l\le i$, $\pi_l(x)$ is well-defined;
\item[(b)]. For any $1\leq l\le i_1\le i_2\le k-1$, $y\in \mathfrak{N}_{i_1}, z\in \mathfrak{N}_{i_2}$, we have
\begin{align*}
\big|d(y, z)^2-|\rho_l(y),\rho_l(z)|^2-d(\pi_l(y),\pi_l(z))^2 \big|\le \tilde{\epsilon}_{l}(r_{q_l})^2;
\end{align*}
\item[(c)]. $\displaystyle \#(\mathfrak{N}_i)\le \frac{C(n)}{\epsilon^{n}}$;
 \item[(d)]. For any $1\le l\le i$, $\mathfrak{N}_i\subset T^{r_{q_l}}_{\eta, l}:=T_\eta^{r_{q_l}}(|\nabla^2f_l-g|)$ with respect to the function $\rho_l$;
\end{enumerate}
then we say $\{\mathfrak{N}_i\}_{i=1}^{k-1}$ is a \textbf{model $\{\epsilon r_{q_i}\}_{i=1}^{k-1}$-dense net with respect to $\Big\{\{ B_{r_{q_i}}(p), f_i, \tilde{\epsilon}_i\}_{i= 1}^{k- 1}, \eta\Big\}$}.
}
\end{definition}

The philosophy of finding $n$ direction points is starting from $\{\mathfrak{N}_i\}_{i=1}^{k-1}$, then we prove the corresponding stratified Gou-Gu Theorem, next we find direction point $q_k$ and the corresponding $\{\mathfrak{N}_i\}_{i=1}^{k}$. The final conclusion follows from the induction method based on the above argument.

Now we obtain the distance estimate between any points in a model $(\epsilon r)$-dense net of the ball $B_r(p)$ as follows by induction method.

\begin{prop}\label{prop dist of points in geodesic balls 1}
{For $\tau\in (0, c(n)]$, $\tau^{-c(n)}\geq \beta\ge 3$ and $r>0$, assume $B_{r}(p)\subset (M,g)$ satisfying $Rc\ge 0$ and $\frac{V(B_{r}(p))}{V(B_{r})}\ge 1-\tau$.  For some fixed $k$ with $1\leq k\leq n$, assume there are $\{q_i\}_{i= 1}^k$ satisfying $\displaystyle d(p, q_1)\leq \beta^{-C(n)}r$ and
\begin{align}
C(n)^{-1}\beta^{\nu^{n+ 1- i}}\leq \frac{d(p, q_i)}{d(p, q_{i+ 1})}\leq C(n)\beta^{\nu^{n+ 1- i}}, \quad \quad \quad \quad 1\leq i\leq k- 1 .\nonumber
\end{align}
Let $\{f_i\}_{i= 1}^k$ be the model functions with respect to $(p, q_i, \beta^{\nu^{n+ 1- i}})$ and
\begin{align}
\epsilon= \beta^{-\nu^{n+ 1}}, \quad \quad \quad \quad  \tilde{\epsilon}_{i}= C(n) \beta^{-\nu^{n+ 1- i}}. \nonumber
\end{align}
If there is a \textbf{model $\{\epsilon r_{q_i}\}_{i=1}^{k-1}$-dense net $\{\mathfrak{N}_i\}_{i=1}^{k-1}$ with respect to $\Big\{\{B_{r_{q_i}}(p), f_i, \tilde{\epsilon}_{i}\}_{i= 1}^{k- 1}, \beta^{-C(n)}\Big\}$}, then there is $\mathfrak{N}_{k}\subseteq B_{r_{q_k}}(p)$ such that $\{\mathfrak{N}_i\}_{i=1}^{k}$ is a \textbf{model $\{\epsilon r_{q_i}\}_{i=1}^{k}$-dense net with respect to $\Big\{\{B_{r_{q_i}}(p), f_i, \tilde{\epsilon}_{i}\}_{i= 1}^{k}, \beta^{-C(n)}\Big\}$}.
}
\end{prop}

\pf
\textbf{Step (1)}. From Lemma \ref{lem difference controlled by volume ratio}, for any $1\leq l\leq k$ we have
\begin{align}
\frac{\fint_{B_{2r_{q_l}}(p)}|\nabla (f_l- \rho_l^2)|}{ r_{q_l}}+ \fint_{B_{4r_{q_l}}(p)} |\nabla^2 f_l- 2g|+\frac{\sup\limits_{B_{2r_{q_l}}(p)}|\rho_l^2-f_l|}{(r_{q_l})^2}\leq C(n)\beta^{-C(n)} .\label{diff inequality together}
\end{align}

In the rest, let $\eta= \beta^{-C(n)}$. We firstly choose $B_{\epsilon r_{q_k}}(z_j)$ with $z_j\in B_{r_{q_k}}(p)$ such that
\begin{align}
B_{r_{q_k}}(p)\subseteq \bigcup_{j= 1}^{m(\epsilon)}B_{\epsilon r_{q_k}}(z_j) , \nonumber
\end{align}
where $\displaystyle m(\epsilon)\leq C(n)\epsilon^{-n}$ from the Bishop-Gromov comparison Theorem.

From Lemma \ref{lem dist between points and segment}, for $1\leq i\leq k$, we know
\begin{align}
\frac{V(B_{\epsilon r_{q_k}}(z_1)\cap  \mathcal{B}(q_i))}{V(B_{\epsilon r_{q_k}}(z_1))}\leq C(n)\frac{\frac{d(q_i, p)}{r}}{\big(\epsilon\frac{r_{q_k}}{r_{q_i}}\big)^{n-1}} \le C(n)\frac{\beta^{-C(n)}}{\epsilon^{n- 1}}. \label{star ineq}
\end{align}

Denote $T_{\eta,  i}^{r_{q_i}}= T_\eta^{r_{q_i}}(|\nabla^2f_i- 2g|)$ with respect to the function $\rho_i$, where $1\leq i\leq k$.  Similarly let $Q_{\eta, i}^{r_{q_i}}=Q_\eta^{r_{q_i}}(|\nabla(f_i-\rho_i^2)|)$ and $T^{r_{q_i}}_{\eta,i}(x)=T^{r_{q_i}}_\eta(|\nabla^2f_i-2g|,x)$ for $x\in T^{r_{q_i}}_{\eta,i}$ with respect to $\rho_i$, then by Lemma \ref{lem lower bound of good local integra set inner ball}, (\ref{star ineq}) and Bishop-Gromov volume comparison theorem, we know
\begin{align}
&\frac{V\Big(B_{\epsilon r_{q_k}}(z_1)- \bigcup_{i= 1}^k \Big[\mathcal{B}(q_i)\cup \check{T}_{\eta,  i}^{r_{q_i}}\cup \check{Q}_{\eta, i}^{r_{q_i}}\Big]-\bigcup_{i=1}^{k-1}\Big[ \cup_{x\in \mathfrak{N}_i\atop l\le i} \check{T}_{\eta,  l}^{r_{q_l}}(x)\Big]\Big)}{V(B_{\epsilon r_{q_k}}(z_1))\big)} \nonumber \\
&\ge \frac{V\Big(B_{\epsilon r_{q_k}}(z_1)- \bigcup_{i= 1}^k \mathcal{B}(q_i)\Big)}{V(B_{\epsilon r_{q_k}}(z_1))}- \frac{V\Big( \bigcup_{i= 1}^k \big(\check{T}_{\eta,  i}^{r_{q_i}}\cup \check{Q}_{\eta, i}^{r_{q_i}}\big)\Big)}{V(B_{\epsilon r_{q_k}}(z_1))}- \frac{V\Big( \bigcup_{i=1}^{k-1}\Big[ \cup_{x\in \mathfrak{N}_i\atop l\le i} \check{T}_{\eta,  l}^{r_{q_l}}(x)\Big]\Big)}{V(B_{\epsilon r_{q_k}}(z_1))}\nonumber \\
&\ge 1- C(n)\frac{\beta^{-C(n)}}{\epsilon^{n- 1}}-  \frac{C(n)}{\epsilon^n}\big(\frac{r_{q_1}}{r_{q_k}}\big)^n\sum_{i= 1}^k\{\frac{\fint_{B_{4r_{q_i}}(p)}|\nabla^2f_i-g|}{\eta}+\frac{\fint_{B_{2r_{q_i}}(p)}|\nabla(f_i-\rho_i^2)}{\eta r_{q_i}}|+ \epsilon^{-n}\sqrt{\eta}\}. \nonumber
\end{align}

It is direct to get $\displaystyle \frac{r_{q_1}}{r_{q_k}}\leq C(n)\beta^{\nu^{n+ 1}}$, then from the above and (\ref{diff inequality together}), we have
\begin{align}
&\frac{V\Big(B_{\epsilon r_{q_k}}(z_1)- \bigcup_{i= 1}^k \Big[\mathcal{B}(q_i)\cup \check{T}_{\eta,  i}^{r_{q_i}}\cup \check{Q}_{\eta, i}^{r_{q_i}}\Big]-\bigcup_{i=1}^{k-1}\Big[ \cup_{x\in \mathfrak{N}_i\atop l\le i} \check{T}_{\eta,  l}^{r_{q_l}}(x)\Big]\Big)}{V(B_{\epsilon r_{q_k}}(z_1))\big)} \nonumber \\
&\geq 1- C(n)\beta^{-C(n)}\epsilon^{1- n}- C(n)\epsilon^{-2n}\beta^{-C(n)} \geq 1- C(n)\epsilon^{-2n}\beta^{-C(n)} >0. \nonumber
\end{align}
where we use the choice of $\epsilon$ in the last inequality.

Thus there exists $\displaystyle x_1^{(k)}\in \bigcap_{i= 1}^k \big(T_{\eta, i}^{r_{q_i}}\cap Q_{\eta, i}^{r_{q_i}}\big)\cap \bigcap_{i=1}^{k-1}\big(\cap_{x\in \mathfrak{N}_i\atop l\le i}T^{r_{q_l}}_{\eta,l}(x)\big) \cap B_{\epsilon r_{q_k}}(z_1)\backslash \bigcup_{i= 1}^k \mathcal{B}(q_i)$.  From the definition of $\mathcal{B}(q_i)$, we know that $\pi_i(x_1^{(k)})$ is well defined for all $1\leq i\leq k$. Moreover, for $1\le l\le i\le k-1$ and $x\in \mathfrak{N}_i$, by Lemma \ref{lem property 2 of G-H appr} (in this case $\beta$ there will be $\beta^{\nu^{n+ 1- l}}$), we get
\begin{align*}
&\Big|d(x_1^{(k)}, x)^2- \big[\rho_l(x_1^{(k)})- \rho_l(x)\big]^2- d\big(\pi_l(x_1^{(k)}), \pi_l(x)\big)^2\Big| \\
&\le C(n)\{\int_{0}^{d(x_1^{(k)},\pi_l(x_1^{(k)}))}|\nabla(\rho_l^2-f_l)|(\sigma_{x_1^{(k)}}(s))ds
+\int_{0}^{d(x,\pi_l(x))}|\nabla(\rho_l^2-f_l)|(\sigma_{x}(s))ds\\
&+\int_0^{d(x,\pi_l(x))}\int_{\gamma_{\sigma_{x(s)}\tilde{\sigma}_{x_1^{(k)}}(s)}}|\nabla^2f_l-2g|ds
+ \sup_{B_{2r_{q_l}}(p)}|\rho_l^2-f_l|+ \beta^{-\nu^{n+ 1- l}}r_{q_l}^2\}\\
&\le C(n)(\sup_{B_{2r_{q_l}}(p)}|\rho_l^2-f_l|+ \beta^{- \nu^{n+ 1- l}}r_{q_l}^2+\sqrt{\eta}r_{q_l}^2)\le C(n)\beta^{-\nu^{n+ 1- l}}r_{q_l}^2.
\end{align*}

\textbf{Step (2)}. Assume $\displaystyle \mathfrak{N}_{k}= \{x_i^{(k)}\}_{i= 1}^{m(\epsilon)}$, we choose $x_j^{(k)}$ by induction on $j$. For $s\leq m(\epsilon)\leq C(n)\epsilon^{-n}$, assume $\{x_1^{(k)},  \cdots,  x_s^{(k)}\}\subset \mathfrak{N}_{k}$ are chosen such that
\begin{enumerate}
\item[(a)].  For each $1\leq j\leq s$,  the point $\displaystyle x_j^{(k)}\in \bigcap_{i= 1}^k \big(T_{\eta, i}^{r_{q_i}}\cap Q_{\eta, i}^{r_{q_i}}\big)\cap B_{\epsilon r_{q_k}}(z_j)$.
\item[(b)].  For $1\leq j\leq s, 1\leq l\leq k$, the point $\pi_l(x_j^{(k)})$ is well defined.
\item[(c)].  For any $1\le l\le i\le k$, $x\in \mathfrak{N}_i$ and  $x'\in \{x_1^{(k)},  \cdots,  x_s^{(k)}\}$, we have
\begin{align}
 \Big|d(x, x')^2- \big[\rho_l(x)- \rho_l(x')\big]^2- d\big(\pi_l(x), \pi_l(x')\big)^2\Big| \leq \tilde{\epsilon}_lr_{q_l}^2  \nonumber
\end{align}
\end{enumerate}

By Lemma \ref{lem lower bound of good local integra set inner ball}, (\ref{diff inequality together}) and Bishop-Gromov comparison theorem, we know
\begin{align*}
&\frac{V\Big(B_{\epsilon r_{q_k}}(z_{s+ 1})- \bigcup\limits_{i= 1}^k \Big[\mathcal{B}(q_i)\cup \check{T}_{\eta,  i}^{r_{q_i}}\cup \check{Q}_{\eta, i}^{r_{q_i}}\cup \big(\bigcup\limits_{j= 1}^s \check{T}_{\eta,  i}^{r_{q_i}}(x_j^{(k)})\big)\Big]-\bigcup\limits_{i=1}^{k-1}\Big[ \bigcup\limits_{x\in \mathfrak{N}_i\atop l\le i} \check{T}_{\eta,  l}^{r_{q_l}}(x)\Big]\Big)}{V(B_{\epsilon r_{q_k}}(z_{s+ 1}))\big)} \nonumber \\
&\geq 1- C(n)\beta^{-C(n)}\epsilon^{-2n} >0, \nonumber
\end{align*}
where the choice of $\epsilon$ is used in the last inequality. Thus there exists
\begin{align}
x_{s+ 1}^{(k)}\in B_{\epsilon r_{q_k}}(z_{s+ 1})- \bigcup_{i= 1}^k \Big[\mathcal{B}(q_i)\cup \check{T}_{\eta,  i}^{r_{q_i}}\cup \check{Q}_{\eta, i}^{r_{q_i}}\cup \big(\bigcup_{j= 1}^s \check{T}_{\eta,  i}^{r_{q_i}}(x_j^{(k)})\big)\Big]-\bigcup_{i=1}^{k-1}\Big[ \bigcup\limits_{x\in \mathfrak{N}_i\atop l\le i} \check{T}_{\eta,  l}^{r_{q_l}}(x)\Big]. \nonumber
\end{align}

For $1\le i\le k, x\in \mathfrak{N}_i$ and $l\le i$,  by Lemma \ref{lem property 2 of G-H appr}, we get
\begin{align*}
&\Big|d(x, x_{s+1}^{(k)})^2- \big[\rho_l(x)- \rho_l(x_{s+1}^{(k)})\big]^2- d\big(\pi_l(x), \pi_l(x_{s+1}^{(k)})\big)^2\Big| \\
&\le C(n)(\sup_{B_{2r_{q_l}}(p)}|\rho_l^2-f_l|+ \beta^{- \nu^{n+ 1- l}}r_{q_l}^2+\sqrt{\eta}r_{q_l}^2)\le \tilde{\epsilon}_lr_{q_l}^2.
\end{align*}
Finally,  the conclusion follows by the induction method.
\qed

By induction on the number of points $q_i$,  using the volume ratio and the almost Gou-Gu Theorem, we obtain the lower bound of the finite projection set's diameter, which yields the points $q_k$. Also using the above integral of difference, combining Proposition \ref{prop dist of points in geodesic balls 1},  we find suitable $(\epsilon r)$-dense net and the related almost Gou-Gu Theorem at the same time.

Once the above $\{q_i\}_{i=1}^k$ and $\{\mathfrak{N}_{i}\}_{i=1}^{k}$ are constructed, for $0\le i\le k-1$, we can define $\mathscr{P}_{0}^{(i)}:B_{r_{q_{i+1}}}\to \mathfrak{N}_{i+1}$ such that  for each $x\in B_{r_{q_{i+1}}}(p)$,

\begin{align}
 d(x, \mathscr{P}_0^{(i)}(x))\leq C(n)\beta^{-\nu^{n+ 1}} r_{q_{i+ 1}}. \nonumber
\end{align}
Especially, we use $\mathscr{P}_0$ to denote $\mathscr{P}_0^{(0)}$ for simplicity. And we define $\pi_0(x)=x$,
\begin{align}
  \hat{\pi}_i= \mathscr{P}_0^{(i)}\circ \pi_i,  \text{ and } \mathscr{P}_s= \hat{\pi}_s\circ \cdots \circ \hat{\pi}_1\circ \hat{\pi}_0, \text{ for } 0\le i,s\le k-1. \nonumber
\end{align}
We further define $\displaystyle \check{\mathscr{P}}_i= \mathscr{P}_0^{(i- 1)}\circ \pi_i\circ \mathscr{P}_{i- 1}$ for $1\le k$ and
\begin{align}
\phi^{(k- 1)}(x)= (\phi_1(x), \cdots, \phi_{k- 1}(x)), \quad \quad and \quad \quad \phi_j(x)= \rho_j(\mathscr{P}_{j- 1} (x))- \rho_j(p). \nonumber
\end{align}
We also use the notation $r_{q_0}=\beta^{-C(n)} r$ and $\mathscr{P}_{-1}(x)= \mathscr{P}_0^{(-1)}(x)=x$.

\begin{remark}\label{rem the meaning of the proj}
{Note when $x\in \mathfrak{N}_1$, if $\pi_{i- 1}\circ \cdots \circ\pi_1(x)$ exists and $\displaystyle \pi_{j}\circ \cdots \circ\pi_1(x)\in \mathfrak{N}_{j+ 1}$ for $1\leq j\leq i- 1$, then we can choose $\mathscr{P}_{i- 1}(x)= \pi_{i- 1}\circ \cdots \circ\pi_1(x)$.
}
\end{remark}

\begin{theorem}\label{thm dist of points in geodesic balls induc-dim}
{For $\tau\in \big(0, c(n)\big)$, $3\leq \beta\leq \tau^{-c(n)}$, assume $\frac{V(B_r(p))}{V(B_r(0))}\geq 1- \tau$,  then we can find $\{q_k\}_{k= 1}^n$ with $\displaystyle d(p, q_1)= \beta^{-C(n)} r, q_{k}\in \check{\mathscr{P}}_{k- 1} (B_{r_{q_{k- 1}}}(p))$ and
\begin{align}
C(n)^{-1}\beta^{\nu^{n+ 1- k}} \leq \frac{d(q_k, p)}{d(q_{k+ 1}, p)}\leq C(n)\beta^{\nu^{n+ 1- k}},  \nonumber
\end{align}

Furthermore, there is a model $(\beta^{-\nu^{n+ 1}} r_{q_i})$-dense net $\{\mathfrak{N}_i\}_{i=1}^{n}$ with respect to
\begin{align}
\Big\{\{B_{r_{q_i}}(p), f_i, \beta^{-\nu^{n+ 1- i}}\}_{i= 1}^{n}, \beta^{-C(n)}\Big\}; \nonumber
\end{align}
and for all $1\leq k\leq n$,
\begin{align}\label{Gougu}
\sup_{x, y\in B_{r_{q_{k}}}(p)}\Big|d(x, y)^2&- \sum\limits_{j= 1}^{k}\big[\rho_j(\mathscr{P}_{j- 1}(x))- \rho_j(\mathscr{P}_{j- 1}(y))\big]^2 \nonumber   \\
&- d(\check{\mathscr{P}}_{k}(x), \check{\mathscr{P}}_{k}(y))^2\Big| \leq C(n) \beta^{-\nu^{n+ 1- k}} r_{q_{k}}^2.
\end{align}
}
\end{theorem}

\begin{remark}
{In the rest of the paper, for simplicity, we use $\displaystyle \gamma_i= \beta^{-2^{-1}\nu^{n+ 1- i}}r_{q_i}$. For any $i, j$, we have the following facts used often:
\begin{align}
r_{q_{i+ 1}}\leq r_{q_i}, \quad \quad \quad \gamma_i\leq \gamma_{i+ 1}, \quad \quad \quad \beta^{-\nu^{n+ 1}}r_{q_i}\leq \gamma_i \quad \quad \quad \gamma_{j- 2}r_{q_{j- 2}}\leq \gamma_{j- 1}^2. \nonumber
\end{align}
}
\end{remark}

\pf
{\textbf{Step (1)}. We can choose $q_1$ such that $d(p, q_1)= \beta^{-C(n)} r$. From Proposition \ref{prop dist of points in geodesic balls 1}, there is $\displaystyle \mathscr{P}_{0} (x), \mathscr{P}_{0} (y)\in \mathfrak{N}_{1}$ satisfying
\begin{align}
&d(\pi_1\mathscr{P}_0(x), \check{\mathscr{P}}_{1}(x))+ d(\pi_1\mathscr{P}_0(y), \check{\mathscr{P}}_{1}(y))\leq C(n)\beta^{-\nu^{n+ 1}} r_{q_1}, \nonumber \\
&\Big|d(\mathscr{P}_{0}(x), \mathscr{P}_{0}(y))^2- \big[\rho_1(\mathscr{P}_{0}(x))- \rho_1(\mathscr{P}_{0}(y))\big]^2\nonumber\\
&\quad \quad \quad \quad \quad \quad \quad  - d(\pi_1\circ \mathscr{P}_{0}(x), \pi_1\circ\mathscr{P}_{0}(y))^2\Big|\leq C(n)\beta^{-\nu^{n}}r_{q_1}^2.\nonumber
\end{align}
Thus (\ref{Gougu}) holds for $k= 1$.

Now we prove the conclusion by induction on $k$. Assume $k\le n$ and the conclusion holds for $1\le i\le k-1$. For $1\leq i\leq k-1$ and $\displaystyle x, y\in B_{r_{q_i}}(p)$, we have
\begin{align}\label{j assumption Gougu}
\Big|d(x, y)^2&- \sum\limits_{j= 1}^{i}\big[\rho_j(\mathscr{P}_{j- 1}(x))- \rho_j(\mathscr{P}_{j- 1}(y))\big]^2 \nonumber   \\
&- d(\check{\mathscr{P}}_{i}(x), \check{\mathscr{P}}_{i}(y))^2\Big| \leq C(n) \beta^{-\nu^{n+ 1- i}} r_{q_{i}}^2.
\end{align}


Define $\displaystyle \tilde{r}_1= r_{q_{k- 1}}$, then $\displaystyle \tilde{r}_1 \in [C(n)^{-1}\beta^{n- k+ 1} r_1, C(n)\beta^{n- k+ 1} r_1]$ by the induction assumption. Note
\begin{align}
\phi^{(k- 1)}(B_{\tilde{r}_1}(p))\subseteq [-\tilde{r}_1, \tilde{r}_1]\times [-2\tilde{r}_1, 2\tilde{r}_1]\times \cdots \times [- 2^{k- 2}\tilde{r}_1, 2^{k- 2}\tilde{r}_1]\subseteq [-2^n \tilde{r}_1, 2^n \tilde{r}_1]^{k- 1}. \nonumber
\end{align}
Then for $\lambda=2^{-n(2n+3)}$,  we can find $\displaystyle \{z_j\}_{j= 1}^{({\frac{2^{n+2}}{\lambda}})^{k- 1}}\subseteq B_{\tilde{r}_1}(p)$ such that for any $x\in B_{\tilde{r}_1}(p)$, there is $z_{j_0}$ satisfying $\displaystyle |\phi^{(k- 1)}(x)- \phi^{(k- 1)}(z_{j_0})|\leq \lambda \tilde{r}_1$.

We define $\displaystyle r_0\vcentcolon= \max_{x\in \check{\mathscr{P}}_{k- 1} (B_{\tilde{r}_1}(p))}d(x, p)$. From (\ref{j assumption Gougu}), for any $x\in B_{\tilde{r}_1}(p)$, we get
\begin{align*}
& \quad d(x, z_{j_0}) \\
&\leq \sqrt{|\phi^{(k- 1)}(x)- \phi^{(k- 1)}(z_{j_0})|^2+ d(\check{\mathscr{P}}_{k- 1}(x), \check{\mathscr{P}}_{k- 1} (z_{j_0}))^2} \nonumber \\
&\quad + C(n)\beta^{-\frac{1}{2}\nu^{n+ 2- k}}\tilde{r}_1 \\
&\leq \lambda \tilde{r}_1+ 2r_0+ C(n)\beta^{-\frac{1}{2}\nu^{n+ 2- k}}\tilde{r}_1 + C(n)\beta^{-14\cdot \nu^{n- 1}} \tilde{r}_1\le 2(\lambda+\frac{r_0}{\tilde{r}_1})\tilde{r}_1,
\end{align*}
where in the last line we use
\begin{align}
r_{q_1}= \beta^{-\nu^n}d(p, q_1)\leq \beta^{-\nu^n}(\beta^{\nu^n}\cdots \beta^{\nu^{n+ 1- (k- 2)}})d(p, q_{k- 1}) \leq \beta^{2\cdot \nu^{n- 1}}\tilde{r}_1. \nonumber
\end{align}

This implies $x\in B_{2(\lambda+\frac{r_0}{\tilde{r}_1})\tilde{r}_1}(z_{j_0})$,  and hence $\displaystyle B_{\tilde{r}_1}(p)\subseteq \bigcup_{i= 1}^{(\frac{2^{n+2}}{\lambda})^{k-1}} B_{2(\lambda+\frac{r_0}{\tilde{r}_1})\tilde{r}_1}(z_i)$.

From the volume comparison Theorem, we obtain
\begin{align}
(1- \tau)\omega_n\tilde{r}_1^n\leq V(B_{\tilde{r}_1}(p))\leq (\frac{2^{n+2}}{\lambda})^{k-1} \omega_n \big(2(\lambda+\frac{r_0}{\tilde{r}_1})\tilde{r}_1\big)^n \nonumber
\end{align}
Using $k\leq n$, we get $\frac{1}{2}\le (\frac{2^{n+2}}{\lambda})^{n-1}\cdot 2^n\cdot (\lambda+\frac{r_0}{\tilde{r}_1})^n$. That is,
\begin{align}
\frac{r_0}{\tilde{r}_1}\ge 2^{-2(n+1)}\lambda^{\frac{n-1}{n}}-\lambda=(2^{-2(n+1)}-\lambda^{\frac{1}{n}})=2^{-n(2n+3)}. \nonumber
\end{align}

Now we choose $q_{k}\in \check{\mathscr{P}}_{k- 1}(B_{\tilde{r}_1}(p))$ with $\displaystyle d(p,  q_{k})\geq 2^{-n(2n+3)} \tilde{r}_1$.

\textbf{Step (2)}. Note we can use Proposition \ref{prop dist of points in geodesic balls 1} to find a model $(\beta^{-\nu^{n+ 1}} r_{q_i})$-dense net $\{\mathfrak{N}_i\}_{i=1}^{k}$ with respect to $\Big\{\{B_{r_{q_i}}(p), f_i, \tilde{\epsilon}_i\}_{i= 1}^{k}, \eta\Big\}$ by induction method.
Thus for $1\le i\le k$,  there exists $\mathscr{P}_0^{(i-1)}: B_{r_{q_i}}(p)\to \mathfrak{N}_i$ such that
\begin{align}
d(x, \mathscr{P}^{(i-1)}_0 (x))\leq C(n)\beta^{-\nu^{n+ 1}} r_{q_i} , \quad \quad \quad \quad \forall 1\leq i\leq k.
\end{align}

We show (\ref{j assumption Gougu}) also holds for $i=k$. Now for $x, y\in B_{r_{q_{k}}}(p)$, from $\mathscr{P}_{k-1} (x),  \mathscr{P}_{k-1} (y)\in \mathfrak{N}_k$, we can obtain
\begin{align}
&\Big|d(\mathscr{P}_{k-1} (x), \mathscr{P}_{k-1} (y))^2- \big[\rho_{k}(\mathscr{P}_{k-1} (x))- \rho_{k}(\mathscr{P}_{k-1} (y))\big]^2 \nonumber \\
&\quad \quad \quad \quad \quad \quad  - d(\pi_{k}\circ\mathscr{P}_{k-1} (x),  \pi_{k}\circ\mathscr{P}_{k-1} (y))^2\Big| \leq C(n)\beta^{-\nu^{n- k+1}} r_{q_{k}}^2.\label{j+1 terms}
\end{align}

We have
\begin{align}
&\quad |d(\mathscr{P}_{k-1}(x),\mathscr{P}_{k-1}(y))-d(\check{\mathscr{P}}_{k- 1}(x),\check{\mathscr{P}}_{k- 1}(y))|\nonumber \\
&\le d(\mathscr{P}_{k-1}(x),\check{\mathscr{P}}_{k- 1}(x))+d(\mathscr{P}_{k-1}(y),\check{\mathscr{P}}_{k- 1}(y))\le 2C(n)\beta^{-\nu^{n+1}}r_{q_{k- 1}},\label{perturb term}
\end{align}

Combining (\ref{j assumption Gougu}), (\ref{perturb term}) with (\ref{j+1 terms}), we get
\begin{align}
\Big|d(x, y)^2- \sum\limits_{j= 1}^{k}\big[\rho_j(\mathscr{P}_{j- 1}(x))- \rho_j(\mathscr{P}_{j- 1}(y))\big]^2 - d(\check{\mathscr{P}}_{k}(x), \check{\mathscr{P}}_{k}(y))^2\Big| \leq C(n) \beta^{-\nu^{n+ 1- k}} r_{q_{k}}^2.\nonumber
\end{align}

From the induction method,  the conclusion follows.
}
\qed

\section{The distance map is quasi-isometry}\label{sec quasi-isom}

The main result of this section is Theorem \ref{thm direction points imply Lipschitz less than 1+ep}, which says the distance map is a pseudo-isometry. The fact that upper bound of discrete Lipschitz constant is close to $1$, provides the distance counterpart of the almost orthogonality of distance functions determined by direction points. The angle counterpart of the almost orthogonality will be provided by Colding's integral Toponogov Theorem, which will be addressed in Subsection \ref{subsec almost og of dist map}.

There are two key points in the proof of Theorem \ref{thm direction points imply Lipschitz less than 1+ep}, which are Proposition \ref{prop n direction points} and Proposition \ref{prop the diam upper bound of (n+ 1)-proj}. We will prove Proposition \ref{prop n direction points} in Subsection \ref{subsec error of proj}, and this will be used to prove the upper bound of discrete Lipschitz constant of the distance map. Proposition \ref{prop the diam upper bound of (n+ 1)-proj} will be proved in Subsection \ref{subsec n+1 proj}, and this result provide the lower bound of discrete Lipschitz constant of the distance map.

\subsection{The error estimate of projections}\label{subsec error of proj}

In this subsection, we will show that the $n$ points $\{q_k\}_{k= 1}^n$ found in Theorem \ref{thm dist of points in geodesic balls induc-dim} combining with the `origin' point $p$ `almost' determines $n$ directions on $M^n$, which is reflected by the main result of this subsection: Proposition \ref{prop n direction points} (also see the proof of Theorem \ref{thm direction points imply Lipschitz less than 1+ep}).

\begin{lemma}\label{lem two key elem ineq}
For $a,b,c, \epsilon,\epsilon_2\ge 0$ and $\epsilon_1 \in [0,1)$, if
$$a+b\le (1+\epsilon_1)\big(\sqrt{a^2-c^2+\epsilon_2^2}+\sqrt{b^2-c^2+\epsilon_2^2}+\epsilon\big),$$
then $\displaystyle c\le 4\sqrt{(\epsilon+\epsilon_2)(a+b+\epsilon_2)}+4\sqrt{\epsilon_1}(a+b+\epsilon_2).$
\end{lemma}

\begin{proof}
Put $\xi=\frac{c}{a+b+2\epsilon_2}$. Then
\begin{align*}
a+b+2\epsilon_2\le (1+\epsilon_1)\sqrt{1-\xi^2}(a+b+2\epsilon_2)+(1+\epsilon_1)\epsilon+2\epsilon_2.
\end{align*}
If $(1+\epsilon_1)\sqrt{1-\xi^2}\ge 1$, then $\xi\le 2\sqrt{\epsilon_1}$ and $\displaystyle c= \xi(a+b+2\epsilon_2)\le 4\sqrt{\epsilon_1}(a+b+\epsilon_2),$
which implies the conclusion. Thus we can assume $(1+\epsilon_1)\sqrt{1-\xi^2}<1$ and from above, we get
$\displaystyle \xi^2(1+\epsilon_1)^2\le \frac{2\epsilon(1+\epsilon_1)+4\epsilon_2}{a+b+2\epsilon_2}+4\epsilon_1,$
which implies $\xi\le \frac{2\sqrt{\epsilon+\epsilon_2}}{\sqrt{a+b+2\epsilon_2}}+2\sqrt{\epsilon_1}$ and hence
$\displaystyle c\le 4\sqrt{(\epsilon+\epsilon_2)(a+b+\epsilon_2)}+4\sqrt{\epsilon_1}(a+b+\epsilon_2).$
\end{proof}

In Euclidean case, if two points lie on the projection image, then the line determined by the two points will also lie on the projection image. We will show the almost version of this fact on manifolds in Lemma \ref{lem crucial trig ineq two} and Lemma \ref{lem crucial trig ineq involved}.

There are two error estimates of projection, we firstly deal with the case with the projection point on the segment $\overline{x_1, x_2}$. The similar argument of the following lemma will be used in Subsection \ref{subsec n+1 proj} too. For $A\subseteq M$, we define $\displaystyle \mathfrak{U}_{\delta}(A)= \{x\in M: d(x, A)< \delta\}$.

\begin{lemma}\label{lem crucial trig ineq involved}
Assume $1\le l<i\le n$, $x_1\in B_{r_{q_i}}(p)$, $x_2\in B_{r_{q_{i-1}}}(p)$ and $j=i \text{ or } i-1$. If there are $\displaystyle\hat{x}_1, \hat{x}_2\in \bigcup_{s\ge l}\mathfrak{N}_{s}$ such that $\displaystyle \sum_{k=1}^2d(\hat{x}_k,x_k)+d(\pi_l(\hat{x}_k),\hat{x}_k)\le \tilde{\delta}$ and $d(x_1,x_2)\le r_{q_{j}}, $ where $\tilde{\delta}\in  (\gamma_{i-1}, r_{q_j})$ .
Then $\displaystyle \sup_{\hat{w}\in \mathfrak{U}_{\tilde{\delta}}(\overline{x_1, x_2})\cap \mathfrak{N}_i}d(\pi_l(\hat{w}),\hat{w})\le C(n)\big(\frac{\tilde{\delta}}{r_{q_j}}\big)^{\frac{1}{4}}r_{q_{j}}$.

\end{lemma}
\begin{figure}[H]
\begin{center}
\includegraphics{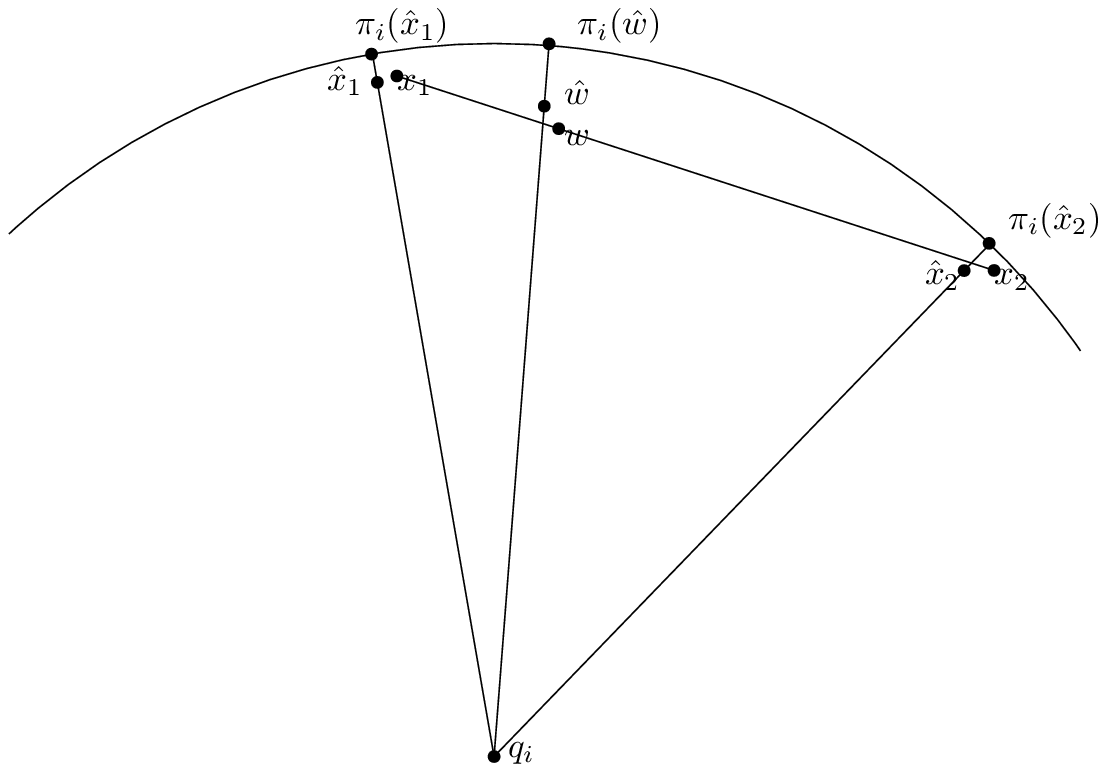}
\caption{Lemma \ref{lem crucial trig ineq involved}}
\label{figure: lemtrig}
\end{center}
\end{figure}

\begin{proof}
There is $w\in \overline{x_1, x_2}$ with $d(w, \hat{w})\leq \tilde{\delta}$. From $\overline{x_1, w, x_2}$ is a segment, we have
\begin{align*}
d(\hat{x}_1,\hat{w})+d(\hat{x}_2,\hat{w})
&\le d(x_1,w)+d(x_2,w)+4\delta=d(x_1,x_2)+4\tilde{\delta}\\
&\leq d(\pi_l(\hat{x}_1), \pi_l(\hat{w}))+ d(\pi_l(\hat{x}_2), \pi_l(\hat{w}))+ 8\tilde{\delta}.
\end{align*}

Applying the definition $(b)$ of $\mathfrak{N}_{s}$ on $\hat{x}_1, \hat{w}$ with respect to $q_l$, we obtain
\begin{align*}
d(\pi_l(\hat{x}_1), \pi_l(\hat{w}))\le \sqrt{d(\hat{x}_1,\hat{w})^2-|\rho_l(\hat{x}_1)-\rho_l(\hat{w})|^2+C(n)\gamma_l^2}.
\end{align*}
Similarly, we have
\begin{align*}
d(\pi_l(\hat{x}_2), \pi_l(\hat{w}))\le \sqrt{d(\hat{x}_2,\hat{w})^2-|\rho_l(\hat{x}_2)-\rho_l(\hat{w})|^2+C(n)\gamma_l^2}.
\end{align*}
Note that
\begin{align*}
\big | |\rho_l(\hat{w})-\rho_l(\hat{x}_2)|^2-|\rho_l(\hat{w})-\rho_l(\hat{x}_1)|^2\big|\le 8\tilde{\delta}(\tilde{\delta}+ d(x_1,x_2)).
 \end{align*}
 We have
 \begin{align*}
 d(\hat{w},\hat{x}_2)+d(\hat{w},\hat{x}_1)\le& \sqrt{d(\hat{w},\hat{x}_2)^2-|\rho_l(\hat{w})-\rho_l(\hat{x}_2)|^2+8\tilde{\delta}(\tilde{\delta}+ d(x_1,x_2))+C(n)\gamma_l^2}\\
 &+\sqrt{d(\hat{w},\hat{x}_1)^2-|\rho_l(\hat{w})-\rho_l(\hat{x}_2)|^2+C(n)\gamma_l^2}+8\tilde{\delta}.
  \end{align*}
  Thus by Lemma \ref{lem two key elem ineq} and note $\gamma_l\le \gamma_{i-1}\le \tilde{\delta}\le r_{q_j}$  and $d(x_1,x_2)\le r_{q_j}$, we get
  \begin{align*}
  |\rho_l(\hat{w})-\rho(\hat{x}_2)|\le C(n) \big(\frac{\tilde{\delta}}{r_{q_j}}\big)^{\frac{1}{4}}r_{q_j},
  \end{align*}
  which implies
  \begin{align*}
  d(\hat{w},\pi_l(\hat{w}))\le  C(n) \big(\frac{\tilde{\delta}}{r_{q_j}}\big)^{\frac{1}{4}}r_{q_j}+\tilde{\delta}\le C(n) \big(\frac{\tilde{\delta}}{r_{q_j}}\big)^{\frac{1}{4}}r_{q_j}.
  \end{align*}
\end{proof}

The second case deals with the projection point on the extension of the segment $\overline{x_1, x_2}$, which will be only used in the proof of Lemma \ref{lem crucial dist est between point and proj of pt}.

\begin{lemma}\label{lem crucial trig ineq two}
Assume $1\le l<i\le n$, $x_1\in B_{r_{q_i}}(p)$ and $x_2\in B_{r_{q_{i-1}}}(p)$. If there are $\hat{x}_1,\hat{x}_2\in \mathfrak{N}_{i-1}$ such that
$$\sum_{k=1}^2d(\hat{x}_k,x_k)+d(\pi_l(\hat{x}_k),\hat{x}_k)\le \delta r_{q_{i-1}}  \text{ and } d(x_1,x_2)\ge c(n) r_{q_{i-1}} ,$$
where $\delta \in (\epsilon_{i-1},\frac{1}{20}c(n))$ for $\epsilon_{i-1}=\frac{\gamma_{i-1}}{r_{q_{i-1}}}$ and $c(n)\in (0,1)$.
Then for $w\in B_{r_{q_i}}(p)$ with a segment  $\overline{w,x_1,x_2}$ and $\hat{w}\in B_{\delta r_{q_{i-1}}}(w)\cap \mathfrak{N}_i$, we have $\displaystyle d(\pi_l(\hat{w}),\hat{w})\le C(n)\delta^{\frac{1}{4}}r_{q_{i-1}}$.
\end{lemma}
\begin{figure}[H]
\begin{center}
\includegraphics{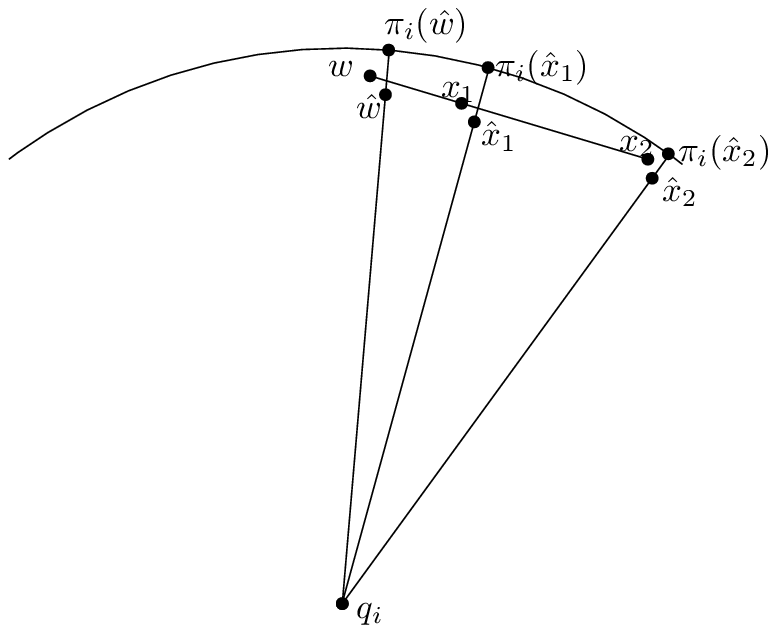}
\caption{Lemma \ref{lem crucial trig ineq two}}
\label{figure: lemtrigtwo}
\end{center}
\end{figure}

\begin{remark}
The lower bound of $d(x_1,x_2)$ is necessary even for Euclidean space $\mathbb{R}^2$. For example, for $q_1=(0,-1)$ and $q_2=(r_1,0)$ and $r_2<r_1$. There exist $x_1^k,x_2^k,w^k\in B_{r_2}(0)$ such that $\overline{w^k,x_1^k,x_2^k}$ is a segment and
$$d(x_1^k,\pi_1(x_1^k))+d(x_2^k,\pi(x_2^k))\le \frac{1}{k}r_1\to 0 \text{ but } \inf_{k}\frac{d(w^k,\pi_1(w^k))}{r_1}\ge \frac{r_2}{2r_1}>0.$$
For example, we can take $x_1^k=(0,-\frac{r_2}{4k}), x_2^k=(0,-\frac{r_2}{2k})$ and $w^k=(0,\frac{r_2}{2})$.
\end{remark}
\begin{proof}
By triangle inequality, we know
\begin{align}\label{triangle ineq}
d(\pi_l(\hat{w}),\pi_l(\hat{x}_2))-d(\pi_l(\hat{w}),\pi_l(\hat{x}_1))\le d(\pi_l(\hat{x}_1),\pi_l(\hat{x}_2))\le d(\hat{x}_1,\hat{x}_2)+2\delta r_{q_{i-1}}.
\end{align}
Since $\hat{x}_1,\hat{x}_2,\hat{w}\in \mathfrak{N}_i$ and $l<i$, we have
\begin{align*}
\big| d(\hat{w},\hat{x}_k)^2-d(\pi_l(\hat{w}),\pi_l(\hat{x}_k))^2-|\rho_l(\hat{w})-\rho_l(\hat{x}_k)|^2 \big|\le C(n)\gamma_l^2, \quad \quad k=1,2.
\end{align*}
Thus
\begin{align}\label{close projection}
d(\pi_l(\hat{w}),\pi_l(\hat{x}_1))\le \sqrt{d(\hat{w},\hat{x}_1)^2-|\rho_l(\hat{w})-\rho_l(\hat{x}_1)|^2+C(n)\gamma_l^2}
\end{align}
and
\begin{align}\label{far projection}
d(\pi_l(\hat{w}),\pi_l(\hat{x}_2))^2\ge d(\hat{w},\hat{x}_2)^2-|\rho_l(\hat{w})-\rho_l(\hat{x}_2)|^2-C(n)\gamma_l^2.
\end{align}
We claim the right hand term in the last inequality is positive, i.e.,
$$d(\hat{w},\hat{x}_2)^2-|\rho_l(\hat{w})-\rho_l(\hat{x}_2)|^2-C(n)\gamma_l^2>0.$$
Otherwise, we know $d(\hat{w},\hat{x}_2)^2-|\rho_l(\hat{w})-\rho_l(\hat{x}_2)|^2\le C(n)\gamma_l^2$  and hence
$$d(\pi_l(\hat{w}),\pi_l(\hat{x}_2))\le \sqrt{2C(n)}\gamma_l.$$
But when noting that $d(\hat{w},\pi_l(\hat{w}))=|\rho_l(\hat{w})-\rho_l(\pi_l(\hat{w}))|=|\rho_l(\hat{w})-\rho_l(p)|\le d(\hat{w},p)$ and  that $\overline{w,x_1,x_2}$ is a segment implies $d(x_1,x_2)\le d(x_2,w)$, we know
\begin{align*}
d(\pi_l(\hat{w}),\pi_l(\hat{x}_2))
&\ge d(\hat{x}_2,\hat{w})- d(\pi_l(\hat{w}),\hat{w})-\delta r_{q_{i-1}}\ge d(x_2,w)-d(\hat{w},p)-3\delta r_{q_{i-1}}\\
&\ge d(x_2,x_1)-d(w,p)-4\delta r_{q_{i-1}}\ge 5c(n)r_{q_{i-1}}\gg C(n) \gamma_{l},
\end{align*}
which is a contradiction.

So, substitute (\ref{close projection}) and the square root of (\ref{far projection}) into (\ref{triangle ineq}), we get
\begin{align}
\frac{A^2-B^2}{A+B}=A-B\le d(x_1,x_2)+4\delta r_{q_{i-1}}, \label{A B ineq}
\end{align}
where $A=\sqrt{d(\hat{w},\hat{x}_2)^2-|\rho_l(\hat{w})-\rho_l(\hat{x}_2)|^2-C(n)\gamma_l^2}$
and $$B=\sqrt{d(\hat{w},\hat{x}_1)^2-|\rho_l(\hat{w})-\rho_l(\hat{x}_1)|^2+C(n)\gamma_l^2}.$$
   Note $d(x_2, w)-d(x_1,w)=d(x_1,x_2)$ and
 \begin{align*}
\big | |\rho_l(\hat{w})-\rho_l(\hat{x}_2)|^2-|\rho_l(\hat{w})-\rho_l(\hat{x}_1)|^2\big|
&=|\rho_l(\hat{x}_1)-\rho_l(\hat{x}_2)|\cdot |\rho_l(\hat{x}_1)+\rho_l(\hat{x}_2)-2\rho_l(\hat{w})|\\
&\le (\sum_{k=1}^2d(\hat{x}_k,\pi_l(\hat{x}_k)))\cdot(\sum_{k=1}^2d(\hat{x}_k,\hat{w}))\\
(\text{ using }20\delta r_{q_{i-1}} \le d(x_1,x_2))&\le 18\delta r_{q_{i-1}} d(x_1,x_2) .
 \end{align*}
we have
 \begin{align*}
 A^2-B^2
 &=d(\hat{w},\hat{x}_2)^2-d(\hat{w},\hat{x}_1)^2+|\rho_l(\hat{w})-\rho_l(\hat{x}_1)|^2-|\rho_l(\hat{w})-\rho_l(\hat{x}_2)|^2-2C(n)\gamma_l^2\\
 &\ge \big( d(x_1,x_2)-4\delta r_{q_{i-1}}\big)\big((\hat{w},\hat{x}_2)+d(\hat{w},\hat{x}_1)\big)-18\delta r_{q_{i-1}} d(x_1,x_2) -C(n)\gamma_l^2 \\
 A+B&\le \sqrt{d(\hat{w},\hat{x}_2)^2-|\rho_l(\hat{w})-\rho_l(\hat{x}_2)|^2}\\
 &+\sqrt{d(\hat{w},\hat{x}_1)^2-|\rho_l(\hat{w})-\rho_l(\hat{x}_2)|^2+C(n)\gamma_l^2+18\delta r_{q_{i-1}} d(x_1,x_2) }.
 \end{align*}
Substituting these two inequalities into (\ref{A B ineq}), we get
 \begin{align*}
 &\big(1-\frac{4\delta r_{q_{i-1}}}{d(x_1,x_2)}\big)\big(d(\hat{w},\hat{x}_2)+d(\hat{w},\hat{x}_1)\big)\\
 \le& 18\delta r_{q_{i-1}}+C(n)\frac{\gamma_l^2}{d(x_1,x_2)}+\big(1+\frac{4\delta r_{q_{i-1}}}{d(x_1,x_2)}\big)\bigg(\sqrt{d(\hat{w},\hat{x}_2)^2-|\rho_l(\hat{w})-\rho_l(\hat{x}_2)|^2}\\
 +&\sqrt{d(\hat{w},\hat{x}_1)^2-|\rho_l(\hat{w})-\rho_l(\hat{x}_2)|^2+C(n)\gamma_l^2+18\delta r_{q_{i-1}}d(x_1,x_2)}\bigg)
 \end{align*}
 Note $d(x_1,x_2)\ge \max\{20\delta r_{q_{i-1}},\gamma_l\}$. We have
 \begin{align*}
 &d(\hat{w},\hat{x}_2)+d(\hat{w},\hat{x}_1)\le \big(1+\frac{16\delta r_{q_{i-1}}}{d(x_1,x_2)}\big)\bigg(\sqrt{d(\hat{w},\hat{x}_2)^2-|\rho_l(\hat{w})-\rho_l(\hat{x}_2)|^2}\\
 &+\sqrt{d(\hat{w},\hat{x}_1)^2-|\rho_l(\hat{w})-\rho_l(\hat{x}_2)|^2+C(n)\gamma_l^2+18\delta r_{q_{i-1}} d(x_1,x_2)}\bigg)+C(n)\big(\gamma_l+\delta r_{q_{i-1}}\big).
  \end{align*}
  Thus by Lemma \ref{lem two key elem ineq} and note $d(\hat{w},\hat{x}_2)+d(\hat{w},\hat{x}_1)\le 8d(x_1,x_2)$ , we get
  \begin{align*}
  |\rho_l(\hat{w})-\rho(\hat{x}_2)|\le C(n)\sqrt{\frac{\gamma_l}{d(x_1,x_2)}+\sqrt{\frac{\delta r_{q_{i-1}}}{d(x_1,x_2)}}}d(x_1,x_2),
  \end{align*}
  which implies
 \begin{align*}
  d(\hat{w},\pi_l(\hat{w}))
  &=|\rho_l(\hat{w})-\rho_l(\pi_l(\hat{w}))|=|\rho_l(\hat{w})-\rho_l(\pi_l(\hat{x}_2))|\\
  &\le d(\hat{x}_2,\pi_l(\hat{x}_2))+ |\rho_l(\hat{w})-\rho(\hat{x}_2)|\\
  &\le \delta r_{q_{i-1}}+ C(n)\sqrt{\frac{\gamma_l}{d(x_1,x_2)}+\sqrt{\frac{\delta r_{q_{i-1}}}{d(x_1,x_2)}}}d(x_1,x_2)\\
  &\le C(n)\delta^{\frac{1}{4}}r_{q_{i-1}},
  \end{align*}
  where we use  $\gamma_l\le \gamma_{i-1}\le \delta r_{q_{i-1}}\le \frac{1}{20}d(x_1,x_2)\le r_{q_{i-1}}$ in the last line.
\end{proof}

In Euclidean case, the image of the composition of different projections will lie on the image of any the projection and the composition of projections. The following lemma and Lemma \ref{lem proj and the composition of proj} show that this fact is almost true on manifolds.
\begin{lemma}\label{lem crucial dist est between point and proj of pt}
{For $2\leq j\leq n$, we have $\displaystyle \sup_{\hat{w}\in \check{\mathscr{P}}_{j- 1}(B_{r_{q_{j- 1}}}(p))\atop 1\leq i\leq j- 1} d(\hat{w},  \pi_i(\hat{w}))\leq C(n)\gamma_{j- 1}$.
}
\end{lemma}

\begin{remark}\label{rem reverse induction}
{In the proof of the above lemma,  we use the induction method.  Unlike the usual induction method from $k$ to $k+ 1$,  we do induction from $k+ 1$ to $k$ here and similar induction argument is also used in the later argument for other results involving projection. The reason is that the later projection depends on the later direction points, which lie on the image of the former projection by the choice of the direction points. And this asymmetric choice of direction points pushes us to use the above reverse induction argument.
}
\end{remark}

\begin{figure}[H]
\begin{center}
\begin{tabular}{c c }
		\includegraphics[width=0.42\linewidth]{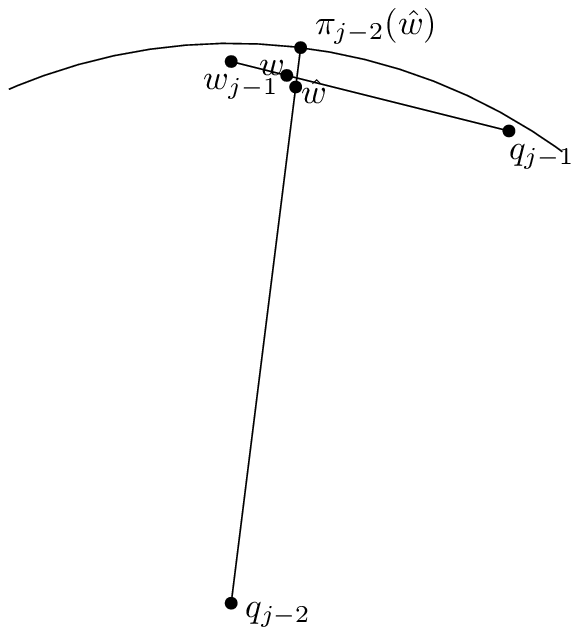}  &
        \hspace{0.05in}
		
		\includegraphics[width=0.50\linewidth]{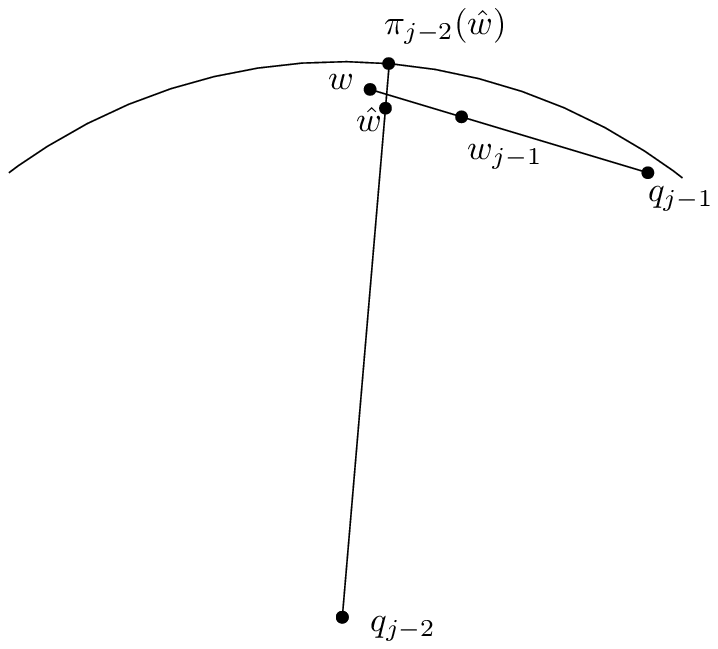}\\

		(Case $(1)$) &(Case $(2)$)\\
	\end{tabular}
\caption{Lemma \ref{lem crucial dist est between point and proj of pt}}
\label{figure: lem4.5}
\end{center}
\end{figure}

\pf
{ We argue by induction for $s=1,2, \ldots, n-1$.

\textbf{Step (1)}.  For $s=1$, we  will prove that for any $2=s+1\le j\le n$,
\begin{align*}
\sup_{\hat{w}\in \check{\mathscr{P}}_{j- 1}(B)\atop j- s\leq i\leq j- 1} d(\hat{w}, \pi_{i}(\hat{w}))=\sup_{\hat{w}\in \check{\mathscr{P}}_{j- 1}(B)} d(\hat{w}, \pi_{j-1}(\hat{w}))\leq C(n)\gamma_{j- 1} .
\end{align*}
   For simplicity, we use $B$ to denote $\displaystyle B_{r_{q_{j- 1}}}(p)$.  We firstly note there is $\displaystyle w\in \pi_{j- 1}\mathscr{P}_{j- 2}(B)$ such that $\displaystyle d(w,  \hat{w})\leq  \gamma_{j- 1}$.

From $\displaystyle d(w,  q_{j- 1})= d(\pi_{j- 1}(\hat{w}), q_{j- 1})= d(p, q_{j- 1})$, we get
\begin{align}
d(\hat{w}, \pi_{j- 1}(\hat{w}))&= |d(\hat{w}, q_{j- 1})- d(\pi_{j- 1}(\hat{w}), q_{j- 1})|=| d(\hat{w}, q_{j- 1})- d(w,  q_{j- 1})| \nonumber \\
&\leq d(\hat{w}, w)\leq C(n)\gamma_{j- 1}. \nonumber
\end{align}

\textbf{Step (2)}. Assume for  some $s\in \{1,2,\ldots, n-2\}$,  the following holds for any $j\in \{s+1,\ldots,n\}$:
\begin{align}
\sup_{\hat{w}\in \check{\mathscr{P}}_{j- 1}(B)\atop j- s\leq i\leq j- 1} d(\hat{w}, \pi_{i}(\hat{w}))\leq C(n)\gamma_{j- 1} .\label{induction assumption for proj}
\end{align}
By induction, to prove the conclusion we only need to show
\begin{align}
d(\hat{w}, \pi_{j- s- 1}(\hat{w}))\leq C(n)\gamma_{j- 1} \text{ for  all }  j\in \{s+2,\ldots, n\}.\nonumber
\end{align}

Assume $w= \pi_{j- 1}(w_{j- 1})$, where $w_{j- 1}\in \mathscr{P}_{j- 2}(B)$. By the definition of $\pi_{j-1}$, we know either $\overline{q_{j- 1}, w, w_{j- 1}}$ or $\overline{q_{j-1},w_{j-1},w}$ is a segment.  We discuss the two cases respectively.

Case (1). If $\overline{q_{j- 1}, w, w_{j- 1}}$ is a segment, then from the definition of $\mathscr{P}_{j- 2}, \check{\mathscr{P}}_{j- 2}$, we can find $\hat{w}_{j- 1}\in \check{\mathscr{P}}_{j- 2}$ such that
\begin{align}
d(w_{j- 1}, \hat{w}_{j- 1})\leq \beta^{-\nu^{n+ 1}}r_{q_{j- 2}} \le \gamma_{j-2}. \label{small perturbation}
\end{align}

Apply the induction assumption (\ref{induction assumption for proj}) for $j- 1$ and $\hat{w}_{j- 1}, q_{j- 1}\in \check{\mathscr{P}}_{j- 2}(B_{r_{q_{j- 2}}}(p))= \check{\mathscr{P}}_{(j- 1)- 1}(B_{r_{q_{j- 2}}}(p))$, we get
\begin{align}
d(\hat{w}_{j- 1}, \pi_{j- 1- s}(\hat{w}_{j- 1}))+ d(q_{j- 1}, \pi_{j- 1- s}(q_{j- 1}))\leq C(n)\gamma_{j- 2}. \label{ineq from induc assumption we have}
\end{align}

Applying Lemma \ref{lem crucial trig ineq involved} on $w_{j- 1}, q_{j- 1}$, note $\displaystyle \hat{w}\in \big(\mathfrak{U}_{\gamma_{j- 2}} \overline{w_{j- 1}q_{j- 1}}\big)\cap \mathfrak{N}_{j- 1}$, we get
\begin{align}
d(\hat{w}, \pi_{j- 1- s}(\hat{w}))\leq C(n)(\frac{\gamma_{j- 2}}{r_{q_{j- 2}}})^{\frac{1}{4}} r_{q_{j- 2}}\leq C(n)\gamma_{j- 1}. \nonumber
\end{align}

Case (2). If $\overline{q_{j- 1}, w_{j-1}, w_{j}}$ is a segment, we still choose $\hat{w}_{j-1}$ as above.  Applying Lemma \ref{lem crucial trig ineq two} on $w_{j-1}$ and $q_{j-1}$,we also get
\begin{align}
d(\hat{w}, \pi_{j- 1- s}(\hat{w}))\leq C(n)\gamma_{j- 1}. \nonumber
\end{align}
Then the conclusion follows from the induction method.
}
\qed

\begin{lemma}\label{lem proj and the composition of proj}
{For $\hat{w}\in \mathfrak{N}_k$ where $j\leq k\leq n$,  if $\displaystyle \sup_{1\leq i\leq j- 1}d(\hat{w}, \pi_i(\hat{w}))\leq \delta$, then
\begin{align}
\sup_{1\leq i\leq j- 1}d(\mathscr{P}_i(\hat{w}), \hat{w})\leq C(n)(\gamma_k+ \delta). \nonumber
\end{align}
}
\end{lemma}

\begin{remark}\label{rem choice of region}
{In the above lemma,  we require $\hat{w}\in \mathfrak{N}_k$ for some $k\geq j$,  to assure that $\mathscr{P}_{j- 1}(\hat{w})$ is well defined.
}
\end{remark}

\pf
{We prove the conclusion by induction on $i$.  When $i= 1$,
\begin{align}
d(\mathscr{P}_1(\hat{w}), \hat{w})\leq d(\pi_1(\hat{w}), \hat{w})+ \beta^{-\nu^{n+ 1}}r_{q_1}\leq C(n)(\gamma_k+ \delta). \nonumber
\end{align}

Assume $\displaystyle d(\mathscr{P}_m(\hat{w}), \hat{w})\leq C(n)(\gamma_k+ \delta)$ for some $i=m\leq j- 2$, then
\begin{align}
d(\mathscr{P}_{m+ 1}(\hat{w}), \hat{w})&\leq d(\pi_{m+ 1}\mathscr{P}_m(\hat{w}), \hat{w})+ \beta^{-\nu^{n+ 1}}r_{q_{m+ 1}} \nonumber \\
&\leq d(\pi_{m+ 1}\mathscr{P}_m(\hat{w}), \pi_{m+ 1}(\hat{w}))+ \delta+ C(n)\gamma_k \nonumber \\
&\leq \sqrt{d(\mathscr{P}_m(\hat{w}), \hat{w})^2+ C(n)\gamma_{m+ 1}^2}+ C(n)(\delta+ \gamma_k) \nonumber \\
&\leq  d(\mathscr{P}_m(\hat{w}), \hat{w})+ C(n)(\delta+ \gamma_k) \leq C(n)(\delta+ \gamma_k). \nonumber
\end{align}
The conclusion follows from the induction method.
}
\qed

After establishing the general results about projection error estimates, we can apply them on $q_j$ as follows.
\begin{cor}\label{cor comp of proj of qi}
{For $j\leq n$, we have $\displaystyle \sup_{1\leq i\leq j-1}d(q_j, \pi_i(q_j))+ \sup_{1\leq i\leq j- 2}d(\mathscr{P}_i(q_j), q_j)\leq C(n)\gamma_{j- 1}$.
}
\end{cor}

\pf
{From $q_j\in \check{\mathscr{P}}_{j- 1}(B)\subseteq \mathfrak{N}_{j- 1}$,  by Lemma \ref{lem crucial dist est between point and proj of pt}, we get
\begin{align}
\sup_{1\leq i\leq j-1}d(q_j, \pi_i(q_j))\leq C(n)\gamma_{j- 1}. \nonumber
\end{align}
Applying Lemma \ref{lem proj and the composition of proj},  we have
\begin{align}
\sup_{1\leq i\leq j- 2}d(\mathscr{P}_i(q_j), q_j)\leq C(n)(\gamma_{j- 1}+ C(n)\gamma_{j- 1})\leq C(n)\gamma_{j- 1}  \nonumber
\end{align}
}
\qed

\begin{prop}\label{prop n direction points}
{If $1\leq  i< j\leq n$,  we have $|\phi_i(q_j)|\leq C(n) \gamma_n$.
}
\end{prop}

\pf
{If $m\leq j- 1$,  using Corollary \ref{cor comp of proj of qi}, we have
\begin{align}
|\phi_m(q_j)|&= |\rho_m(\mathscr{P}_{m- 1}(q_j))- \rho_m(p)| \leq |\rho_m(q_j)- \rho_m(\pi_{m- 1}(q_j))|+ d(q_j, \mathscr{P}_{m- 1}(q_j))\nonumber \\
&\leq d(q_j, \pi_{m- 1}(q_j))+ C(n)\gamma_n \leq C(n)\gamma_{n}. \nonumber
\end{align}
}
\qed

\subsection{The diameter upper bound after $(n+ 1)$ projections}\label{subsec n+1 proj}

In this subsection, we show that there are exactly $n$ directions, which is implied by the following Proposition. This Proposition is only used to prove the lower bound in Theorem \ref{thm direction points imply Lipschitz less than 1+ep}, which implies the equivalence between G-H approximation and splitting map. The results of this subsection will not be used in our existence proof of the splitting map.

The key of the proof of Proposition \ref{prop the diam upper bound of (n+ 1)-proj} is the density of the $(n+ 1)$ projections image, that is,
 \begin{align*}
 \phi^{(n+ 1)}(B_{r_{q_{n+ 1}}}(p)) \text{ is } C(n)\beta^{-1}r_{q_{n+ 1}}\text{-dense  in } [-\frac{r_{q_{n+ 1}}}{2},0]^{n+1}\subset \mathbb{R}^{n+1},
 \end{align*}
which is the content of Proposition \ref{prop (n+1) comp proj is dense}. We prove Proposition \ref{prop the diam upper bound of (n+ 1)-proj} by firstly assuming Proposition \ref{prop (n+1) comp proj is dense}.

\begin{prop}\label{prop the diam upper bound of (n+ 1)-proj}
{$\displaystyle \sup_{x\in \check{\mathscr{P}}_n(B_{r_{q_n}}(p))} d(p, x)\leq \beta^{-1}r_{q_n}$.
}
\end{prop}

\begin{proof}
Assume there exists a point $q_{n+1}$ such that $r_0= d(q_{n+1},p)\ge \beta^{-1} r_{q_n}$. From Proposition \ref{prop (n+1) comp proj is dense}, the set $\phi^{(n+1)}(B_{r_{q_{n+1}}}(p))$ is $C(n)\beta^{-1} r_{q_{n+1}}$-dense in $[-\frac{r_{q_{n+1}}}{2},0]^{n+ 1}\subset \mathbb{R}^{n+1}$.

For $\lambda\ge C(n)\beta^{-1} r_{q_{n+1}}$ to be determined, by Proposition \ref{prop (n+1) comp proj is dense}, we can choose maximal disjoint family of balls $\{B_{\frac{\lambda}{2}}(y_i)\}_{i=1}^{N}$ in $\mathbb{R}^{n+1}$ such that $y_i\in \phi^{(n+1)}(B_{r_{q_{n+1}}}(p))$. Then  $\phi^{(n+1)}(B_{r_{q_{n+1}}}(p))\subset \cup_{i=1}^{N} B_{\lambda}(y_i)$ and hence
$$[-\frac{r_{q_{n+1}}}{2},0]^{n+ 1}\subset B_{C(n)\beta^{-1} r_{q_{n+1}}}\big( \phi^{(n+1)}(B_{r_{q_{n+1}}}(p))\big)\subset \bigcup_{i= 1}^NB_{\lambda+C(n)\beta^{-1} r_{q_{n+1}}}(y_i)\subset \bigcup_{i=1}^{N} B_{2\lambda}(y_i).$$

Thus we get
\begin{align*}
N\ge \frac{\big(\frac{r_{q_{n+1}}}{2}\big)^{n+ 1}}{\omega_{n+1}(2\lambda)^{n+1}}=\frac{r_{q_{n+1}}^{n+1}}{2^{2(n+1)}\omega_{n+1}\lambda^{n+1}}.
\end{align*}
Since $y_i\in \phi^{(n+1)}(B_{r_{q_{n+1}}}(p))$, there exists $x_i\in B_{r_{q_{n+1}}}(p)$ such that $\phi^{(n+1)}(x_i)=y_i$. By the construction of $\phi^{(n+1)}$, we know for $1\le i\neq j\le N$,
$$\lambda\le |y_i-y_j|=|\phi^{(n+1)}(x_i)-\phi^{(n+1)}(x_j)|\le \sqrt{n+1}d(x_i,x_j),$$
which means $B_{\frac{\lambda}{2\sqrt{n+1}}}(x_i)\cap B_{\frac{\lambda}{2\sqrt{n+1}}}(x_j)=\emptyset.$ Thus by the volume comparison theorem,
\begin{align*}
\omega_n(2r_{q_{n+1}})^n\ge V(B_{2r_{q_{n+1}}}(p))&\ge \sum_{i=1}^{N}V(B_{\frac{\lambda}{2\sqrt{n+1}}}(x_i))\\
&\ge \frac{r_{q_{n+1}}^{n+1}}{2^{2(n+1)}\omega_{n+1}\lambda^{n+1}} \cdot (1-(n+1)\tau)\omega_n \big(\frac{\lambda}{2\sqrt{n+1}}\big)^n,
\end{align*}
which implies $r_{q_{n+1}}\le  2^{10n}n^n\omega_{n+1}\lambda.$ So, if we take $\lambda=C(n)\beta^{-1} r_{q_{n+1}}$, then we get
$$1\le 2^{10n}n^n\omega_{n+1}C(n)\beta^{-1}.$$ This is the contradiction if we take $\beta> C(n)$ is big enough.
\end{proof}

In the rest of this subsection,  we always assume that
\begin{align}
\sup_{x\in \check{\mathscr{P}}_n(B_{r_{q_n}}(p))} d(p, x)\geq \beta^{-1}r_{q_n},\nonumber
\end{align}
we will prove Proposition \ref{prop (n+1) comp proj is dense} under this assumption. Let $\epsilon=\beta^{-\nu^{n+1}}$ in the rest argument, unless otherwise mentioned.

To prove Proposition \ref{prop (n+1) comp proj is dense}, we need set up $(n+ 1)$ projections and the corresponding net $\mathfrak{N}_{n+ 1}$, which will be used to obtain the corresponding almost Gou-Gu formula. All these results (Lemma \ref{lem one more point implies one more net} and Corollary \ref{cor one more gougu}) are obtained by the similar argument of Proposition \ref{prop dist of points in geodesic balls 1} and Theorem \ref{thm dist of points in geodesic balls induc-dim}.

\begin{lemma}\label{lem one more point implies one more net}
Assume there is a point $q_{n+1}\in \check{\mathscr{P}}_n(B_{r_{q_n}}(p))$ such that
$$d(p,q_{n+1})=\alpha r_{q_n}$$
for some $\alpha\ge \frac{1}{\beta}$. Then for $r_{q_{n+1}}:=\beta^{-2}d(p,q_{n+1})$, there exists an $\epsilon r_{q_{n+1}}$-dense net $\mathfrak{N}_{n+1}$ of $B_{r_{q_{n+1}}}(p)$ such that
\begin{enumerate}
\item[(a)]. For any $x\in \mathfrak{N}_{n+1}$ and $1\le l\le n+1$, $\pi_l(x)$ is well-defined;
\item[(b)]. For any $1\leq l\le i\le n+1$, $y\in \mathfrak{N}_{i}, z\in \mathfrak{N}_{n+1}$, we have
\begin{align*}
\big|d(y, z)^2-|\rho_l(y),\rho_l(z)|^2-d(\pi_l(y),\pi_l(z))^2 \big|\le \tilde{\epsilon}_{l}(r_{q_l})^2\le C(n)\beta^{-2}r_{q_{n+1}}^2,
\end{align*}
where $\tilde{\epsilon}_l=C(n)\beta^{-\nu^{n+1-l}}$ for $l\le n$ and $\tilde{\epsilon}_{n+1}=C(n)\beta^{-2}$;
 \item[(c)]. For any $1\le l\le {n+1}$, $\mathfrak{N}_{n+1}\subset T^{r_{q_l}}_{\eta, l}:=T_\eta^{r_{q_l}}(|\nabla^2f_l-g|)$ with respect to the function $\rho_l(\cdot)=d(\cdot, q_l)$.
\end{enumerate}
\end{lemma}

\begin{proof}
The proof is similar to the proof of Proposition \ref{prop dist of points in geodesic balls 1}.
\end{proof}

With the net $\mathfrak{N}_{n+1}$ constructed above, we define the projection $\mathscr{P}_0^{(n)}: B_{r_{q_{n+1}}}(p)\to \mathfrak{N}_{n+1}$ such that
$$d(\mathscr{P}_0^{(n)}(x),x)\le C(n)\beta^{-\nu^{n+1}}r_{q_{n+1}}, \quad  \forall x\in B_{r_{q_{n+1}}}(p).$$
Furthermore we define $$\hat{\pi}_n=\mathscr{P}_0^{n}\circ \pi_n,\quad \quad \mathscr{P}_n=\hat{\pi}_n\circ \ldots \circ \hat{\pi}_1\circ \hat{\pi}_0=\hat{\pi}_n\circ \mathscr{P}_{n-1},$$ $$\pi_{n+1}(x)=\overline{q_{n+1}x}\cap \partial B_{d(p,q_{n+1})}(q_{n+1})\text{ for } x\in B_{r_{q_{n+1}}}(p),$$
$$\check{\mathscr{P}}_{n+1}=\mathscr{P}_0^n\circ \pi_{n+1}\circ \mathscr{P}_n.$$

 Then we have the following corollary.
 \begin{cor}\label{cor one more gougu}
 Assume $\tau \in (0, c(n))$ and $1\ll \beta_0(n)\le \beta \le \tau^{-c(n)}$. Assume $\frac{V(B_r(p))}{V(B_r(0))}\ge 1-\tau$ and $\{q_k\}_{k=1}^n$ are constructed as in Theorem \ref{thm dist of points in geodesic balls induc-dim} and there exists a point $q_{n+1}\in \check{\mathscr{P}}_{n}(B_{r_{q_n}}(p))$ such that $\displaystyle d(q_{n+1},p)=\alpha r_{q_n}\ge \frac{1}{\beta} r_{q_n}$. Then, for the net $\mathfrak{N}_{n+1}$ constructed in Lemma \ref{lem one more point implies one more net} and $\check{\mathscr{P}}_{n+1}$ defined above, we have
 \begin{align*}
 \sup_{x,y\in B_{r_{q_{n+1}}}(p)}\big| d(x,y)^2&-\sum_{j=1}^{n+1}[\rho_j(\mathscr{P}_{j-1}(x))-\rho_j(\mathscr{P}_{j-1}(y))]^2\\
 &-d(\check{\mathscr{P}}_{n+1}(x),\check{\mathscr{P}}_{n+1}(y))^2\big|\le C(n)\beta^{-2}r_{q_{n+1}}^2.
 \end{align*}
 \end{cor}

 \begin{proof}
 For any $x,y\in B_{r_{q_{n+1}}}(p)$,  by Theorem \ref{thm dist of points in geodesic balls induc-dim}, we know
 \begin{align*}
 \big| d(x,y)^2&-\sum_{j=1}^{n}[\rho_j(\mathscr{P}_{j-1}(x))-\rho_j(\mathscr{P}_{j-1}(y))]^2\\
 &-d(\check{\mathscr{P}}_{n}(x),\check{\mathscr{P}}_{n}(y))^2\big|\le C(n)\beta^{-\nu}r_{q_{n}}^2\le C(n)\beta^{-2}r_{q_n}^2.
 \end{align*}
 By the definition of $\check{\mathscr{P}}_n, \mathscr{P}_0^{(n-1)}, \mathscr{P}_n$ and  $\mathscr{P}_0^{(n)}$, we know
 $$\big|d(\check{\mathscr{P}}_n(x),\check{\mathscr{P}}_n(y))^2-d(\mathscr{P}_n(x),\mathscr{P}_n(y))^2\big|\le C(n)\beta^{-\nu^{n+1}}(r_{q_n}^2+r_{q_{n+1}}^2)\le C(n)\beta^{-2}r_{q_{n+1}}^2.$$
 Note $\mathscr{P}_n(x),\mathscr{P}_n(y)\in \mathfrak{N}_{n+1}$. By Lemma \ref{lem one more point implies one more net}, we know
 \begin{align*}
 \big|d(\mathscr{P}(x),\mathscr{P}(y))^2&-|\rho_{n+1}(\mathscr{P}_n(x))-\rho_{n+1}(\mathscr{P}_n(y))|^2\\
 &-d(\pi_{n+1}\mathscr{P}_n(x),\pi_{n+1}\mathscr{P}_n(y))^2\big|\le C(n)\beta^{-2}r_{q_{n+1}}^2.
 \end{align*}
 Again by the definition of $\mathscr{P}_0^{(n)}$ and $\check{\mathscr{P}}_n$, we know
 $$\big|d(\pi_{n+1}\mathscr{P}_n(x),\pi_{n+1}\mathscr{P}_n(y))^2-d(\check{\mathscr{P}}_n(x),\check{\mathscr{P}}_n(y))^2\big|\le C(n)\beta^{-\nu^{n+1}}r_{q_{n+1}}^2.$$
 Combining the four inequalities above together, we get the conclusion.
 \end{proof}

Let $\gamma_{n+ 1}= \beta^{-1}r_{q_{n+ 1}}$ in the rest argument.

\begin{lemma}\label{lem crucial dist est-(n+1)}
{For $\hat{w}\in \mathfrak{N}_k$ where $j\leq k= n+ 1$,  if $\displaystyle \sup_{1\leq i\leq j- 1}d(\hat{w}, \pi_i(\hat{w}))\leq \delta$, then
\begin{align}
&\sup_{1\leq i\leq j- 1}d(\mathscr{P}_i(\hat{w}), \hat{w})\leq C(n)(\gamma_k+ \delta), \nonumber \\
&\sup_{1\leq i\leq n}d(q_{n+ 1}, \pi_i(q_{n+ 1}))+ \sup_{1\leq i\leq n- 1}d(\mathscr{P}_i(q_{n+ 1}), q_{n+ 1})\leq C(n)\gamma_n . \nonumber
\end{align}
}
\end{lemma}

\pf
{For the first conclusion, by Lemma \ref{lem proj and the composition of proj},  we only need to prove the case of $j=n+1$. Also by Lemma \ref{lem proj and the composition of proj}, we know in this case
 \begin{align}\label{an inequality needed}
 \sup_{1\le i\le n-1} d(\hat{w},\mathscr{P}_{i}(w))\le C(n)(\delta+\gamma_n).
 \end{align}
Now, by Lemma \ref{lem one more point implies one more net} and Corollary \ref{cor one more gougu}, we know
\begin{align*}
d(\hat{w},\mathscr{P}_n(\hat{w}))
&\le \delta +\beta^{-\nu^{n+1}r_{q_{n+1}}}+d(\pi_n\mathscr{P}_{n-1}(\hat{w}),\pi_n(\hat{w}))\\
&\le \delta +\gamma_{n+1}+\sqrt{2d(\hat{w},\mathscr{P}_{n-1}(\hat{w}))^2+C(n)\gamma_{n+1}^2}\\
&\le C(n)(\delta +\gamma_n+\gamma_{n+1}).
\end{align*}
By the choosing of $\gamma_{n+1}$, we know $\gamma_{n+1}\ge \gamma_n$. Thus the first half conclusion follows.  For  the second conclusion, we follow the proof of Lemma \ref{lem crucial dist est between point and proj of pt}. Denote $B_{r_{q_{n}}}(p)=B$. Since $q_{n+1}\in \check{\mathscr{P}}_n(B)$, there exists a $\check{q}_{n+1}\in \pi_n\mathscr{P}_{n-1}(B)$ such that $\displaystyle d(\check{q}_{n+1},q_{n+1})\le \beta^{-\nu^{n+1}}r_{q_n}\le \gamma_n$.

Thus we know
$$d(q_{n+1},\pi_{n}(q_{n+1}))=|d(q_{n+1},q_n)-d(\pi_n(q_{n+1}),q_n)|=|d(q_{n+1},q_n)-d(\check{q}_{n+1},q_n)|\le \gamma_n.$$
Now, choose $w_n\in \mathscr{P}_{n-1}(B)$ such that $\check{q}_{n+1}=\pi_n(w_n)$. Then by the definition of $\mathscr{P}_{n-1}$ and $\check{\mathscr{P}}(B)$ such that $d(w_n,\hat{w}_n)\le C(n)\gamma_{n-1}$. Note $q_n\in \check{\mathscr{P}}_{n-1}(B_{r_{q_{n-1}}}(p))$. By Lemma \ref{lem crucial dist est between point and proj of pt}, we know
\begin{align}
d(\hat{w}_n,\pi_i(\hat{w}_n))+d(q_n,\pi_i(q_n))\le C(n)\gamma_{n-1}, \quad \quad \quad 1\le i\le n-1. \nonumber
\end{align}

 Since either $\overline{q_n,w_n,\check{q}_{n+1}}$ or $\overline{q_n,\check{q}_{n+1},w_n}$ is a segment. By Lemma \ref{lem crucial trig ineq involved} or Lemma \ref{lem crucial trig ineq two}, we get
 $$d(\pi_i(q_{n+1}),q_{n+1})\le C(n)\big(\frac{\gamma_{n-1}}{r_{q_{n-1}}}\big)^{\frac{1}{4}}r_{q_{n-1}}\le C(n)\gamma_n.$$

Combining the above we know $\displaystyle \sup_{1\leq i\leq n}d(q_{n+ 1}, \pi_i(q_{n+ 1}))\leq C(n)\gamma_n$. Thus by (\ref{an inequality needed}) we get $\displaystyle \sup_{1\leq i\leq n- 1}d(\mathscr{P}_i(q_{n+ 1}), q_{n+ 1})\leq C(n)\gamma_n$. The conclusion follows.
}
\qed

Define
\begin{align}
\mathscr{A}_{n+ 1}= \Big(\mathfrak{U}_{\epsilon r_{q_{n+ 1}}} (\overline{p, q_{n+ 1}})\Big)\cap \mathfrak{N}_{n+ 1}, \quad \quad \quad \mathscr{A}_i= \bigcup_{y\in \mathscr{A}_{i+ 1}}\Big(\big(\mathfrak{U}_{\epsilon r_{q_{n+ 1}}}(\overline{q_i, y})\big)\cap \mathfrak{N}_{n+ 1}\Big) , \nonumber
\end{align}

\begin{lemma}\label{lem proj of points in Ai}
{We have $\displaystyle \sup_{1\leq k\leq i\leq n+ 1\atop z\in \mathscr{A}_i} d(\mathscr{P}_{k- 1}(z), z)+ d(\pi_{k- 1}(z), z)\leq C(n)\epsilon_n^{4^{i-n-1}}r_{q_{n}}$, where $\displaystyle \epsilon_n= \frac{\gamma_n}{r_{q_n}}$.
}
\end{lemma}

\pf
{\textbf{Step (1)}.  We firstly show the conclusion for $i= n+ 1$. For $z\in \mathscr{A}_{n+ 1}$ and $k\leq n+ 1$, note there is $\hat{p}\in \mathfrak{N}_{n}$ with
\begin{align}
d(\pi_{k- 1}(\hat{p}), \hat{p})+ d(\hat{p}, p)\leq \epsilon r_{q_n}; \nonumber
\end{align}
and also note $\displaystyle d(\pi_{k- 1}(q_{n+ 1}), q_{n+ 1})\leq C(n)\gamma_n$ by Lemma \ref{lem crucial dist est-(n+1)}. Applying Lemma \ref{lem crucial trig ineq involved} on $p, q_{n+ 1}$,  note $z\in \mathfrak{U}_{\epsilon r_{q_{n+ 1}}}(\overline{pq_{n+ 1}})\cap \mathfrak{N}_{n+ 1}$, we obtain
\begin{align}
d(\pi_{k- 1}(z), z)\leq C(n)\epsilon_n^{\frac{1}{4}}r_{q_n}. \label{proj est needed now}
\end{align}

From Lemma \ref{lem crucial dist est-(n+1)} and (\ref{proj est needed now}), for $z\in \mathscr{A}_{n+ 1}$, we have
\begin{align}
d(\mathscr{P}_{k- 1}(z), z)\leq C(n)(\gamma_n+ \epsilon_n^{\frac{1}{4}}r_{q_{n}})\leq C(n)\epsilon_n^{\frac{1}{4}}r_{q_n}. \nonumber
\end{align}
Hence the conclusion holds for $i= n+ 1$.

\textbf{Step (2)}.  Assume the conclusion holds for $i$,  hence
\begin{align}
 \sup_{1\leq k\leq i\atop z\in \mathscr{A}_i} d(\mathscr{P}_{k- 1}(z), z)+ d(\pi_{k- 1}(z), z)\leq C(n)\epsilon_n^{4^{i-n-1}}r_{q_n}. \nonumber
\end{align}

We will show the conclusion holds for $i- 1$.  For any $z\in \mathscr{A}_{i- 1}, k\leq i- 1$, we note $z\in \mathfrak{U}_{\epsilon r_{q_{n+ 1}}}(\overline{y q_{i- 1}})$ where $y\in \mathscr{A}_i$.  From the induction assumption,  we have
\begin{align}
d(\pi_{k- 1}(y), y)\leq C(n)\epsilon_n^{4^{i-n-1}}r_{q_n}. \nonumber
\end{align}

Note $\displaystyle d(\pi_{k- 1}(q_{i- 1}), q_{i- 1})\leq C(n)\gamma_n$ by Corollary \ref{cor comp of proj of qi},  from Lemma \ref{lem crucial trig ineq involved}, we obtain
\begin{align}
d(\pi_{k- 1}(z), z)\leq C(n)\epsilon_n^{4^{i-n-2}}r_{q_n}. \nonumber
\end{align}
From the above and Lemma \ref{lem crucial dist est-(n+1)},  we also get
$\displaystyle d(\mathscr{P}_{k- 1}(z), z)\leq C(n)\epsilon_n^{4^{i-n-2}}r_{q_n}$. The conclusion follows by the induction method.
}
\qed

\begin{lemma}\label{lem comp proj of z and comp proj of y}
{For $y\in \mathscr{A}_{k+ 1}, \hat{z}\in \Big(\mathfrak{U}_{\epsilon r_{q_{n+ 1}}}(\overline{y, q_k})\cap \mathfrak{N}_{n+ 1}\Big)$, we have
\begin{align}
\sup_{k\leq i\leq n} d(\mathscr{P}_i(y), \mathscr{P}_i(\hat{z}))\leq C(n)\epsilon_n^{4^{-n}}r_{q_n} . \nonumber
\end{align}
}
\end{lemma}

\pf
{When $i= k$, since $\hat{z}\in \Big(\mathfrak{U}_{\epsilon r_{q_{n+ 1}}}(\overline{y, q_k})$, there exists $z
\in \overline{y,q_k}$ such that $d(z,\hat{z})\le \epsilon r_{q_{n+1}}$. This implies
\begin{align*}
\big||d(y,q_k)-d(q_k,\hat{z})|^2-d(y,\hat{z})^2\big|
&\le \big||d(y,q_k)-d(q_k,\hat{z})|-d(y,\hat{z})\big|\big(|d(y,q_k)-d(q_k,\hat{z})|+d(y,\hat{z})\big)\\
&\le 2\epsilon r_{q_{n+1}}(2\epsilon r_{q_{n+1}}+2d(y,\hat{z}))\le C(n)\epsilon r_{q_{n+1}}^2.
\end{align*}
 Thus by Lemma \ref{lem proj of points in Ai}, we have
\begin{align}
d(\mathscr{P}_k(y), \mathscr{P}_k(\hat{z}))&\leq d(y,  \mathscr{P}_k(\hat{z}))+ C(n)\epsilon_n^{2^{-n}}r_{q_1} \nonumber \\
&\leq d(\pi_k\mathscr{P}_{k- 1}(\hat{z}), \pi_k(\hat{z}))+ d(\pi_k(\hat{z}), y)+ C(n)\epsilon_n^{4^{-n}}r_{q_n} \nonumber \\
&\leq d(\mathscr{P}_{k-1}(\hat{z}), \hat{z})+ d(\pi_k(\hat{z}),  \pi_k(y))+ C(n)\epsilon_n^{4^{-n}}r_{q_n}\nonumber \\
&\leq \sqrt{d(\hat{z}, y)^2- |d(\hat{z}, q_k)- d(y, q_k)|^2+C(n)\gamma_k^2}+ C(n)\epsilon_n^{4^{-n}}r_{q_n} \nonumber \\
&\leq C(n)\epsilon_n^{4^{-n}}r_{q_n} .\nonumber
\end{align}

Assume the conclusion holds for some $i\geq k$, we show the conclusion holds for $i+ 1$ as follows:
\begin{align}
d(\mathscr{P}_{i+ 1}(y), \mathscr{P}_{i+ 1}(\hat{z}))&\leq d(\pi_{i+ 1}\mathscr{P}_i(y), \pi_{i+ 1}\mathscr{P}_i(\hat{z}))+ C(n)\epsilon_n^{4^{-n}}r_{q_n}\nonumber \\
&\leq d(\mathscr{P}_i(y),  \mathscr{P}_i(\hat{z}))+ C(n)\epsilon_n^{4^{-n}}r_{q_n}\leq C(n)\epsilon_n^{4^{-n}}r_{q_n}. \nonumber
\end{align}
The conclusion follows by the induction method.
}
\qed

\begin{prop}\label{prop (n+1) comp proj is dense}
The set $\phi^{(n+ 1)}(B_{r_{q_{n+ 1}}}(p))$ is $C(n)\beta^{-1}r_{q_{n+ 1}}$-dense in $[-\frac{r_{q_{n+ 1}}}{2},0]^{n+1}\subset \mathbb{R}^{n+1}$.
\end{prop}

\pf
{Define $\psi^{(k)}= (\phi_k, \cdots, \phi_{n+ 1})$, we claim that $\psi^{(k)}(\mathscr{A}_k)$ is $C(n)\epsilon_n^{4^{-n}}r_{q_n}$-dense in $I^{n-k+ 2}= [-2^{-1}r_{q_{n+ 1}}, 0]^{n- k+ 2}$. Note $C(n)\beta^{-1}r_{q_{n+ 1}}\geq C(n)\epsilon_n^{4^{-n}}r_{q_n}$ by $\nu\geq 4^{n+ 1}$, then the conclusion follows from the case $k= 1$.

When $k= n+ 1$, for any $t_{n+ 1}\in I^1$, choose $\check{y}\in \overline{p q_{n+ 1}}$ such that $d(p, \check{y})= |t_{n+ 1}|$.  Choose $\hat{y}\in \mathscr{A}_{n+ 1}$ such that $d(y, \hat{y})\leq \epsilon r_{q_{n+ 1}}$,  from Lemma \ref{lem proj of points in Ai},  we get
\begin{align}
|\phi_{n+ 1}(\hat{y})- t_{n+ 1}|&= |d(\mathscr{P}_{n}(\hat{y}), q_{n+ 1})- d(p, q_{n+ 1})- t_{n+ 1}| \nonumber \\
&\leq |d(\check{y}, q_{n+ 1})- d(p, q_{n+ 1})- t_{n+ 1}|+ d(\check{y}, \mathscr{P}_{n}(\hat{y})) \nonumber \\
&\le d(\hat{y}, \mathscr{P}_{n}(\hat{y}))+ d(\hat{y}, \check{y}) \leq  C(n)\epsilon_n^{4^{-1}}r_{q_n}. \nonumber
\end{align}
The claim holds for $k= n+ 1$. We will prove the claim by induction on $k$.

Assume the claim holds for $k+ 1$, we prove the claim for $k$. Now for $(t_k, t_{k+ 1}, \cdots, t_{n+ 1})\in I^{n-k+ 2}$, from the induction assumption we can find $y\in \mathscr{A}_{k+ 1}$ such that
\begin{align}
|\psi^{(k+ 1)}(y)- (t_{k+ 1}, \cdots, t_{n+ 1})|\leq C(n)\epsilon_n^{4^{-n}}r_{q_n} . \nonumber
\end{align}

Choose $\check{z}\in \overline{q_k y}$ such that $d(\check{z}, y)= |t_k|$, choose $\hat{z}\in \mathfrak{N}_{n+ 1}$ such that $d(\hat{z}, \check{z})\leq \epsilon r_{q_{n+ 1}}$, then $\hat{z}\in \mathscr{A}_k$.  Now from Lemma \ref{lem comp proj of z and comp proj of y} we have
\begin{align}
&|\phi_k(\hat{z})- t_k|= |d(\mathscr{P}_{k- 1}(\hat{z}), q_k)- d(p, q_k)- t_k|\nonumber \\
&\leq |d(\check{z}, q_k)- d(p, q_k)- t_k|+ d(\check{z}, \mathscr{P}_{k- 1}(\hat{z})) \nonumber \\
&\leq d(\hat{z}, \mathscr{P}_{k- 1}(\hat{z}))+ d(\check{z}, \hat{z})\leq \epsilon r_{q_k}+ C(n)\epsilon_n^{4^{-n}}r_{q_n}\leq C(n)\epsilon_n^{4^{-n}}r_{q_n}  . \nonumber
\end{align}

For $i\geq k+ 1$,  from Lemma \ref{lem comp proj of z and comp proj of y} we have
\begin{align}
|\phi_i(\hat{z})- t_i|&= |d(\mathscr{P}_{i- 1}(\hat{z}), q_i)- d(p, q_i)- t_i|\nonumber \\
&\leq |d(\mathscr{P}_{i- 1}(y), q_i)- d(p, q_i)- t_i|+ d(\mathscr{P}_{i- 1}(y), \mathscr{P}_{i- 1}(\hat{z})) \nonumber \\
&\leq |\phi_i(y)- t_i|+ C(n)\epsilon_n^{4^{-n}}r_{q_n} \leq C(n)\epsilon_n^{4^{-n}}r_{q_n} .  \nonumber
\end{align}
Note $\phi^{(n+1)}=\psi^{(1)}$. The conclusion follows by the induction method.
}
\qed

\subsection{The quasi-isometry property of distance map}

For the orthogonal basis $\{e_i\}_{i=1}^n$ of $\mathbb{R}^n$, let $q_i=\beta r e_i$, $b_i(x)=d(q_i,x)-d(q_i,0)$ and $\psi(x)=(b_1(x),b_2(x),\ldots, b_n(x))$. Then we have the bi-Lipschitz estimate
$$1-\frac{4\sqrt{n}}{\sqrt{\beta}}\le \frac{|\psi(x)-\psi(y)|}{|x-y|}\le 1+\frac{5\sqrt{n}}{\sqrt{\beta}}, \quad  \forall  x,y\in B_r(0).$$

The following lemma is the almost version of the above result on manifolds, which shows that the distance map is a discrete bi-Lipschitz map with Lipschitz constant close to $1$.
\begin{theorem}\label{thm direction points imply Lipschitz less than 1+ep}
{Let $\displaystyle \psi(x)= (b_1^+(x), \cdots, b_n^+(x))$ and $b_i^+(x)= d(x, q_i)-d(p,q_i)$, then
\begin{align}
1- C(n)\beta ^{-\frac{1}{2}} \le \sup_{x, y\in B_{r_{q_n}}(p), \atop d(x, y)\ge \beta^{-\frac{1}{2}} r_{q_n}}\frac{d(\psi(x), \psi(y))}{d(x, y)}\leq 1+ C(n)\beta ^{-\frac{1}{2}}  . \nonumber
\end{align}
}
\end{theorem}

\begin{remark}\label{rem dist map is G-H map}
{In later argument, we use this distance map to get a `canonical' $\epsilon$-G-H approximation in analytical sense, which is the $\epsilon$-splitting map.
}
\end{remark}

\pf
{Let $r_1= r_{q_n}$, note
\begin{align}
d(\psi(x), \psi(y))^2= \sum_{i= 1}^n |b_i^+(x)- b_i^+(y)|^2= \sum_{i= 1}^n \Big|\frac{d(x, q_i)^2- d(y, q_i)^2}{d(x, q_i)+ d(y, q_i)}\Big|^2. \label{Lipschitz map}
\end{align}

For $i\geq 1$ and any $x\in B_{r_1}(p)$, apply Theorem \ref{thm dist of points in geodesic balls induc-dim} on $B_{r_1}(p)$, we get
\begin{align}
\Big|d(x, q_i)^2- \sum_{j= 1}^{i- 1}|\phi_j(x)- \phi_j(q_i)|^2- d\big(\check{\mathscr{P}}_{i-1}(x), \check{\mathscr{P}}_{i-1}(q_i)\big)^2\Big|\leq C(n)\gamma_{n}^2. \nonumber
\end{align}

Now for any $x, y\in B_{r_1}(p)$, from Proposition \ref{prop n direction points} we get
\begin{align}
&\Big||d(x, q_i)^2- d(y, q_i)^2|- |d\big(\check{\mathscr{P}}_{i-1}(x), \check{\mathscr{P}}_{i-1}(q_i)\big)^2- d\big(\check{\mathscr{P}}_{i-1}(y), \check{\mathscr{P}}_{i-1}(q_i)\big)^2| \Big|\nonumber\\
&\leq C(n)\gamma_{n}^2+ \sum_{j= 1}^{i- 1}|\phi_j(x)- \phi_j(y)|\cdot |\phi_j(x)+ \phi_j(y)- 2\phi_j(q_i)| \leq C(n)r_{q_n}^2\label{upper bound of num}
\end{align}
By Corollary \ref{cor comp of proj of qi} and the definition of $\check{\mathscr{P}}_{i-1}$, we know
$$d(\check{\mathscr{P}}_{i-1}(q_i),q_i)+d(\check{\mathscr{P}}_{i-1}(x),\mathscr{P}_{i-1}(x))\le C(n)\gamma_i,$$
which implies
\begin{align}\label{check to approx}
\big|d\big(\check{\mathscr{P}}_{i-1}(x), \check{\mathscr{P}}_{i-1}(q_i)\big)^2-d(\mathscr{P}_{i-1}(x),q_i)^2\big|\le C(n)\epsilon_n r^2_{q_n}+C(n)\gamma_i d(p,q_i).
\end{align}

From the assumption, for any $x\in B_{r_1}(p)$, we have $d(x, q_i)\geq c(n)\beta^{\nu^{n+ 1- i}} r_{q_i}.$ Thus $\frac{\gamma_id(p,q_i)}{d(x,q_i)+d(y,q_i)}\le C(n) \gamma_i$ and
\begin{align}\label{almost lip}
\frac{|d(\mathscr{P}_{i-1}(x),q_i)^2- d(\mathscr{P}_{i-1}(y),q_i)^2|}{d(x,q_i)+d(y,q_i)}&=\frac{|\phi_i(x)-\phi_i(y)|\big(\rho_i(\mathscr{P}_{i-1}(x))+\rho_i(\mathscr{P}_{i-1}(y))\big)}{d(x,q_i)+d(y,q_i)}\nonumber\\
&\in \big((1-\beta^{-5}), (1+\beta^{-5})\big)|\phi_i(x)-\phi_i(y)|.
\end{align}

From (\ref{upper bound of num}), (\ref{check to approx}) and (\ref{almost lip}), we have
\begin{align}
\Big|\frac{d(x, q_i)^2- d(y, q_i)^2}{d(x, q_i)+ d(y, q_i)}\Big| &\le (1+ \beta^{-5}) |\phi_i(x)- \phi_i(y)|+ C(n)\beta^{-10}r_{q_n} , \nonumber \\
\Big|\frac{d(x, q_i)^2- d(y, q_i)^2}{d(x, q_i)+ d(y, q_i)}\Big| &\ge (1- \beta^{-5}) |\phi_i(x)- \phi_i(y)|- C(n)\beta^{-10}r_{q_n} . \nonumber
\end{align}
Plugging the above into (\ref{Lipschitz map}), note $|\phi_{i}(x)- \phi_{i}(y)|\leq C(n) r_{q_n}$, we get
\begin{align}
d(\psi(x), \psi(y))^2 &\leq  (1+ \beta^{-5})^2\sum_{i= 1}^n |\phi_i(x)- \phi_i(y)|^2+ C(n)\beta^{-10} r^2_{q_n} ,\label{upper bound of numerator-final} \\
d(\psi(x), \psi(y))^2 &\ge  (1- \beta^{-5})^2\sum_{i= 1}^n |\phi_i(x)- \phi_i(y)|^2-C(n)\beta^{-10} r^2_{q_n}.\label{lower bound of numerator-final}
\end{align}
From Theorem \ref{thm dist of points in geodesic balls induc-dim}, one has
\begin{align}
\Big|d(x, y)^2- \sum_{j= 1}^n|\phi_j(x)- \phi_j(y)|^2- d(\check{\mathscr{P}}_{n}(x), \check{\mathscr{P}}_{n}(y))^2\Big|\leq C(n)\gamma_n^2 . \nonumber
\end{align}
Now we get
\begin{align}
\sup_{x, y\in B_p(r_{q_n}), \atop d(x, y)\ge \beta^{-\frac{1}{2}} r_{q_n}}\frac{d(\psi(x), \psi(y))}{d(x, y)}
&\leq \sqrt{\frac{(1+ \beta^{-5})^2\sum\limits_{i= 1}^n |\phi_i(x)- \phi_i(y)|^2+ C(n)\beta^{-10}r^2_{q_n}}{d(x, y)^2}}  \nonumber \\
&\leq \sqrt{\frac{(1+ \beta^{-5})^2d(x, y)^2+ C(n) \beta^{-2}r_{q_n}^2}{d(x, y)^2}} \leq 1+ C(n)\beta^{-\frac{1}{2}}, \nonumber
\end{align}
where in the last line we use $d(x,y)\ge \beta^{-\frac{1}{2}}r_{q_n}.$

On the other hand, for $x, y\in B_{r_1}(p)$ by  Proposition \ref{prop the diam upper bound of (n+ 1)-proj}, we get
\begin{align}
|\sum_{j= 1}^n|\phi_j(x)- \phi_j(y)|^2- d(x, y)^2| \le C(n)\gamma_n^2+ d(\check{\mathscr{P}}_{n}(x),  \check{\mathscr{P}}_{n}(y))^2 \le C(n)\beta^{-2}r_{q_n}^2. \nonumber
\end{align}
Similarly we have
\begin{align}
\sup_{x, y\in B_p(r_{q_n}), \atop d(x, y)\ge \beta^{-\frac{1}{2}} r_{q_n}}\frac{d(\psi(x), \psi(y))}{d(x, y)}
&\ge \sqrt{\frac{(1-\beta^{-5})^2\sum\limits_{i= 1}^n |\phi_i(x)- \phi_i(y)|^2- C(n)\beta^{-2}r^2_{q_n}}{d(x, y)^2}}  \nonumber \\
&\ge 1-C(n) \beta^{-\frac{1}{2}}. \nonumber
\end{align}
}
\qed

\section{The existence of $\epsilon$-splitting maps}\label{sec existence of splitting map}

\subsection{The almost orthogonality of distance maps}\label{subsec almost og of dist map}

The integral Toponogov Comparison Theorem for $Rc\geq 0$ was proved in \cite{Colding-volume} by approximating the distance function by harmonic functions and covering argument. Motivated by the cosine law for triangle (the terms with distance functions are in the form of the square of distance functions), we consider the model function $f$ defined in Subsection \ref{subsec integral est of diff} relating with the square of the distance function, present a different proof of this result here (Corollary \ref{cor Colding-1.28}), which avoids the covering argument. And our argument is more consistent with the former argument for almost Gou-Gu Theorem of distance functions. Then we use the integral Toponogov comparison Theorem to establish the almost orthogonality of distance maps.

Let $SM^n$ be the unit tangent bundle of $M^n$, if $\pi: SM^n\rightarrow M^n$ is the projection map, for any $\Omega\subset SM^n$, the \textbf{Liouville measure} of $\Omega$, denoted by $\mu^L(\Omega)$, is defined by $\mu^L(\Omega)= \int_{\pi(\Omega)}\mathcal{H}^{n-1}(\pi^{-1}(x))d\mu(x)$, where $\mu$ is the volume measure of $(M^n, g)$ determined by the metric $g$.  We begin with an integral version of the law of cosine with respect to the model $\mathbb{R}^n$.

\begin{lemma} \label{lem almost cosine}
Assume $B_{11r}(q)$ is a geodesic ball in $(M^n,g)$ with $Rc\geq 0$ and $f$ is a smooth function. Then for any $s\in [0, 10r]$, we have
  \begin{align*}
    &\fint_{SB_r(q)}\big|\frac{d_q^2(\gamma_v(s))}{2}-\frac{d_q^2(x)}{2}-\langle \nabla \frac{d_q^2(x)}{2},v\rangle s-\frac{s^2}{2}\big|d\mu^L(x, v)\\
    &\le C(n)\Big\{\sup_{B_{11r}(q)}|f-\frac{d_q^2}{2}|+ r\big(\fint_{B_r(q)}|\nabla(f-\frac{d^2_q}{2})|^2\big)^{\frac{1}{2}}+ r^2\cdot \big(\fint_{B_{11 r}(q)}|\nabla^2 f-g|^2\big)^{\frac{1}{2}}\Big\},
  \end{align*}
  where $d\mu^L(x, v)$ is the Liouville measure on the sphere bundle $SM$.
\end{lemma}
\begin{proof}
  For any $v\in SB_r(q)$, we have
  $$(f\circ \gamma_v)'(t)-(f\circ \gamma_v)'(0)=\int_0^t\nabla^2f(\gamma_v'(\tau),\gamma_v'(\tau))d\tau.$$
  So, we have
  \begin{align*}
  \sup_{0\le t\le 10r} |(f\circ \gamma_v)'(t)-(f\circ \gamma_v)'(0)-t|\le \int_0^t |\nabla^2 f-g|_{\gamma_v(\tau)}d\tau.
  \end{align*}
  Integrating on $SB_r(q)$ with respect to the Louisville measure $d\mu^L(x, v)$, we get
  \begin{align*}
  &\int_{SB_r(q)}\sup_{0\le t\le 10r} |(f\circ \gamma_v)'(t)-\langle \nabla f(x),v\rangle -t|d\mu^L(x, v)\\
  &\le \int_0^{10r}\int_{SB_r(q)}|\nabla^2f-g|(\gamma_v(\tau))d\mu^L(x, v)d\tau\\
  &\le  \int_0^{10r}\int_{SB_{11r}(q)}|\nabla^2f-g|(y) dV(y,v)d\tau\\
  &=10n\omega_n r\int_{B_{11r}(q)}|\nabla^2f-g|(y)dV(y).
  \end{align*}
  Note that for $0\le s\le 10r$, we have
    \begin{align*}
  \big|(f\circ \gamma_v)(s)-(f\circ \gamma_v)(0)-\langle \nabla\frac{d_q^2}{2},v\rangle s-\frac{s^2}{2}\big|&=\big|\int_0^2(f\circ \gamma_v)'(t)-\langle \frac{d_q^2(x)}{2},v\rangle-t\big|dt\\
  &\le 10r\sup_{0\le t\le 10r}|(f\circ \gamma_v)'(t)-\langle \frac{d_q^2(x)}{2},v\rangle-t|,
  \end{align*}
  \begin{align*}
  \int_{SB_r(q)}|\langle \nabla (f-\frac{d_q^2}{2}),v\rangle|d\mu^L(x, v)\le n\omega_n(V(B_r(q)))^{\frac{1}{2}}\big(\int_{B_r(q)}|\nabla(f-\frac{d^2_q}{2})|^2dV(x)\big)^{\frac{1}{2}},
  \end{align*}
  \begin{align*}
  \int_{SB_r(q)}|f-\frac{d_q^2}{2}|d\mu^L(x, v)\le n\omega_nV(B_r(q))\sup_{B_r(q)}|f-\frac{d_q^2}{2}|
  \end{align*}
  and
  \begin{align*}
  \int_{SB_r(q)} |f\circ \gamma_v(s)-\frac{d_q^2(\gamma_v(s))}{2}|d\mu^L(x, v)\le n\omega_nV(B_r(q))\sup_{B_{11r}(q)}|f-\frac{d_q^2}{2}|.
  \end{align*}
  We know
  \begin{align*}
  &\int_{SB_r(q)} \big|\frac{d_q^2(\gamma_v(s))}{2}-\frac{d_q^2(x)}{2}-\langle \nabla \frac{d_q^2(x)}{2},v\rangle s-\frac{s^2}{2}\big|  \\
    &\le  2n\omega_nV(B_r(q))\sup_{B_{11r}(q)}|f-\frac{d_q^2}{2}|
    +10r n\omega_nV(B_r(q))^{\frac{1}{2}}\big(\int_{B_r(q)}|\nabla(f-\frac{d^2_q}{2})|^2dV(x)\big)^{\frac{1}{2}}\\
    & \quad \quad \quad \quad \quad \quad \quad \quad \quad \quad \quad \quad \quad \quad + 100r^2n\omega_n\cdot\big(\int_{B_{11 r}(q)}|\nabla^2 f-g|dV(x)\big)^{\frac{1}{2}}.
  \end{align*}
  Dividing both sides by $n\omega_n V(B_r(q))$ and using  H\"older inequality and Bishop-Gromov volume comparison theorem, we get the conclusion.
\end{proof}
The integral law of cosine implies the Busemann function is almost linear. More precisely, the directional derivative is close to the difference quotient in an integral sense.
\begin{lemma}\label{lem almost linear}
  Assume $B_{11r}(q)\subseteq (M^n,g)$ with $Rc\geq 0$ and $f$ is a smooth function.  Then for any $B_{r_1}(p)\subset B_r(q)$ satisfying $d_q(p)\ge 2 r_1$ and $s\le 10 r$, there holds
\begin{align*}
  &\fint_{SB_{r_1}(p)}\big|\langle \nabla b^+(x),v\rangle -\frac{b^+(\gamma_v(s))-b^+(x)}{s}\big|\\
  &\le  \frac{2s}{d_q(p)}+ \frac{C(n)}{sd_q(p)}\cdot (\frac{r}{r_1})^n\big\{ \sup_{B_{11r}(q)}|f-\frac{d_q^2}{2}|+r\big(\fint_{B_r(q)}|\nabla(f-\frac{d^2_q}{2})|^2dV(x)\big)^{\frac{1}{2}}\\
  &\quad \quad\quad\quad\quad\quad\quad\quad\quad\quad\quad\quad \quad\quad\quad\quad+ r^2\cdot\big(\fint_{B_{11 r}(q)}|\nabla^2 f-g|^2dV(x)\big)^{\frac{1}{2}}
  \big\}
  \end{align*}
where $b^+(x)=d(q,x)-d(q,p)$.
\end{lemma}
\begin{proof}
  By the definition of $b^+$, we know
  \begin{align*}
  &\frac{1}{2sd_q(x)}\big| d_q^2(\gamma_v(s))-d_q^2(x)-\langle d_q^2(x),v\rangle s-s^2\big|\\
  &=\big|(1+\frac{b^+(\gamma_v(s))-b^+(x)}{2d_q(x)})(\frac{b^+(\gamma_v(s))-b^+(x)}{s})-\langle \nabla b^+(x),v\rangle -\frac{s}{2d_q(x)}\big|.
  \end{align*}
  Note that $|b^+(\gamma_v(s))-b^+(x)|\le s$ and $d_q(x)\ge \frac{d_q(p)}{2}$ for $x\in B_{r_1}(p)$ and $2r_1\le d_q(p)$. We have
  \begin{align*}
  &\big|\langle \nabla b^+(x),v\rangle -\frac{b^+(\gamma_v(s))-b^+(x)}{s}\big|\\
  &\le \frac{1}{2sd_q(x)}\big| d_q^2(\gamma_v(s))-d_q^2(x)-\langle d_q^2(x),v\rangle s-s^2\big|+\frac{s}{d_q(x)}\\
  &\le \frac{1}{sd_q(p)}\big| d_q^2(\gamma_v(s))-d_q^2(x)-\langle d_q^2(x),v\rangle s-s^2\big|+\frac{2s}{d_q(p)}
  \end{align*}
  So, by Lemma \ref{lem almost cosine}, we know
  \begin{align*}
  &\fint_{SB_{r_1}(p)}\big|\langle \nabla b^+(x),v\rangle -\frac{b^+(\gamma_v(s))-b^+(x)}{s}\big|\\
  &\le \frac{2s}{d_q(p)}+\frac{V(B_r(q))}{sd_q(p) V(B_{r_1}(p))}\fint_{SB_{r}(q)}\big| d_q^2(\gamma_v(s))-d_q^2(x)-\langle d_q^2(x),v\rangle s-s^2\big|d\mu^L(x, v)\\
  &\le  \frac{2s}{d_q(p)}+ C(n)\frac{r^n}{sd_q(p) r_1^n}\big\{ \sup_{B_{11r}(q)}|f-\frac{d_q^2}{2}|+r\big(\fint_{B_r(q)}|\nabla(f-\frac{d^2_q}{2})|^2dV(x)\big)^{\frac{1}{2}}\\
  &\quad \quad\quad\quad\quad\quad\quad\quad\quad\quad\quad\quad \quad\quad\quad\quad+ r^2\cdot\big(\fint_{B_{11 r}(q)}|\nabla^2 f-g|^2dV(x)\big)^{\frac{1}{2}}
  \big\}
  \end{align*}
\end{proof}

\begin{cor}\label{cor Colding-1.28}
  For $0<\tau\le \tau_0(n)\ll 1$, $0<\beta_0(n)\le \beta \le \tau^{-c(n)}$ and $r>0$.  Assume $B_r(p)\subseteq (M,g)$ with $Rc\ge 0$ satisfying $\displaystyle \frac{V(B_r(p))}{V(B_r(0))}\ge 1-\tau.$
Choose $\{q_i\}_{i=1}^{n}$ as in Theorem \ref{thm dist of points in geodesic balls induc-dim} and $b_i^+(x)=d(q_i,x)-d_(q_i,p)$. Then for $r_1=r_{q_n}=c(n)\beta^{-(C(n)+\frac{\nu^{n+1}-\nu}{\nu-1})}r$
 and $s\le r_1$, there holds
 \begin{align*}
  \fint_{SB_{r_1}(p)}\big|\langle \nabla b_i^+(x),v\rangle &-\frac{b_i^+(\gamma_v(s))-b_i^+(x)}{s}\big|d\mu^L(x, v)\\
  &\le  2\beta^{-\frac{\mu^{n+2-i}-\nu^2}{\nu-1}}\frac{s}{r_1}+C(n)\frac{r_1}{s}\beta^{-(C(n)-\frac{(n+1)(\nu^{n+2-i}-\nu^2)}{\nu-1})}.
  \end{align*}
 \end{cor}

 \begin{proof}
 By the choose of $q_i$,  we know $C(n)\beta^{\frac{v^{n+2-i}-v^2}{v-1}}r_1= d_{q_i}(p)\le \beta^{-C(n)}r$.  So, we know
 $$\frac{V(B_{200d_{q_i}(p)}(q_i))}{V(B_{200d_{q_i}(p)}(0))}\ge 1-C(n)\beta^{-C(n)}$$
  and $B_{r_1}(p)\subset B_{d_{q_i}(p)}(q_i)$ with $2r_1\le d_{q_i}(p)$. Thus by Proposition \ref{prop L2 diff of gadient of F}, \ref{prop C0 diff of F} and \ref{prop integral Hessian small}, we know there exists $f_i$ such that
  \begin{align*}
  \frac{1}{d^2_{q_i}(p)}\sup_{B_{11 d_{q_i}(p)}(q_i)}|f-d_{q_i}^2|&+\frac{1}{d_{q_i}(p)}\big(\int_{B_{d_{q_i}(p)}(q_i)}|\nabla(f-d_{q_i}^2(p))|^2\big)^{\frac{1}{2}}\\
  &+\big(\int_{B_{11d_{q_i}(p)}(q_i)}|\nabla^2f-2g|^2\big)^{\frac{1}{2}}
  \le C(n)\beta^{-\frac{C(n)}{2n+4}}.
  \end{align*}
  So, by Lemma \ref{lem almost linear}, we know
  \begin{align*}
  \fint_{SB_{r_1}(p)}\big|\langle \nabla b_i^+(x),v\rangle &-\frac{b_i^+(\gamma_v(s))-b_i^+(x)}{s}\big|\le  \frac{2s}{d_{q_i}(p)}+ \frac{C(n)}{sd_{q_i}(p)} (\frac{d_{q_i}(p)}{r_1})^n \beta^{-C(n)}d_{q_i}^2(p)\\
  &\le 2\beta^{-\frac{\nu^{n+2-i}-\nu^2}{\nu-1}}\frac{s}{r_1}+C(n)\frac{r_1}{s}\beta^{-(C(n)-\frac{(n+1)(\nu^{n+2-i}-\nu^2)}{\nu-1})}.
  \end{align*}
 \end{proof}

 Theorem \ref{thm direction points imply Lipschitz less than 1+ep} means the functions $\{b_i^{+}\}_{1\le i\le n}$ are almost orthogonal in the difference quotient sense. Combining it with the above corollary, we get the orthogonality of $\{b_i^{+}\}_{1\le i\le n}$ in the following sense.

\begin{prop}\label{prop approximation component is almost orthogonal}
{For $0<\tau\le \tau_0(n)\ll 1$, $0<\beta_0(n)\le \beta \le \tau^{-c(n)}$ and $r>0$.  Assume $B_r(p)$ is a geodesic ball in $(M,g)$ with $Rc\ge 0$ satisfying $\displaystyle \frac{V(B_r(p))}{V(B_r(0))}\ge 1-\tau$. Then for  $\{q_i\}_{i=1}^{n}$ as in Theorem \ref{thm dist of points in geodesic balls induc-dim},  $r_1=r_{q_n}=c(n)\beta^{-(C(n)+\frac{\nu^{n+1}-\nu}{\nu-1})}r$  and $b_i^+(x)=d(q_i,x)-d_(q_i,p)$, there holds
\begin{align}
\fint_{B_{r_1}(p)} \big|\langle \nabla b_i^+, \nabla b_j^+\rangle\big|^2\leq C(n)(\beta^{-\frac{1}{2}}+\tau^{\frac{2}{n+1}}) , \quad \quad \quad \quad 1\leq i\neq j\leq n .\nonumber
\end{align}
}
\end{prop}

\pf
{We only need to show the conclusion for fixed $i, j$ with $i\neq j$. For $\theta\in (0, \frac{\pi}{2})$ (to be determined later), set $\displaystyle C_{\theta}(x)\vcentcolon= \Big\{\nu\in S(T_xM)|\ \langle \nu, \nabla b_i^+(x)\rangle\geq  \cos\theta \Big\}$ and $\displaystyle C_\theta=\cup_{x\in B_{r_1}(p)}C_\theta(x)$.

Moreover, recall $l_v$ is the largest number such that $\exp_x(tv)$ is a segment for any $t\le l_v$.  $\mathcal{B}(x)=\{v\in S(T_xM)| l_v\le  r_1\}$, , and $\mathcal{B}_\theta=\cup_{x\in B_{r_1}(p)} \mathcal{B}(x)\cap C_{\theta}$.  Then for any $\nu\in C_{\theta}$,
\begin{align}
|\nu- \nabla b_i^+|^2= |\nu|^2+ |\nabla b_i^+|^2- 2\langle \nu, \nabla b_i^+\rangle \leq 2- 2\cos\theta. \label{need length difference}
\end{align}
Using  $\frac{V(B_{r}(p))}{\omega_n r^n}\ge 1-\tau $, by the same argument as the proof of Lemma \ref{lem dist between points and segment}, we know
$\frac{\mathcal{H}^{n-1}(\mathcal{B}(x))}{\mathcal{H}^{n-1}(S(T_xM))}\le 2\tau,$
which implies
\begin{align}\label{bad ratio}
\frac{V(\mathcal{B}_\theta)}{V(SB_{r_1}(p))}\le 2\tau.
\end{align}
 From Corollary \ref{cor Colding-1.28}, let $s=r_1$, then we have
\begin{align}
\frac{1}{V\big(SB_{r_1}(p)\big)} \int_{SB_{r_1}(p)} \Big|\langle \nabla b_k^+, \nu \rangle- f_s(\nu)\Big|\leq \epsilon_1 \ , \quad \quad \quad \quad k= 1, \cdots, n\label{need 2.11}
\end{align}
where $f_k(\nu)\vcentcolon= \frac{(b_k^+\circ \gamma_{\nu})(r_1)- (b_k^+\circ \gamma_{\nu})(0)}{r_1}$ for $\nu\in SB_{r_1}(p)$ and $\epsilon_1=C(n)\beta^{-C(n)}$. So,
\begin{align}
&\quad \frac{1}{V\big(SB_{r_1}(p)\big)}\int_{C_{\theta}} \Big|\big|f_i(\nu)\big|^2- 1\Big|  \nonumber \\
&\leq \frac{1}{V\big(SB_{r_1}(p)\big)}\int_{C_{\theta}} \Big|\big|f_i(\nu)\big|^2- \langle \nabla b_i^+, \nu \rangle^2\Big| +  \frac{1}{V\big(SB_{r_1}(p)\big)}\int_{C_{\theta}} \big|1- \langle \nabla b_i^+, \nu \rangle^2 \big|  \nonumber \\
&\leq 2\epsilon_1+ (1- \cos^2\theta)\cdot \frac{V(C_{\theta})}{V\big(SB_{r_1}(p)\big)}\leq 2\epsilon_1+ \theta^2 \frac{V(C_{\theta})}{V\big(SB_{r_1}(p)\big)} \label{need lem 4.4.1}
\end{align}

Now from Theorem \ref{thm direction points imply Lipschitz less than 1+ep}, Corollary \ref{cor Colding-1.28}, (\ref{bad ratio})and (\ref{need lem 4.4.1}), for $j\neq i$, we have
\begin{align}
&\quad \frac{1}{V\big(SB_{r_1}(p)\big)} \int_{C_{\theta}} \big|\langle \nabla b_j^+, \nu \rangle\big|^2\nonumber \\
&\leq \frac{2}{V\big(SB_{r_1}(p)\big)} \int_{SB_{r_1}(p)} \Big|\langle \nabla b_j^+, \nu \rangle- f_j(\nu)\Big|^2+ \frac{2}{V\big(SB_{r_1}(p)\big)} \int_{C_{\theta}} \big|f_j(\nu)\big|^2 \nonumber \\
& \leq 4\epsilon_1 + \frac{2}{V\big(SB_{r_1}(p)\big)}\Big(\int_{C_{\theta}} \Big|\frac{d\big(\psi\circ \gamma_{\nu}(r_1), \psi\circ \gamma_{\nu}(0)\big)}{r_1}\Big|^2- \big|f_i(\nu)\big|^2\Big) \nonumber \\
&\leq 4\epsilon_1+ \frac{2}{V\big(SB_{r_1}(p)\big)}\Big\{3\delta V(C_{\theta}\backslash \mathcal{B}_\theta)+nV(\mathcal{B}_\theta)+ \Big(\int_{C_{\theta}} \Big|1- \big|f_i(\nu)\big|^2\Big|\Big)\Big\} \nonumber \\
&\leq 8\epsilon_1 +4n\tau + (6\delta + 2\theta^2)\frac{V(C_{\theta})}{V\big(SB_{r_1}(p)\big)}, \nonumber
\end{align}
where $\delta =\beta^{-\frac{1}{2}}$ comes from Theorem \ref{thm direction points imply Lipschitz less than 1+ep}.
Note (\ref{need length difference}), then
\begin{align}
\int_{C_{\theta}} \big|\langle \nabla b_i^+, \nabla b_j^+\rangle\big|^2 &\leq 2\int_{C_{\theta}} \big|\langle \nabla b_{j}^+, v- \nabla b_i^+\rangle\big|^2+ 2\int_{C_{\theta}} \big|\langle \nabla b_j^+, \nu\rangle \big|^2 \nonumber \\
&\leq 2\int_{C_{\theta}} \big|v- \nabla b_i^+\big|^2+ 2\int_{C_{\theta}} \big|\langle \nabla b_j^+, \nu\rangle \big|^2 \nonumber \\
&\leq 2\theta^2 \cdot V(C_{\theta})+ 16n(\epsilon_1+ \tau)V\big(SB_{r_1}(p)\big)+ 12(\delta+\theta^2)V(C_{\theta}) .\nonumber
\end{align}

Now note $\langle \nabla b_i^+, \nabla b_j^+\rangle$ is constant on $T_xM^n$ for any fixed $x\in M^n$, we have
\begin{align}
\fint_{B_{r_1}(p)} \big|\langle \nabla b_i^+, \nabla b_j^+\rangle \big|^2&= \frac{1}{V\big(C_{\theta}\big)}\int_{C_{\theta}} \big|\langle \nabla b_i^+, \nabla b_j^+\rangle\big|^2 \nonumber \\
&\leq 14(\theta^2+\delta)+ 16n(\epsilon_1+ \tau)\cdot \frac{V\big(SB_{r_1}(p)\big)}{V(C_{\theta})} \nonumber
\end{align}
Note that $\frac{V(C_\theta)}{V(SB_{r_1}(p))}=\frac{\int_0^\theta \sin^{n-2}(\alpha)d\alpha}{\int_0^{\pi} \sin^{n-2}(\alpha)d\alpha}\ge c(n)\theta^{n-1}$. We can choose $\theta=c(n)(\epsilon_1+\tau)^{\frac{1}{n+1}}$ such that $\theta^2=\frac{\epsilon_1+\tau}{c(n)\theta^{n-1}}$. Then
$$\fint_{B_{r_1}(p)} \big|\langle \nabla b_i^+, \nabla b_j^+\rangle \big|^2\le C(n)(2\theta^2+\delta)\leq C(n)(\beta^{-\frac{1}{2}}+\tau^{\frac{2}{n+1}}).$$
}
\qed

\subsection{The $\epsilon$-splitting map induced by distance map}\label{subsec splitting map}

\begin{definition}\label{def excess est}
{For $q^+, q^-, p\in \mathbf{X}$, where $\mathbf{X}$ is a metric space, we say that $[q^+, q^-, p]$ is an \textbf{AG-triple on $\mathbf{X}$ with the excess $s$ and the scale $t$} if
\begin{align}
\mathbf{E}(p)= s \quad\quad and \quad \quad
\min\big\{d(p, q^+), d(p, q^-)\big\}= t \nonumber
\end{align}
where $\mathbf{E}(\cdot)= d(\cdot, q^+)+ d(\cdot, q^-)- d(q^+, q^-)$.
}
\end{definition}

We recall the following Abresch-Gromoll lemma (\cite{AG}, also see \cite[Lemma $2.2$]{Xu-group}).
\begin{lemma}\label{lem general Abresch-Gromoll}
{On complete Riemannian manifold $(M^n, g)$ with $Rc\geq 0$, assume that $[q^+, q^-, p]$ is an AG-triple with the excess $\leq \frac{1}{n}\frac{r^2}{R}$ and the scale $\geq R$, furthermore assume $R\geq 2^{2n}r$, then $\displaystyle \sup_{B_{r}(p)} \mathbf{E}\leq 2^6\cdot \big(\frac{r}{R}\big)^{\frac{1}{n- 1}}r$.
}
\end{lemma}\qed

On Riemannian manifolds, if there is a segment $\gamma_{p, q}$ between two points $p, q$, we can choose the middle point of the segment $\gamma_{p, q}$, denoted as $z$. Then $[p, q, z]$ is an AG-triple with the excess $0$ and the scale $\frac{1}{2}d(p, q)$. The following result is the almost version of the above fact.

\begin{prop}\label{prop diameter pair points}
{For $\tau\in (0,\frac{1}{6(n+1)}), 3\leq \beta\leq \tau^{-c(n)}$ and $r>0$, if $B_{r}(p)\subset (M,g)$ satisfies $Rc\ge 0$ and $\displaystyle \frac{V(B_{r}(p))}{V(B_r(0))}\ge 1-\tau$. Then for $q\in B_{\beta^{-C(n)} r}(p)- B_{\frac{1}{2}d(p,q_n)}(p)$, there is $q^{-}$  such that $[q,q^-,p]$ is an AG-triple with the excess $\le C(n)\beta^{-C(n)}d_q(p)$ and the scale $\ge d_q(p)$.
}
\end{prop}

\pf
{By (the proof of) Lemma \ref{lem dist between points and segment}, we know
\begin{align*}
\frac{V(A_{d_{q}(p)-\epsilon r_{q},d_{q}(p)+\epsilon r_{q}}(q)\backslash \mathcal{C})}{V(B_{\epsilon r_{q}}(p))}\le C(n)\frac{\beta^{-C(n)}}{\epsilon^{n-1}},
\end{align*}
where
$$\mathcal{C}=\{y\in A_{d_q(p)-\epsilon r_q, d_q(p)+\epsilon r_q}(q)|  \exists \theta(y)\in S(T_{q}M) s.t.  y= \exp_q(\rho(y)\theta(y)), l_{\theta(y)}\ge 3 d_q(p)\}.$$
Taking $\epsilon=(2 C(n)\beta^{-C(n)})^{\frac{1}{n-1}}$, we know  $C(n)\frac{\beta^{-C(n)}}{\epsilon^{n-1}}=\frac{1}{2}$, which means there exists $\tilde{p}\in B_{(2 C(n)\beta^{-C(n)})^{\frac{1}{n-1}}r_{q}}(p)$ such that $\tilde{p}=\exp_{q}(\rho(\tilde{p})\theta(\tilde{p}))$ and $l_{\theta(\tilde{p})}\ge 3d_q(p)$. That is,
$\exp_{q}(t\theta(\tilde{p})):[0,3d_q(p)]\to M$ is a segment.  So,
$$d(\exp_{q}(3d_q(p)\theta(\tilde{p})),p)\ge d(\exp_{q}(3d_q(p)\theta(\tilde{p})),q)-d(p,q)=2d_q(p)$$
and
$$d(\exp_{q}(\rho(\tilde{p})\theta(\tilde{p})),p)=d(\tilde{p},p)\le C(n)\beta^{-C(n)}r_q\ll d_q(p).$$
By continuity, there exists $t_1\in [\tilde{\rho},3d_q(p)]$ such that $q^-=\exp_{q}(t_1\theta(\tilde{p}))$ satisfies $\displaystyle d(q^-,p)=d_q(p)$ and
$$d(q^{-},q)=d(q^{-},\tilde{p})+d(\tilde{p},q)\ge d(q^-,p)+d(p,q^-)-2d(p,\tilde{p})\ge (2-C(n)\beta^{-C(n)})d_q(p).$$
Moreover, we have
$$\mathbf{E}(p)=d(p,q)+d(p,q^{-})-d(q,q^-)\le 2d(p,\tilde{p})\le C(n)\beta^{-C(n)}r_q.$$
}
\qed

Recall we have the following existence result of splitting function $\mathbf{b}$ with respect to the local Busemann function $b^+$.

\begin{lemma}\label{lem existence of harmonic function-r}
{On complete Riemannian manifold $M^n$ with $Rc\geq 0$, assume that $[q^+, q^-, p]$ is an AG-triple with the excess $\leq \frac{4}{n}\frac{r^2}{R}$ and the scale $\geq R$, also assume $R\geq 2^{2n+ 1}r$. Then there exists harmonic function $\mathbf{b}$ defined on $B_{2r}(p)$ such that
\begin{align}
&\sup_{B_r(p)}|\mathbf{b}- b^+|\leq C(n)(\frac{r}{R})^{\frac{1}{n- 1}}r, \nonumber \\
&\sup_{B_{r}(p)} |\nabla \mathbf{b}|\leq 1+ 2^{51n^2}\big(\frac{r}{R}\big)^{\frac{1}{4(n- 1)}} \quad \quad and \quad \quad
\fint_{B_{r}(p)} \big|\nabla (\mathbf{b}- b^+)\big|^2 \leq C(n)\big(\frac{r}{R}\big)^{\frac{1}{n- 1}} , \nonumber
\end{align}
where $b^{+}(x)= d(x, q^{+})- d(p, q^+)$.
}
\end{lemma}

\pf
{see \cite[Lemma $2.6$]{Xu-group}, especially \cite[(2.12), (2.15)]{Xu-group} there.
}
\qed

Now we establish our main theorem about the construction of $\epsilon$-splitting map.
\begin{theorem}\label{thm AG-triple imply one more splitting-pre}
{If $Rc(M^n)\geq 0$ and $V\big(B_r(p)\big)\geq (1- \tau)V\big(B_r(0)\big)$, then for $\beta=\tau^{-c(n)}$ and $r_1=c(n)\beta^{-(C(n)+\frac{\nu^{n+1}-\nu}{\nu-1})}r$, there are harmonic functions $\big\{\mathbf{b}_i\big\}_{i= 1}^{n}$ defined on $B_{r_1}(p)\subset B_{r}(p)$, such that
\begin{align}
\sup_{B_{r_1}(p)\atop i= 1, \cdots, n} |\nabla \mathbf{b}_i|\leq 1+ C(n)\beta^{-1} \quad and \quad \fint_{B_{r_1}(p)} \big|\langle \nabla \mathbf{b}_i, \nabla \mathbf{b}_j\rangle- \delta_{ij}\big|^2\leq C(n)\tau^{c(n)}.\nonumber
\end{align}
}
\end{theorem}

\pf
{We can apply Proposition \ref{prop approximation component is almost orthogonal} to get
\begin{align}
\fint_{B_{r_1}(p)}\sum_{i, j= 1}^{n} \big|\langle \nabla b_i^+, \nabla b_j^+\rangle- \delta_{ij}\big|\leq C(n)(\beta^{-\frac{1}{2}}+\tau^{\frac{2}{n+1}}).\label{almost o.n.-need-1}
\end{align}

From Proposition \ref{prop diameter pair points},  it is easy to see that there exist  AG-triple $[q_i,q_i^{-},p]$ with the excess $\leq C(n)\beta^{-C(n)}d_{q_i}(p)$ and the scale $\geq d_{q_i}(p)$ for $i= 1, \cdots, n$.

Note that
$$
\max\{C(n)\frac{r_1}{d_{q_i}(p)},\frac{C(n)\beta^{-C(n)}d_{q_i}(p)}{r_1}\}\le C(n)\beta^{-\nu}+C(n)\beta^{-(C(n)-\frac{\nu^{n+1}}{\nu-1})}\le C(n)\beta^{-\nu}.$$
 Now we can apply Lemma \ref{lem existence of harmonic function-r} to obtain harmonic functions $\big\{\mathbf{b}_i\big\}_{i= 1}^{n}$ satisfying
\begin{align}
&\sup_{B_{r_1}(p)\atop i= 1, \cdots, n} |\nabla \mathbf{b}_i|\leq 1+ 2^{51n^2}\Big(C(n)\beta^{-\nu}\Big)^{\frac{1}{4(n- 1)}}\leq 1+ C(n)\beta^{-\frac{\nu}{4(n-1)}}  \label{gradient bound of b need} \\
&\sup_{i= 1, \cdots, n}\fint_{B_{r_1}(p)} \big|\nabla (\mathbf{b}_i- b_i^+)\big|^2 \leq 2^{4n}\Big(C(n)\beta^{-\nu}\Big)^{\frac{1}{n- 1}}\leq  C(n)\beta^{-\frac{\nu}{n-1}}.\label{harmonic is almost linear}
\end{align}

From (\ref{gradient bound of b need}) and (\ref{harmonic is almost linear}), we get
\begin{align}
&\quad \fint_{B_{r_1}(p)} \big|\langle \nabla \mathbf{b}_i, \nabla\mathbf{b}_j\rangle- \delta_{ij}\big|^2 \nonumber \\
&\leq 3\fint_{B_{r_1}(p)} \big|\nabla (\mathbf{b}_i- b_i^+)\big|^2\cdot |\nabla \mathbf{b}_j|^2+
 \Big|\big\langle \nabla b_i^+, \nabla(\mathbf{b}_j- b_j^+)\big\rangle\Big|^2+ \big|\langle \nabla b_i^+, \nabla b_j^+ \rangle- \delta_{ij}\big|^2 \nonumber \\
 &\leq C(n)\beta^{-\frac{\nu}{n-1}}+ 2\fint_{B_{r_1}(p)} \big|\langle \nabla b_i^+, \nabla b_j^+ \rangle- \delta_{ij}\big| .\nonumber
\end{align}

From (\ref{almost o.n.-need-1}) and the above inequality, we have
\begin{align}
\fint_{B_{r_1}(p)} \sum_{i, j= 1}^{n}\big|\langle \nabla \mathbf{b}_i, \nabla\mathbf{b}_j\rangle- \delta_{ij}\big|^2 \leq \epsilon_1=C(n)(\beta^{-\frac{1}{2}}+\beta^{-\frac{\nu}{n-1}}+\tau^{\frac{2}{n+1}}).\nonumber
\end{align}
Note that we have fixed $\nu>4(n-1)$, we get the conclusion.
}
\qed

\section*{Acknowledgments}
We thank Zuoqin Wang for arranging our visit to University of Science and Technology of China, part of the work was done during the visit. Also we thank Zichang Liu for sharing his preprint \cite{Liu} and some comments on the earlier version of this paper.

\end{document}